\documentclass[11pt]{article}
\usepackage[utf8]{inputenc}
\usepackage{setspace}
\usepackage[margin=1.0in]{geometry}
\usepackage{amssymb, amsmath}
\usepackage{amsthm}
\usepackage{colordvi,verbatim,hyperref}
\usepackage{color,enumerate}
\usepackage{graphicx}
\usepackage{caption}
\usepackage{algorithm,algorithmic}
\usepackage{subcaption}
\usepackage[justification=centering]{caption}
\usepackage{authblk}
\usepackage[toc,page]{appendix}
\usepackage{lineno}
\usepackage{rotating}
\usepackage{tablefootnote}
\usepackage{mathtools}
\usepackage{natbib}

\theoremstyle{definition}

\newtheorem{cor}{Corollary}

\newtheorem{lm}{Lemma}
\newtheorem{thrm}{Theorem}

\allowdisplaybreaks

\title{\textbf{When Should the FDA Inspect Pharmaceutical Manufacturing Facilities to Better Mitigate Drug Shortages?}}
\author[1]{Daniel Kosmas}
\author[1]{{\"O}zlem Ergun}
\affil[1]{\footnotesize Department of Mechanical and Industrial Engineering, Northeastern University, Boston, MA 02115, USA}
\date{}

\begin{document}

\maketitle

\doublespacing
\vspace{-1cm}
\begin{abstract}
    \noindent\textbf{Problem Definition:} Drug shortages have been a persistent problem in American healthcare for decades, and the resulting lack of access to necessary drugs has been disastrous to patient health. A majority of these shortages were caused by quality issues related to problems in the manufacturing process. More frequent inspections can help reduce quality concerns, but deciding when to inspect is a complex problem; strict regulation enforcement can force low-profit facilities to close due to excessive maintenance costs, while lax enforcement allows for regulation violations to persist, both of which can cause drug shortages.
    
    \noindent\textbf{Methodology/Results:} We propose a novel model to assist the FDA in determining when to inspect manufacturing facilities. We formulate this problem as a finite-horizon partially observable Markov decision process (POMDP) based on the classifications the FDA assigns to each facility after inspection, as well two  disruptive events: a manufacturing failure occurring or the facility closing for non-mandatory maintenance. We theoretically show that this problem can be reduced to only needing to consider whether or not to inspect immediately, which is independent of the time horizon. We additionally determine the sensitivity of the optimal inspection time on the penalty incurred for an unexpected disruptive event occurring. 
    
    \noindent\textbf{Managerial Implications:} 
    Our computational study demonstrates a quadratic relationship between the relative difference in average value accumulated between inspecting based on the optimal inspection time produced by our model and inspecting based on the expected time to an unexpected disruptive event, highlighting the importance of allocating more inspection resources to high-risk facilities that produce drugs that highly impact public health. We additionally find that optimal inspection time is more sensitive to changes in the penalty incurred from a disruptive event occurring the longer it has been since the last inspection.
\end{abstract}

\section{Introduction}
\label{sec:intro}
% paragraph on drug shortages
Drug shortages have plagued American health care for decades, with no signs of improvement in sight \citep{tucker2020drug}. Over 680 active drug shortages were reported between 2011 and 2020 \citep{fda2020report}. The number of active drug shortages increased by 30\% between 2021 and 2022, reaching a 5-year record high of 295 active drug shortages at the end of 2022. The average drug shortage lasts 1.5 years, and at least 15 critical drugs have been in shortage for over a decade \citep{hsgac2023shortsuppy}. Drug shortages are disastrous for patients and health care providers, as a lack of access to necessary drugs can result in medication errors and treatment delays or cancellations that can lead to adverse patient health outcomes or even death \citep{hsgac2023shortsuppy, ventola2011drug}. In 2019, the Food and Drug Administration (FDA) found that quality issues were the leading cause of drug shortages, accounting for 62\% of shortages \citep{fda2019shortages}. These quality issues are often a result of manufacturing problems in pharmaceutical manufacturing facilities \citep{francas2023drivers, jensen2002fda}.

Pharmaceutical manufacturing facilities can experience numerous problems that lead to unexpected disruptions in the manufacturing process, such as quality issues with ingredients, outdated manufacturing equipment and natural disasters damaging the facilities. These disruptions are costly for both the manufacturer, who is responsible for resolving these regulatory issues in order to resume production, and the consumers, who receive a lower quality product or no product at all \citep{fda2019shortages}. The true cost and length of an 
unexpected disruption are often incredibly difficult to predict, making it difficult for the FDA to determine if they need to intervene \citep{nasem2022building}. In the worst case scenario, the capacity reduction caused by the unexpected disruption can result in a drug shortage \citep{fda2019shortages}. Since facilities are only required to report when production is interrupted or permanently discontinued, inspections are essential to ensuring that risk factors that could lead to disruptions are appropriately mitigated \citep{Ergun2022SupplyChains,fda2020notifying,nasem2022building}. Through routine inspections, the FDA can help prevent unexpected disruptions from occurring, whether it be from preventing low quality products entering the market that may later need to be recalled or from identifying issues in the production line that may result in the facility needing to shut down for maintenance \citep{anand2012decay}. 

If the FDA inspects a facility and finds that production conditions are not in line with the current regulations, they are equipped with various regulatory tools to encourage companies to address these concerns \citep{woodcock2019testimony}. However, the decision to ramp down production when issues are identified is not without drawbacks. If the FDA is too penalizing with their decisions or the manufacturing process issues are too costly for the manufacturer to address, the facility may be forced to stop production for maintenance or decide to discontinue production of the product altogether, which is another disruptive event that reduces production capacity. This may unintentionally result in a shortage due to a number of factors, such as the combination of a sudden ramp down in production and limited FDA resources delaying re-inspection or other manufacturers of similar products being unable to ramp up their production to account for the facility's closure \citep{kimball2023temporary, ventola2011drug, woodcock2020quality}. However, delayed re-inspections are also problematic when managing unexpected disruptions, so it is better for the FDA to be aware of the severity of the facility's problems. If the FDA is aware of the severity of the facility's problems, they can take action to mitigate the impact of the disruption they cause, such as working with the facility to address their concerns without fully halting production or expediting the process of enabling new manufacturers to enter the market \citep{jensen2002fda, kimball2023temporary, ventola2011drug}. Additionally, it is better for the FDA to force a facility to halt production if the defective product can negatively impact patient health \citep{buckley2013countering,pew2012,white2020generic}. A product negatively impacting patient heath is especially a concern when it has few (or no) substitute-able alternatives since the lack of alternatives forces patients to use a product that may be more harmful than not receiving treatment at all \citep{edney2019cheap, wang2023barriers}. Henceforth, we distinguish between a disruption caused by the FDA inspection resulting in the facility needing to close for mandatory maintenance and a disruption caused by an event that the FDA is not aware of, such as manufacturing failures and non-mandatory maintenance, which we refer to as an \emph{unexpected disruptive event}. 

If the FDA does not inspect frequently enough, manufacturers have less incentive to maintain high quality standards \citep{wu2020monitoring}. Declining quality can lead to products that need to be recalled, which, if caught, can cause shortages if there are few (or no) alternative products. If not caught, the product can negatively impact patient health. Thus, the FDA must strike a careful balance of inspecting frequently enough to incentivize pharmaceutical manufacturing facilities to maintain high production quality and identify production issues before they can cause a shortage, while not inspecting too frequently, wasting their resources and potentially causing shortages through enforcing regulatory requirements so strictly that the manufacturer deems production to be unprofitable. The FDA must additionally balance how inspection resources are allocated across all pharmaceutical manufacturing facilities, forcing the FDA to identify which facilities are more critical for public health, and thus need more resources allocated to monitoring them.

We propose a finite horizon partially observable Markov decision process (POMDP) \citep{krishnamurthy2016pomdp} to help support the FDA in determining when to inspect pharmaceutical manufacturing facilities that balances the reward of allowing a facility to continue production uninhibited with the risks of unexpected disruptive events. We analytically show that the optimal inspection time suggested by our model can be deduced by simple computations each time period. We additionally analyze how the sensitivity in the severity of the impact caused by an unexpected disruptive events will change the optimal inspection time. 

The remainder of this paper is organized as follows. Section \ref{sec:lit} reviews relevant literature around applying tools from operations research to decision-making for inspections. Section \ref{sec:model} describes the model of determining how much time should pass between sequential inspections of manufacturing facilities, and derives closed-form equations for the value functions of the POMDP. Section \ref{sec:vf} demonstrates theoretical properties of this model that allow us to determine the optimal inspection time only by determining if we should inspect now or wait one time period, and then analyzes the sensitivity of the optimal inspection time suggested by the model on the penalties incurred by the facility experiencing an unexpected disruptive event. Section \ref{sec:modelvar} proposes a variant to this model to account for the outcome of the inspection, additionally resulting in a penalty for the inspector if the facility needs to close for mandatory maintenance. Section \ref{sec:numres} computationally explores the trade-offs of different inspection rules for varying transition matrices and shows numerically the sensitivity of the optimal inspection time and total value accumulated on various parameters. Section \ref{sec:con} concludes the paper and proposes avenues for future research.

\section{Literature Review}
\label{sec:lit}
%In this section, we first review literature related to FDA inspections. We then review literature on how POMDPs have been applied to maintenance and inspection problems. We also discuss mission abort problems, which have also been modeled using POMDPs and have similar features to the model we propose.

% part 1 pharma literature
The FDA inspects pharmaceutical manufacturing facilities to verify compliance with applicable laws and regulations, known as the current good manufacturing practices (CGMP) \citep{fda2022cgmp}. When the FDA inspects a facility, they assign one of three classifications: No Action Indicated, Voluntary Action Indicated, or Official Action Indicated. No Action Indicated means that no violations of regulations were identified during the inspection. Voluntary Action Indicated means that the inspector has identified violations, but the FDA chooses to not take regulatory action. Official Action Indicated means that the violations identified are severe enough for the FDA to take regulatory action \citep{fda2020faq}. Under the Federal Food, Drug and Cosmetic Act, the FDA was required to inspect pharmaceutical manufacturing facilities at least once every two years. However, the FDA struggled to meet these requirements for all facilities due to limited resources \citep{woodcock2019testimony}. In 2005, a risk-based site-selection model (SSM) was developed by the Center for Drug Evaluation and Research (CDER) to help support the FDA's decision-making around determining which facilities to inspect \citep{nasem2022building, woodcock2019testimony}. The SSM incorporates factors such as facility type, time since last inspection, compliance history and inherent product risk to determine a ``risk score" for each facility, which helps the FDA determine which facilities are more in need of being inspected. While the specifics of how the scores are determined is not publicly known, it is known that these risk scores are based on factors intrinsic to the facility and product(s) being produced and do not account for how a disruption would impact public health \citep{cder2018manual}. The two-year minimum requirement was replaced in the Food and Drug Administration Safety and Innovation Act (FDASIA) of 2012. Under FDASIA, the inspection requirement was expanded to include both domestic and foreign pharmaceutical manufacturing facilities, but the fixed minimum inspection interval was fully removed, leaving inspection scheduling solely dictated by the SSM \citep{woodcock2019testimony}. Our analytical results allow us to directly supplement the SSM since the optimal inspection time suggested by our model can be deduced by simple computations each time period. These computations can be converted into a score which can be integrated into the SSM.

There has been some quantitative work to support the FDA in determining when to inspect different pharmaceutical manufacturing facilities. \cite{klimberg1992improving} developed an integer program where facilities are either inspected in year one or year two, and the set of facilities is partitioned into these two sets. \cite{macher2006exploring} investigated how different factors, such as days between inspection, compliance history, and the quality of training of the inspectors, impact the classification a facility receives. Similarly, \cite{pazhayattil2022quantitative} explored how violations of different sections of the CGMP were correlated with the inspection classifications. \cite{hein2019comparing} compared various data-driven methods to predict the classification facilities would receive. \cite{ahuja2021enhancing} used data-driven methods to determine the safety of pharmaceutical products, helping determine if FDA intervention for specific products is necessary. \cite{anand2012decay} demonstrated how FDA inspections act as an inhibitor to operational decay. They were also proponents of implementing targeting high-risk facilities for inspection. \cite{ball2017investigator} showed that inspection frequency is predictive of future facility quality and highlighted that inspector complacency can be detrimental ensuring that inspections enforce quality standards. \cite{wu2020monitoring} explore the impact of FDA surveillance inspections on domestic pharmaceutical manufacturing facilities producing generic drugs. They suggest that more frequent monitoring of facilities improved product quality, subject to diminishing returns. Their analysis found that more frequent inspections were most beneficial in reducing recalls for products from facilities that were previously inspected less than once every two years. They additionally note that cost concerns can cause high-risk facilities to close, rather than perform the necessary maintenance to comply with the CGMP, suggesting that the FDA must consider how manufacturers will react to inspections. \cite{zhang2019unannounced} proposed a game theoretic model to determine how strictly a supervisor should monitor a pharmaceutical manufacturing facility to prevent drug fraud. While pharmaceutical inspection has received some academic attention, there has been more work on supporting the FDA in inspecting food manufacturing facilities which have focused more on data-driven models for predicting inspection results and determining how increased inspection improved self-regulation amongst food vendors \citep{babich2012managing,chang2017risk,jin2021food,levi2019supply,zhou2022effects}. Our work differs from past work on FDA inspection by analytically determining when a facility should be inspected, as opposed to using data-driven methods or determining a rate of inspection frequency. Our model may also be applied to the inspection of food manufacturing facilities as well, as FDA uses the same classification system to determine their compliance with relevant laws and regulations.

While quantitative models supporting FDA inspection decision making has primarily focused on game theoretic models or data driven models, determining the optimal inspection and maintenance time for general maintenance has also been well studied through the lens of condition-based maintenance models \citep{ahmad2012condition, zhu2021cbm}. Condition-based maintenance has been popularized for its ability to react to the observed conditions of the system, preventing system failure when degradation is observed sooner than expected. Markov decision processes (MDPs) have commonly been leveraged for condition-based maintenance due to their ability to characterize multiple maintenance actions and multiple degradation states \citep{papakonstantinou2014inspection,papakonstantinou2014implementation}. In particular, POMDPs have been utilized in the optimal maintenance of numerous systems (eg., \citep{byon2010wind,ge2007optimum,ghandali2020sustainable, papakonstantinou2014pomdp,zhang2017continuous,Zouch2011road}). More recently, \cite{morato2022dynamic} proposed integrating Bayesian networks with POMDPs to provide a joint framework for inspection and maintenance planning. \cite{xu2022constrained} introduces risk-averse variants to the traditional MDP model for maintenance, and proves that there exists an optimal maintenance policy when optimizing Value-at-Risk (VaR) and Conditional Value-at-Risk (CVaR). \cite{liu2023sequential} proposed a framework for the sequential inspection and maintenance of multi-state systems, where a fixed number of inspection or maintenance actions can be taken each break and the actions are implemented sequentially, allowing for the decisions to be updated as inspections occur. Our problem differs from traditional inspection and maintenance problems because the inspector does not have the capabilities to perform the maintenance needed. As such, the inspector can only plan when they will inspect and try to account for whether the outcome of the inspection will determine that maintenance is needed. Additionally, most work on applying POMDPs to maintenance and inspection focuses on numerical experiments. We additionally provide a theoretical analysis of our problem to determine optimal inspection time without the need of more complex solvers.

% part 4 mission abort literature
Our problem is also related to the single-component mission abort problem \citep{levitin2019cost}. The goal of mission abort problems is to determine if a degrading system is able to complete its mission before failure occurs, and if not, when to terminate the mission. The operator only receives a reward if the mission is successfully completed. If the operator chooses to abort the mission, the penalty of failure is less than that of if the system fails. In this problem, the system is observed at regular time intervals. The state of the system is healthy, unhealthy, or failed. However, when inspected, the only information about the state that is revealed is whether or not it is operational, leaving ambiguity between whether or not the system is healthy or unhealthy. This is congruent with current regulations for pharmaceutical manufacturing facilities; no information regarding the quality of these facilities is shared with FDA unless a manufacturing failure has occurred \citep{nasem2022building}. Typically, if some issue has come to the attention of the FDA, they will conduct a ``for-cause" inspection to immediately investigate the severity of the issue \citep{fda2018inspect}. \cite{myers2009probability} was the first to study the mission abort problem, which considers a multi-component system. Later, \cite{levitin2019cost} explored the single-component mission abort problem, where they model the mission abort as the initiation of the rescue operation to prevent the destruction of the system. \cite{qiu2019gamma} derives degradation-based and duration-based abort policies for systems that degrade according to a Gamma process. \cite{qiu2023optimal} formulates the mission abort problem as a stochastic dynamic program and compares the effectiveness of mission abort decisions against various heuristic measures, such as the length of time spent in the unhealthy state. They numerically show that heuristic abort strategies were less economically efficient then strategies provided by solving the stochastic dynamic program.

To the best of our knowledge, very few works formulate the mission abort problem as a MDP. \cite{zhao2021optimal} determines the optimal mission abort time for a system that degrades according to a stationary Gamma process. This degradation process is then converted to a continuous-state MDP. \cite{cheng2023optimal} utilizes a hidden Markov model to model the system's degradation, which is then transformed into a continuous-state MDP. Both of these works solve their models via value iteration \citep{krishnamurthy2016pomdp}, whereas we derive a closed-form equation for our value functions. Moreover, our results allow us to perform simple calculations to determine the optimal inspection time. Our problem differs from existing literature on MDP formulations of the mission abort problems in the following ways. First, these models only receive a reward if the mission is fully completed, whereas we accumulate value each time period the manufacturing facility is operating. Second, while some works generalize the number of quality states the system has, previous works on mission abort problems only have theoretical results for two operational states and one absorbing state. The model we consider has three operational states and two absorbing states, and we provide theoretical results on the optimal inspection time for our model.

\section{POMDP Model}
\label{sec:model}
We consider a finite time POMDP model to determine when the FDA should inspect a pharmaceutical manufacturing facility. We focus on determining when to inspect a single facility independently of other facilities since the SSM generates scores for each individual facility based on risk factors intrinsic to that facility and does not consider them interdependently \citep{cder2018manual}. In our model, the FDA must inspect a facility by a pre-specified time. While in reality the formal requirement to inspect in by a certain time has been removed, it is reasonable to assume that the FDA would not allow a facility to go uninspected after a certain amount of time has passed. We consider two unexpected disruptive events in our model: a manufacturing failure occurring and the facility closing for non-mandatory maintenance. A manufacturing failure is an unexpected event that results in a reduction or elimination of production capacity due to issues in the manufacturing process, while a facility closing to perform non-mandatory maintenance is a proactive event by the facility to shut down production for a set timeline, but the details of this event are unknown to the FDA. While we only consider two unexpected disruptive events, this model can easily be generalized to multiple unexpected disruptive events.

In each time period, the facility may be in operation, may have closed due to a manufacturing failure, or may have closed to perform non-mandatory maintenance. The FDA receives a reward if the facility is in operation.  Otherwise, the FDA incurs a large penalty if the facility has experienced a manufacturing failure, or a small penalty if the facility closes for non-mandatory maintenance. These penalties reflect the impact of facility no longer manufacturing their pharmaceutical product(s) on social well-being. While a facility undergoing maintenance is positive in long-term planning, the immediate halt in production without prior knowledge by the FDA is disruptive and can potentially cause a shortage. We assume that closing for non-mandatory maintenance has a smaller penalty than closing due to a manufacturing failure because the true impact and timeline of a manufacturing failure is difficult to predict, while non-mandatory maintenance will have a more consistent timeline. 

The FDA accrues rewards until they decide to inspect (or must inspect if the end of the time horizon is reached) or one of the  disruptive events occurs. The goal of the FDA is to maximize the total reward accrued. We first consider a model where there is no penalty for inspecting, regardless of quality of the facility. We also consider a model where inspecting can also result in the facility closing for mandatory maintenance. As with closing for non-mandatory maintenance, closing for mandatory maintenance as a result of an FDA inspection can potentially cause a shortage. However, in this case the FDA will be aware that the closure will be occurring, allowing them to better mitigate the impact of the closure. As opposed to defining the states of the POMDP to be arbitrary quality levels of the facility, then determining what classification the FDA should assign to each quality level, we instead define the states to be the classifications that the FDA will give after inspecting a facility with a certain quality level. Since the FDA cannot directly control the quality levels of the facility, their primary concern is identifying the classification the facility should receive. This assumes that facility quality levels are able to be coherently clustered into the classification system, and this clustering is not influenced by external factors.

In each time period the facility is in operations, the FDA must decide whether or not they will inspect the facility at the end of the time period. In order to determine when is best to inspect, the FDA maintains their private belief about how likely the facility is to receive each classification rating in each time period. If they do not inspect, there is a chance that an unexpected disruptive event will occur in the next time period. When they inspect, the true classification of the facility is revealed and no future reward is received. We assume that, when the FDA decides to inspect, they observe the quality of the facility accurately and will assign the appropriate classification to a facility based on the conditions it is in, regardless of the implications on public health. Outside of a manufacturing failure occurring or inspecting, the FDA receives no information on the status of the facility  \citep{nasem2022building}. 

We now mathematically define our model; all notation is summarized in Appendix \ref{app:notation}. Given a finite horizon $T$, the FDA must choose when to inspect the facility. They may do so any time before $T$, but they are required to inspect the facility by $T$. Once an inspection occurs, the time requirement on inspection resets. Thus, the FDA choosing to inspect the facility acts as the end of the game. Additionally, if an unexpected disruptive event occurs, the FDA receives the appropriate penalty for not inspecting before the event occurred, and then must inspect upon the facility reopening. Since we are only considering the length of time between sequential inspections for a single facility, we only need to consider whether an unexpected disruptive event has occurred or not. Additionally, since the time requirement on inspection resets after each inspection, we only need to model the FDA's decision making when they should next inspect, and do not need to determine the long-term decision making over an infinite time horizon. When the FDA inspects the facility, they would gain the information necessary to update model parameters as appropriate before restarting the model.

The classifications the FDA assigns serve as the basis for the states of the Markov chain: No Action Indicated is $N$, Voluntary Action Indicated is $V$, and Official Action Indicated is $O$. We additionally include states indicating that the facility has experienced a manufacturing failure ($D$), has closed for non-mandatory maintenance ($C$) or has been inspected ($I$). Each period, the FDA can receive one of three observations regarding facility quality. If the facility experiences an unexpected disruptive event, it must be reported to the FDA \citep{fda2020notifying}. If the facility has experienced a manufacturing failure, the FDA observes that a manufacturing failure has occurred $(dr)$, and if the facility has closed for non-mandatory maintenance, the FDA observes that the facility is closing for non-mandatory maintenance $(cr)$, indicating that they know if the facility is in states $D$ or $C$ with full certainty. Otherwise, the FDA observes no reported issues from the facility, $(nr)$, indicating that the facility is still operating. The FDA receives a reward $r$ each period that the facility is successfully producing. If the facility experiences a manufacturing failure, then the FDA receives a penalty of $-d$. If the facility closes for non-mandatory maintenance, then the FDA receives a penalty of $-c$. This penalty is to reflect that the facility is no longer producing, which may cause a shortage. After the facility has been inspected, the FDA receives a reward of $0$, indicating that the facility has been inspected before an unexpected disruptive event occurred and the game is over. The focus of this initial model is on inspecting before an unexpected disruptive event occurs, and does not account for the outcome of the inspection. This model represents the ideal situation for the FDA: where inspection focuses on identifying whether or not a facility is in compliance with the CGMP, and does not need to consider how the outcome of the inspection can negatively impact public health. We later extend our model to account for the outcome of the inspection in Section \ref{sec:modelvar}. For simplicity, we set $r=1$ and re-scale other parameters to reflect this.

\begin{figure}
    \centering
    \includegraphics[width=.45\linewidth]{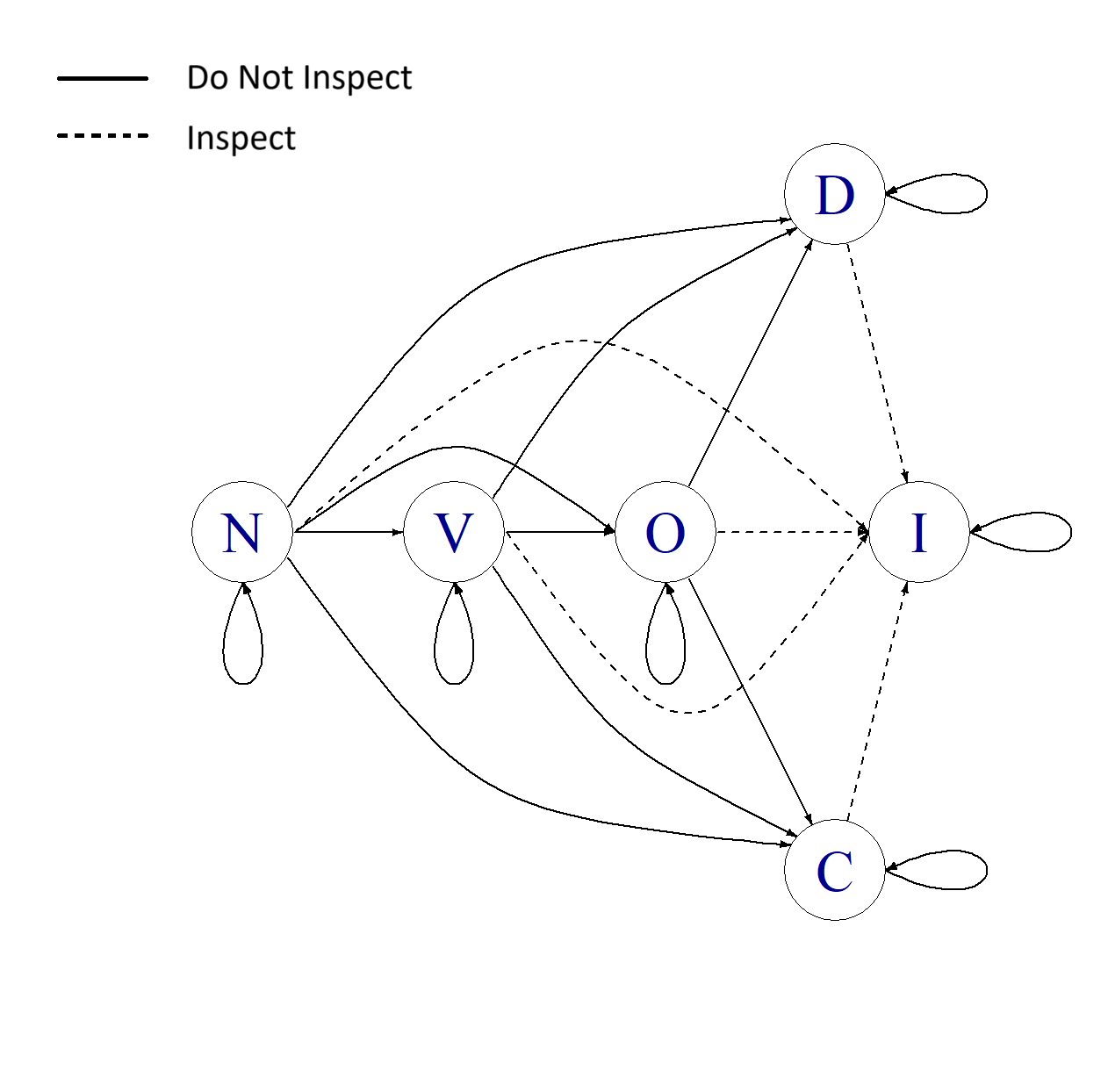}
    \caption{State space of the Markov chain}
    \label{fig:states}
\end{figure}

At each state, the FDA can either choose to inspect $(i)$ or not inspect $(ni)$. We assume the cost of inspection has already been budgeted for, and thus do not include it in our model. Figure \ref{fig:states} presents the state space of the Markov chain, as well as the transition between states. Solid arcs represent the transitions that can occur when the FDA does not inspect, while dashed arcs represent transitions that occur when the FDA does inspect.
If the FDA does not inspect, the quality of the facility can degrade, with the following transition probabilities:

\small
\begin{equation*}
    P(ni) = 
\left[
\begin{array}{c|cccccc}
 & N & V & O & D & C & I \\ \hline
N& p_{NN} & p_{NV} & p_{NO} & p_{ND} & p_{NC} & 0 \\
V& 0 & p_{VV} & p_{VO} & p_{VD} & p_{VC} & 0 \\
O& 0 & 0 & p_{OO} & p_{OD} & p_{OC} & 0 \\
D& 0 & 0 & 0 & 1 & 0 & 0 \\
C& 0 & 0 & 0 & 0 & 1 & 0 \\
I& 0 & 0 & 0 & 0 & 0 & 1 
\end{array}
\right]
\end{equation*}
\normalsize

Once the FDA has inspected, the game is over and no more reward can be collected. Thus, inspecting at any state enforces that the next state will be $I$, which is an absorbing state, i.e., $P(i)_{SI} = 1$ for all $S$ and $P(i)_{SS'} = 0$ for $S' \ne I$.

%\small
%\begin{equation*}
%    P(i) = 
%\left[
%\begin{array}{cccccc}
%0 & 0 & 0 & 0 & 0 & 1 \\
%0 & 0 & 0 & 0 & 0 & 1 \\
%0 & 0 & 0 & 0 & 0 & 1 \\
%0 & 0 & 0 & 0 & 0 & 1 \\
%0 & 0 & 0 & 0 & 0 & 1 \\
%0 & 0 & 0 & 0 & 0 & 1 
%\end{array}
%\right]
%\end{equation*} \normalsize

In each time period, the FDA maintains their private belief about how likely it is for the facility to be in each state. Let $b^t = (b^t_N, b^t_V, b^t_O, b^t_D, b^t_C, b^t_I)$ be the vector indicating the belief probabilities of being in each state at time period $t$. Note that, based on our assumptions, we can classify these belief vectors into four categories. Once the facility has been inspected, we have that $b^t = (0,0,0,0,0,0,1)$. Likewise, if the facility has experienced a manufacturing failure, we have that $b^t = (0,0,0,1,0,0)$, and if the facility has closed for non-mandatory maintenance, $b^t = (0,0,0,0,1,0)$. If no unexpected disruptive event has been reported and we have not inspected, then the facility is still manufacturing their product, and thus $b^t = (b^t_N, b^t_V, b^t_O, 0, 0, 0)$.

We model the FDA determining when they should choose to inspect a facility as them inspecting in the time period that would maximize the expected rewards accrued over the entire time horizon. As opposed to the traditional stochastic-dynamic formulation of the value functions, we focus on the conditional plan formulation. This choice allows us to compute the expected values of sub-optimal decisions, as opposed to just determining that those decisions are sub-optimal. We can reasonably generate all conditional plans since choosing to inspect and two of the three observations we can receive after deciding to not inspect result in the game ending (i.e., no future decision needs to be made in response to those observations). Thus, all conditional plans can be described by when the FDA chooses to inspect the facility or is forced to inspect after an unexpected disruptive event has occurred. Let $\sigma_{j}^t$ be the conditional plan starting at time $t$ where we inspect after $j$ periods of the facility being in operation, where $0 \le j \le T-t$. To be clear, $\sigma_0^{t}$ inspects immediately and $\sigma_{T-t}^t$ does not inspect the facility until the end of the time horizon. Similarly, we can describe a conditional plan as the sequence of actions taken in each time period, i.e., $(a_1, \ldots, a_T)$, where $a_t$ is the action taken at time $t$. Since no actions taken after an inspection occurs, we can reduce this notation to $(a_1, \ldots, a_{k-1}, a_k) = (ni, \ldots, ni, i)$, where $k$ is the time period when we decide to inspect. For example, the conditional plan $(i)$ is the plan where we inspect immediately (i.e., $\sigma^t_0$) and the conditional plan $(ni,i)$ is the plan where we wait a single time period before inspecting.

We now demonstrate how the belief probabilities will update over time in the event that the FDA observes $(nr)$. To do so, we first determine the probabilities of each observation occurring in each time period. Let $P[(dr)]^t$, $P[(cr)]^t$, and $P[(nr)]^t$, be the probabilities that we believe we will observe $(dr)$, $(cr)$ and $(nr)$, respectively, in time period $t$. We can compute these probabilities as 
\begin{align*}
    &P[(dr)]^t = p_{ND}b_N^{t-1} + p_{VD}b_V^{t-1} + p_{OD} b_O^{t-1}\text{,}\\
    &P[(cr)]^t = p_{NC}b_N^{t-1} + p_{VC}b_V^{t-1} + p_{OC} b_O^{t-1}\text{, and} \\
    &P[(nr)]^t = (p_{NO}+p_{NV}+p_{NN})b_N^{t-1}+(p_{VO}+p_{VV})b_V^{t-1}+p_{OO}b_O^{t-1}\text{.}
\end{align*}

Thus, in the event that the FDA observes $(nr)$, the belief probabilities update as $$b^{t+1}_N = \frac{p_{NN} b^t_N}{P[(nr)]^{t+1}}\text{, }b^{t+1}_V = \frac{p_{NV} b^t_N + p_{VV} b^t_V}{P[(nr)]^{t+1}}\text{, and } b^{t+1}_O = \frac{p_{NO} b^t_N + p_{VO} b^t_V + p_{OO} b^t_O}{P[(nr)]^{t+1}}\text{.}$$
We next construct the value functions of the model. Let $V^{t}(b^t, \sigma^t_j)$ to be the expected reward accumulated from time $t$ to $T$ given the current belief vector $b^t$ and implementing conditional plan $\sigma^t_j$. Value functions are defined as the sum of the expected reward for the state the FDA believes they are in at time $t$ and expected future rewards from time periods $t+1$ to $T$ following the actions dictated in $\sigma^t_j$. We start by defining the value functions that represent the end of the game. Recall that, if the facility has already been inspected, i.e., it is in state $I$, then the game is over and the FDA receives a reward of zero. We define $V^t((0,0,0,0,0,0,1), \sigma^t_j) = 0$ for any $\sigma^t_j$. The FDA receives a reward of $1$ if the facility is in state $N$, $V$ or $O$, a reward of of $-d$ if the facility is in state $D$, and a reward of $-c$ if the facility is in state $C$.  Note that, no reward is collected after an inspection has occurred, and inspecting incurs no cost. To match our assumption that a penalty for an unexpected disruptive event is only incurred once, if the facility is in state $D$ or $C$, it will then be optimal to ``inspect" to transition the Markov chain to state $I$, which indicates the game is over. This likewise reflects that a facility must be re-inspected by the FDA in order to restart production. Thus, for any $t$, $V^t((0,0,0,0,1,0,0), (i)) = -d$ and $V^t((0,0,0,0,0,1,0), (i)) = -c$, both of which are the optimal policy for the respective belief states.

What remains to be constructed are the value functions for when the facility is in operations, i.e., $b^t = (b^t_N, b^t_V, b^t_O, 0,0,0)$. Note that, when the facility is in operation and the FDA chooses to inspect, the FDA will receive the reward for the facility being open and then will not receive any future rewards. Thus, $V^t(b^t, (i)) = 1$ for any $t$. To define the remaining $V^t$, we first note that $V^{T-1}(b^{T-1},(ni,i)) = 1 - P[(dr)]^T d - P[(cr)]^T c + P[(nr)]^T V^T(b^{t+1},(i)) = 1 - P[(dr)]^T d - P[(cr)]^T c + P[(nr)]^T$. In this equation, the expected value accumulated in time periods $T-1$ and $T$ is the immediate value received for the facility being in operation in time period $T-1$ ($1$), and the expected reward for choosing to wait to inspect in time period $T$ ($- P[(dr)]^T d - P[(cr)]^T c + P[(nr)]^T$). 

Observe that $p_{NO}+p_{NV}+p_{NN}=1-p_{ND}+p_{NC}$, $p_{VO}+p_{VV}=1-p_{VD}-p_{VC}$ and $p_{OO}=1-p_{OD}-p_{OC}$. Thus, $V^{T-1}(b^{T-1},(ni,i)) = 1 + (1 - (d+1)p_{ND} - (c+1)p_{NC})b_N^{T-1} + (1 - (d+1)p_{VD} - (c+1)p_{VC})b_V^{T-1} + (1 - (d+1)p_{OD} - (c+1)p_{OC})b_O^{T-1}$. By defining $k_S = 1 - (d+1)p_{SD} - (c+1)p_{SC}$ for $S \in \{N,V,O\}$, we can more concisely say $V^{T-1}(b^{T-1},(ni,i)) = 1 + \sum_{S \in \{N,V,O\}}k_S b_S^{T-1}$.

We start by recursively defining the value functions for $\sigma^t_{T-t}$, when the FDA does not inspect until they are required to do so. For $t < T$, the value functions are defined as 
\begin{equation}
\label{eq:rec}
    V^{t}(b^t, \sigma_{T-t}^{t}) = 1 - P[(dr)]^{t+1} d - P[(cr)]^{t+1} c + P[(nr)]^{t+1} V^{t+1}(b^{t+1},\sigma_{T-t-1}^{t+1}). 
\end{equation}
Let $\Theta^i_{SS'}$ be the set of all paths from $S$ to $S'$ in $i-1$ time steps for $i \ge 2$ where $\theta^i \in \Theta^i_{S S'}$ is defined as $\theta^i = (\theta^i_1, \ldots, \theta^i_i)$ with $\theta^i_j$ as the state after $(j-1)$ time steps. Let $\Theta^i = \bigcup_{S, S' \in \{N,V,O\}} \Theta^i_{SS'}$. We can additionally define a closed form equivalent of $V^{t}(b^t, \sigma_{T-t}^{t})$ as:
\begin{equation}
\label{eq:nonrec}
    V^{t}(b^t, \sigma_{T-t}^{t}) = 1 + \sum_{S \in \{N,V,O\}} \sum_{\substack{S' \in \{N,V,O\} \\S' \text{ is reachable from }S}} k_{SS'}^{T-t} b_S^t,
\end{equation}
where 
\begin{minipage}[t]{0.45\linewidth}
    \begin{equation*}
    k_{SS'}^{T-t} = k_{S'} \sum_{\theta^{T-t} \in \Theta^{T-t}_{SS'}} f_{\theta^{T-t}},
\end{equation*}
\end{minipage}%
\begin{minipage}[t]{0.45\linewidth}
    \begin{equation*}
    f_{\theta_1^2,\theta_2^2} = \begin{cases}
        p_{\theta_1^2 \theta_1^2} + 1 \text{ if } \theta_1^2=\theta_2^2, \\
        p_{\theta_1^2 \theta_2^2} \text{ otherwise.} 
    \end{cases}
\end{equation*}
\end{minipage}

and
\begin{equation*}
    f_{\theta^{T-t}} = \begin{cases}
        p_{\theta_1^{T-t} \theta_2^{T-t}} f_{(\theta_2^{T-t}, \ldots, \theta_{T-t}^{T-t})} + 1 \text{ if }\theta_1^{T-t} =\theta_2^{T-t}=\theta_{T-t}^{T-t}, \\
        p_{\theta_1^{T-t} \theta_2^{T-t}} f_{(\theta_2^{T-t}, \ldots, \theta_{T-t}^{T-t})} \text{ otherwise,}
    \end{cases}
\end{equation*}

The proof of this equivalency, as well as all future proofs, appears in Appendix \ref{app:pfvar}.

\begin{thrm}
\label{thrm:equiv}
    Equation \eqref{eq:rec} is equivalent to Equation \eqref{eq:nonrec} for $1 \le t \le T-2$, i.e., there are closed form equations for the value functions that incorporate the dynamic nature of the equations explicitly.
\end{thrm}

We lastly need to define $V^t(b^t, \sigma_j^t)$ for $j \ne T-t$. Since no value is accumulated after an inspection occurs and the reward accumulated each time period is constant over the entire time horizon, the closed form definition of $V^t(b^t, \sigma_j^t)$ will be similar to the closed form definition of $V^{T-j}(b^{T-j}, \sigma_{T-j}^{T-j})$. This is because both conditional plans accumulated the rewards of $j$ time periods. The only adjustment that needs to be made is to use the belief probabilities at time $t$ instead of $T-j$. Thus, the closed form definition is $V^t(b^t, \sigma_j^t) = V^{T-j}(b^t, \sigma_{T-j}^{T-j})$. This completes the definition of the value functions.

\section{Analysis of Conditional Plans}
\label{sec:vf}
Similarly to \cite{cheng2023optimal},\cite{qiu2023optimal} and \cite{zhao2021optimal}, we demonstrate structural properties of our problem when we enforce additional structure on the sequences $(P[(dr)]^t)$, $(P[(cr)]^t)$ and $(P[(nr)]^t)$. We assume that $(P[(dr)]^t)$ and  $(P[(cr)]^t)$ are non-decreasing over $t$, and thus that $(P[(nr)]^t)$ is non-increasing over $t$. This is a reasonable assumption, as the quality of the facility will degrade over time, which in turn will make the facility more susceptible to a manufacturing failure. As the facility becomes more susceptible to a manufacturing failure, the facility is in turn more likely to halt operations to allow for non-mandatory maintenance to occur.

With this additional structure, we can directly compare the values of conditional plans. We first compare two conditional plans where one plan inspects in time period $j$ and the other inspects in time period $j+1$. Lemma \ref{lm:timeinc} states that once inspecting earlier has a greater expected value, it will continue to have a greater expected value for the remainder of the time horizon. Likewise, if inspecting later has a greater expected value, it will also have a greater expected value in earlier time periods.

\begin{lm}
\label{lm:timeinc}
    Suppose $(P[(dr)]^t)$ and $(P[(cr)]^t)$ are non-decreasing sequences and $(P[(nr)]^t)$ is a non-increasing sequence. If, for some $t$ and $j \le T-t$, $V^t(b^t, \sigma_j^t) \ge V^t(b^t, \sigma_{j+1}^t)$, then $V^{t+1}(b^{t+1}, \sigma_j^{t+1}) \ge V^{t+1}(b^{t+1}, \sigma_{j+1}^{t+1})$. If $V^t(b^t, \sigma_j^t) \le V^t(b^t, \sigma_{j+1}^t)$, then $V^{t-1}(b^{t-1}, \sigma_j^{t-1}) \le V^{t-1}(b^{t-1}, \sigma_{j+1}^{t-1})$.
\end{lm}

Lemma \ref{lm:timeinc} allows us to prove the following theorem about the optimal inspection time. Theorem \ref{thrm:nii} shows that, if there exists a time period $t$ where inspecting immediately is better than waiting to inspect until time period $t+1$, then waiting $j$ time periods after $t$ to inspect will be better than waiting $j+1$ time periods after $t$.

\begin{thrm}
\label{thrm:nii}
    Suppose $(P[(dr)]^t)$ and $(P[(cr)]^t)$ are non-decreasing sequences, $(P[(nr)]^t)$ is a non-increasing sequence and there exists a $t$ such that $V^t(b^t,(i)) \ge V^t(b^t, (ni,i))$. Then $V^t(b^t,\sigma_{j}^t) \ge V^t(b^t,\sigma_{j+1}^t)$ for $1 \le j < T-t$. 
\end{thrm}

From Theorem \ref{thrm:nii}, we can further deduce that the optimal inspection time can be determined by comparing $V^t(b^t,(i))$ and $V^t(b^t, (ni,i))$.

\begin{cor}
\label{cor:nii}
    Suppose $(P[(dr)]^t)$ and $(P[(cr)]^t)$ are non-decreasing sequences and $(P[(nr)]^t)$ is a non-increasing sequence and there exists a $t$ such that $V^t(b^t,(i)) \ge V^t(b^t, (ni,i))$. Then \newline$V^t(b^t, (ni,i)) \ge V^t(b^t, \sigma_j^t)$ for all $2 \le j \le T-t$.  
\end{cor}

Corollary \ref{cor:nii} indicates that, at any time period, we only need to compare the value functions for the plans $(i)$ and $(ni,i)$. Hence, making a simple comparison at each time period is sufficient to deduce the optimal inspection time for the overall problem, significantly reducing the complexity of the problem. If $V^t(b^t,(i)) \le V^t(b^t, (ni,i))$, then inspecting now is not optimal, and thus the immediate choice of not inspecting will be the decision in the optimal conditional plan. However, if $V^t(b^t,(i)) \ge V^t(b^t, (ni,i))$, then by Corollary \ref{cor:nii} we know that the optimal conditional plan is to inspect immediately. Lemma \ref{lm:timeinc} further confirms that this is indeed the correct decision, as it implies $V^{\bar{t}}(b^{\bar{t}},(i)) \ge V^{\bar{t}}(b^{\bar{t}}, (ni,i))$ for $t \le \bar{t} \le T$. This result is similar to the identification of critical thresholds in the mission abort literature, despite the difference in reward structures \citep{cheng2023optimal, qiu2023optimal, zhao2021optimal}. For any time period $t$, this reduces the number of conditional plans we need to consider from $T-t+1$ to $2$, drastically reducing the computational cost. Thus, we have an algorithm to identify the optimal inspection time, which requires at most $2T-2$ computations. The maximum number of computations is only required when the FDA does not inspect until time periods $T-1$ or $T$. Additionally, the FDA could integrate the score $V^t(b^t, (ni,i)) - V^t(b^t,(i))$ into their existing SSM model. This score indicates that a facility is more in need to being inspected as the score decreases. Given that the score is easy to compute and does not require the use of any complex algorithms, it will be simple to integrate into the SSM.

From Lemma \ref{lm:timeinc}, we can derive further structure to the value of the conditional plans. The expected value from inspecting in time period $t$ is non-decreasing as $t$ increases until the optimal inspection time $t^*$, then is non-increasing as $t$ increases.
\begin{cor}
\label{cor:structure}
    Let $t^*$ be the optimal inspection time, i.e., $V^1(b^1,\sigma^1_{t^*}) \ge V^1(b^1,\sigma^1_{t})$ for all $t$. Then $V^1(b^1,\sigma^1_{t}) \le V^1(b^1,\sigma^1_{t+1})$ for $1 \le t < t^*$ and $V^1(b^1,\sigma^1_{t}) \ge V^1(b^1,\sigma^1_{t+1})$ for $t^* \le t \le T-1$.
\end{cor}

\subsection{Sensitivity of Optimal Inspection Time to Penalty Parameters}
We additionally explore how varying $d$ and $c$ impact the optimal inspection time $t^*$. We consider two scenarios: one where we are interested in inspecting at time $t^*-1$, and one where $\bar{t}$ is a user-provided input and we want to determine what values $d$ and $c$ can take such that $t^* = \bar{t}$. For the first scenario, we derive an inequality that determines the minimum change in $d$ and $c$ such that inspecting at time $t^* - 1$ is strictly better than waiting until time $t^*$ to inspect. We analytically find that the optimal inspection time is more sensitive to change in penalty parameters as the probabilities of unexpected disruptive events increase. For the second case, we derive a pair of inequalities on $d$ and $c$ that ensure that $\bar{t}$ is the optimal inspection time. Our computational study of these inequalities shows that the range these penalty parameters can lie in decreases as $\bar{t}$ increases. More details about these inequalities are in Appendix \ref{app:sens}.

\section{Modeling Accounting for Inspection Outcomes}
\label{sec:modelvar}
We next consider a model where inspection may result in the facility needing to close for mandatory maintenance, resulting in a small penalty for the FDA. This model considers how FDA intervention can also cause a disruption that the FDA must account for. The FDA would force a facility that receives the OA classification to halt production; receiving the OA classification means that there are too many regulatory violations in the manufacturing process, and the drug produced is either ineffective or could cause harm when used. To incorporate this variation into our model, we add a new state $IC$, which can be transitioned to when the FDA inspects from the states $N$, $V$, or $O$. %Figure \ref{fig:statesIC} presents the updated Markov chain with the state $IC$.

%\begin{figure}[h]
%    \centering
%    \includegraphics[width=.45\linewidth]{statespaceIC_fix.png}
%    \caption{State space of the Markov chain when accounting for inspection outcomes}
%    \label{fig:statesIC}
%\end{figure}

Let $p_{SIC}$ be the probability of closing after being inspected in state $S$, and let $\Tilde{c} \ge 0$ be the penalty for being in state $IC$, where $\Tilde{c} \le c \le d$. We choose $\Tilde{c}$ in this manner since closure does increase vulnerability to shortages, but, unlike in the event of closing for non-mandatory maintenance, this closure is not unexpected by the FDA, and thus the FDA can take appropriate measures to mitigate further vulnerabilities to shortages. %Figure \ref{fig:outcomes} demonstrates how inspecting in states $N$, $V$, or $O$ can result in transitioning to $I$ or $IC$, where the values on each arc are the transition probabilities and the values to the right of $I$ and $IC$ are the value accumulated for being in that state.

%\begin{figure}[h]
%    \centering
%    \includegraphics[width=.45\linewidth]{outcomes_edit.jpg}
%    \caption{Transition probabilities and values accumulated when the FDA inspects a facility in operations}
%    \label{fig:outcomes}
%\end{figure}

With this adaptation, if $b_N^t+b_V^t+b_O^t = 1$, then $V^t(b^t, (i)) = 1 - \Tilde{c}(p_{NIC} b_N^t + p_{VIC} b_V^t + p_{OIC} b_O^t)$. Under this definition, if the facility is in operations in time period $t$, then inspecting will result in a penalty in expectation. Thus,
\begin{align*}
    V^{T-1}(b^{T-1},(ni,i)) &= 1 - P[(dr)]^T d - P[(cr)]^T c + P[(nr)]^T V^T(b^T,(i)) \\
    &= 1 + k_N b_N^{T-1} + k_V b_V^{T-1} + k_O b_O^{T-1} \\&- \Tilde{c}(p_{NIC} (p_{NN} b_N^{T-1}) + p_{VIC} (p_{NV} b_N^{T-1} + p_{VV} b_V^{T-1})\\
    &+ p_{OIC} (p_{NO} b_N^{T-1} + p_{VO} b_V^{T-1} + p_{OO} b_O^{T-1})) \\
    & = 1 + (k_N - \Tilde{c}(p_{NN}p_{NIC} + p_{NV}p_{VIC} + p_{NO}p_{OIC}))b_N^{T-1} \\
    &+ (k_V - \Tilde{c}(p_{VV}p_{VIC}+p_{VO}p_{OIC}))b_V^{T-1} + (k_O -\Tilde{c}p_{OO}p_{OIC})b_O^{T-1}
\end{align*}
Let $\Tilde{c}_{SS'IC}^i = \Tilde{c} p_{S'IC} \sum_{\theta^i \in \Theta^i_{SS'}}\left(\prod_{j=1}^{i-1}p_{\theta_{j}^i \theta_{j+1}^i}\right)$. Under this adaptation, we can again define a closed form equivalent for the value functions as:
\begin{equation}
\label{eq:nonrecAlt}
    V^{t}(b^{t}, \sigma_{t}^{T-t}) = 1 + \sum_{S \in \{N,V,O\}} \sum_{\substack{S' \in \{N,V,O\} \\S' \text{ is reachable from }S}} (k_{SS'}^{T-t} - \Tilde{c}_{SS'IC}^{T-t+1}) b_S^{t}.
\end{equation}

\begin{thrm}
\label{thrm:equivAlt}
    Under the definition of $V^t(b^t, (i)) = 1 - \Tilde{c}(p_{NIC} b_N^t + p_{VIC} b_V^t + p_{OIC} b_O^t)$ when $b^t_N+b^t_V+b^t_O = 1$, Equation \eqref{eq:rec} is equivalent to Equation \eqref{eq:nonrecAlt} for $i \ge 2$.
\end{thrm}

Let $P[(cr)|(i)]^t = p_{NIC} b_N^t + p_{VIC} b_V^t + p_{OIC} b_O^t$  i.e., $P[(cr)|(i)]^t$ is the probability we believe that inspecting in time period $t$ will result in the facility needing to close for mandatory maintenance. To prove an equivalent version of Theorem \ref{thrm:nii}, we need an equivalent version of Lemma \ref{lm:timeinc}. Let $d = \alpha_d \Tilde{c}$ and $c = \alpha_c \Tilde{c}$, and note that $\alpha_d \ge 1$ and $\alpha_c \ge 1$. As before, we assume $(P[(dr)]^t)$, $(P[(cr)]^t)$ and $(P[(cr)|(i)]^t)$ are non-decreasing sequences and $(P[(nr)]^t)$ is a non-increasing sequence. We additionally assume that $(P[(nr)]^{t} P[(cr)|(i)]^{t})$ is a non-decreasing sequence, i.e., the probability of a facility needing to close after an inspection increases faster than probability of no issues being reported decreases. Since $P[(nr)]^{t} + P[(dr)]^{t} + P[(cr)]^{t} = 1$ for all $t$, this assumption is equivalent to saying, the probability of a facility closing after an inspection increases faster than the probability of a facility experiencing an unexpected disruptive event increases. This is a reasonable assumption since we expect the FDA to be as sensitive to CGMP violations as facility managers would be.

\begin{lm}
    \label{lm:timeincAlt}
    Suppose $(P[(dr)]^t)$, $(P[(cr)]^t)$ and $(P[(cr)|(i)]^t)$ are non-decreasing sequences and $(P[(nr)]^t)$ is a non-increasing sequence. Additionally, suppose that $(P[(nr)]^{t} P[(cr)|(i)]^{t})$ is a non-decreasing sequence and $\alpha_d (P[(dr)]^{t+1} - P[(dr)]^t) + \alpha_c (P[(cr)]^{t+1} - P[(cr)]^t) \ge (P[(cr)|(i)]^{t+1} - P[(cr)|(i)]^t)$. If, for some $t$ and $j \le T-t$, $V^t(b^t, \sigma_{j}^t) \ge V^t(b^t,\sigma_{j+1}^t)$, then $V^{t+1}(b^{t+1}, \sigma_{j}^{t+1}) \ge V^{t+1}(b^{t+1},\sigma_{j+1}^{t+1})$.
\end{lm}

We note that the bound $\alpha_d (P[(dr)]^{t+1} - P[(dr)]^t) + \alpha_c (P[(cr)]^{t+1} - P[(cr)]^t) \ge (P[(cr)|(i)]^{t+1} - P[(cr)|(i)]^t)$ is not tight, but instead provides a nice practical interpretation. As long as the expected penalty of closure due to inspection is increasing at a slower rate than the expected penalty of a manufacturing failure or unexpected closure occurring, then a conditional plan $\sigma_{j}^t$, once preferred over $\sigma_{j+1}^t$, will remain preferred as $t$ increases. We now prove an equivalent version of Theorem \ref{thrm:nii}.

\begin{thrm}
\label{thrm:niiAlt}
    Suppose $(P[(dr)]^t)$, $(P[(cr)]^t)$ and $(P[(cr)|(i)]^t)$ are non-decreasing sequences and $\\(P[(nr)]^t)$ is a non-increasing sequence. Additionally, suppose that $(P[(nr)]^{t} P[(cr)|(i)]^{t})$ is a non-decreasing sequence, $\alpha_d (P[(dr)]^{t+1} - P[(dr)]^t) + \alpha_c (P[(cr)]^{t+1} - P[(cr)]^t) \ge (P[(cr)|(i)]^{t+1} - P[(cr)|(i)]^t)$. If there exists a $t$ such that $V^t(b^t,(i)) \ge V^t(b^t, (ni,i))$, then $V^t(b^t, \sigma^t_j) \ge V^t(b^t,\sigma_{j+1}^t)$ for $1 \le j \le T-t$.
\end{thrm}

Theorem \ref{thrm:niiAlt} again allows us to further deduce an equivalent version of Corollary \ref{cor:nii}, and we can thus again determine the optimal inspection time in this variant model solely by comparing $V^t(b^t,(i))$ and $V^t(b^t,(ni,i))$.

\begin{cor}
\label{cor:niiAlt}
    Suppose $(P[(dr)]^t)$, $(P[(cr)]^t)$ and $(P[(cr)|(i)]^t)$ are non-decreasing sequences and $\\(P[(nr)]^t)$ is a non-increasing sequence.  Additionally, suppose that $(P[(nr)]^t P[(cr)|(i)]^{t})$ is a non-decreasing sequence, $\alpha_d (P[(dr)]^{t+1} - P[(dr)]^t) + \alpha_c (P[(cr)]^{t+1} - P[(cr)]^t) \ge (P[(cr)|(i)]^{t+1} - P[(cr)|(i)]^t)$. If there exists a $t$ such that $V^t(b^t,(i)) \ge V^t(b^t, (ni,i))$, then $V^t(b^t, (ni,i)) \ge V^t(b^t, \sigma_j^t)$ for all $2 \le j \le T-t$.  
\end{cor}

As with the base model, we only need to compare the value of inspecting immediately against the value of inspecting after waiting one time period to determine the optimal solution. This again reduces our problem to a series of simple computations to determine when to inspect. Thus, the algorithm proposed in Section \ref{sec:vf} will still identify the optimal inspection time, and the score $V^t(b^t,(ni,i)) - V^t(b^t, (i))$ can still be integrated into the FDA's existing SSM model.

\section{Computational Study}
\label{sec:numres}
We simulate how the quality of a facility will degrade over time to test how effective our model is at maximizing the value of a facility while being cognizant of the penalties of unexpected disruptive events occurring. For these experiments, we assume the facility starts in the No Official Action ($N$) state with full certainty, i.e., $b^1=(1,0,0,0,0,0)$, and we iterate the Markov chain until an unexpected disruptive event occurs. We choose this starting belief probability under the assumption that the facility was just inspected and given the No Action Indicated classification, so the state is known with full certainty. We say an inspection rule \emph{catches an unexpected disruptive event} if the time of inspection calculated by the inspection rule is strictly less than the time the Markov chain reaches either states $D$ or $C$ in the simulation, which we denote to be $t_F$. We follow the procedure outlined in Section \ref{sec:vf} to identify the optimal inspection time. Let $t_V$ be the optimal inspection time when not accounting for inspection outcomes, and let $t_{VC}$ be the optimal inspection time when accounting for inspection outcomes. As a baseline, we compare these inspection times against inspecting before the expected time to reach states $D$ or $C$ from the state $N$. Since we cannot compare against SSM-recommended inspection times, this serves as the best existing means of comparison that accounts for qualities intrinsic to the facility and inherent risks from production. The expected time to reach state $D$ or $C$ from a state $S$ is computed as $\mu_S = 1 + \sum_{S' \in \{N,V,O\}: p_{SS'}>0} p_{SS'} \mu_{S'}$ \citep{bremaud2020markov}. Then given the expected time to reach state $D$ or $C$, in the baseline setting we inspect in the time period $t_E = \lceil \mu_N \rceil - 1$, which is the largest integer that is strictly less than $\mu_N$. We refer to this rule as the ETD inspection rule. We refer to these three inspection rules as the risk-based inspection rules.

We perform $100,000$ experiments for each set of fixed penalty parameters, and we record the percentage of instances where each inspection rule inspects before an unexpected disruptive event occurs and the total accumulated value for each instance. We additionally compare the total accumulated values against the values accumulated when there is no inspection, i.e., when the Markov chain runs until it reaches either states $D$ or $C$, and when the FDA inspects after a fixed number of time periods: $T_1=24$, $T_2=60$ and $T_3=120$. If each time period were to be representative of a month, then these time horizons would be equivalent to two years, five years and ten years, respectively. These comparisons shed light on the value of planning inspections based on risk. Numerical experiments were conducted in R on a laptop with an Intel\textsuperscript{\textregistered} Core\textsuperscript{TM} i7-11370H processor and 16 GB RAM running Windows 11. The iGraph package \citep{csardi2006igraph} was used to conduct the penalty parameter sensitivity analysis.

\subsection{Catching Unexpected Disruptive Events}
\label{ssec:base}
We explore how different transition probability matrices and penalty parameters will impact the percent of unexpected disruptive events caught and the total value accumulated over time. We first analyze the results of our inspection rules with a fixed transition probability matrix and fixed penalty parameters $d$ and $c$. For these experiments, we set the penalty parameters to be $c=5$ and $\Tilde{c}=1$, and we set the probabilities of closure due to inspection to be $p_{NIC}=0$, $p_{VIC}=0.3$ and $p_{OIC}=1$. A facility classified as No Action Needed has no need to close for maintenance,
while a facility classified as Official Action Needed will have regulatory action taken against it if it
does not close for maintenance. We speculate that a facility classified as Voluntary Action Needed
will be more likely to stay in operations, as closing for non-mandatory maintenance will reduce
profits through maintenance costs and lost sales. We choose $\Tilde{c}$ such that the assumptions outlined in Section 5 are satisfied for all tested values of $d$. We test three different probability matrices with different expected times to unexpected disruptive events and various penalty parameters $d$, all of which satisfy the assumptions in Theorems \ref{thrm:nii} and \ref{thrm:niiAlt}. The first matrix we test is 

\footnotesize $$ P^1(ni) = \begin{bmatrix}
0.9125 & 0.0875 & 0 & 0 & 0 & 0 \\
0 & 0.825 & 0.1125 & 0.045 & 0.0175 & 0 \\
0 & 0 & 0.75 & 0.175 & 0.075 & 0 \\
0 & 0 & 0 & 1 & 0 & 0 \\
0 & 0 & 0 & 0 & 1 & 0 \\
0 & 0 & 0 & 0 & 0 & 1 
\end{bmatrix}. $$\normalsize

Table \ref{tbl:inspect_p1} presents the recommended inspection times for each inspection rule with varying values for $d \in \{14, 18, 22, 26, 30\}$. These choices of $d$ provided a wide variety of optimal inspection times.

\begin{table}[h!]
\centering\small 
\begin{tabular}{|c|c|c|c|}
\hline
$d$  & $t_V$ & $t_{VC}$ & $t_E$ \\ \hline
14 & 27 & 31 & 19 \\ \hline
18 & 15 & 16 & 19 \\ \hline
22 & 12 & 12 & 19 \\ \hline
26 & 10 & 10 & 19 \\ \hline
30 & 8 & 8 & 19 \\ \hline
\end{tabular}
\caption{Recommended Inspection Time Period for $P^1(ni)$}
\label{tbl:inspect_p1}
\end{table}\normalsize

Table \ref{tbl:value_p1} presents the average value accumulated with each inspection rule, as well as the percent of unexpected disruptive events caught. We use the results of $t_E$ as a baseline, and, for each other inspection rule, we provide the absolute difference for percent of unexpected disruptive events caught and the relative difference for the average value accumulated. For each rule, we present the value accumulated if we are not accounting for the outcome of the inspection (no IC) and the value accumulated if we are accounting for the outcome of the inspection (IC). We additionally present statistics around the time until an unexpected disruptive event occurs in in Table \ref{tbl:endtime_p1}.

\begin{table}[h!]
\centering
\resizebox{\textwidth}{!}{\begin{tabular}{|c|c|c|c|c|c|c|c|}
\hline
Rule          & $t_E$ & $t_V$ & $t_{VC}$ & 24 & 60 & 120 & NoIns \\ \hline\hline
\% Caught, $d=14$ & 44.297\% & -21.084\% & -27.857\% & -14.459\% & -43.103\% & -44.265\% & -  \\ \hline
Value (no IC), $d=14$ & 8.3402 & +1.580\% & - & +1.501\% & +0.277\% & +0.170\% & +0.171\% \\ \hline
Value (IC), $d=14$ & 8.2012 & - & +2.692\% & +1.980\% & +1.923\% & +1.868\% & - \\ \hline\hline
\% Caught, $d=18$ & 44.438\% & +14.944\% & +10.900\% & -14.534\% & -43.195\% & -44.433\% & -  \\ \hline
Value (no IC), $d=18$ & 6.7869 & +1.798\% & - & -4.114\% & -17.156\% & -17.734\% & -17.739\% \\ \hline
Value (IC), $d=18$ & 6.6488 & - & +1.558\% & -3.621\% & -15.464\% & -16.022\% & - \\ \hline\hline
\% Caught, $d=22$ & 44.351\% & +27.165\% & +27.165\% & -14.587\% & -43.177\% & -44.344\% & -  \\ \hline
Value (no IC), $d=22$ & 5.1509 & +14.801\% & - & -13.884\% & -46.358\% & -47.779\% & -47.787\% \\ \hline
Value (IC), $d=22$ & 5.0123 & - & +14.409\% & -13.382\% & -44.961\% & -46.336\% & - \\ \hline\hline
\% Caught, $d=26$ & 44.700\% & +35.326\% & +35.326\% & -14.583\% & -43.484\% & -44.696\% & -  \\ \hline
Value (no IC), $d=26$ & 3.7237 & +43.271\% & - & -30.341\% & -98.751\% & -101.641\% & -101.646\% \\ \hline
Value (IC), $d=26$ & 3.5834 & - & +43.955\% & -30.404\% & -98.822\% & -101.705\% & - \\ \hline\hline
\% Caught, $d=30$ & 44.029\% & +43.628\% & +43.628\% & -14.329\% & -42.828\% & -44.026\% & -  \\ \hline
Value (no IC), $d=30$ & 1.9024 & +155.435\% & - & -77.991\% & -251.477\% & -259.425\% & -259.446\% \\ \hline
Value (IC), $d=30$ & 1.7651 & - & +165.991\% & -81.984\% & -263.503\% & -271.826\% & - \\ \hline
\end{tabular}}
\caption{Percent of Unexpected Disruptive Events Caught and Average Value Accumulated with Each Inspection Rule for $P^1(ni)$}
\label{tbl:value_p1}
\end{table}

\begin{table}[h!]
\centering \small 
\begin{tabular}{|c|c|c|c|c|c|}
\hline
$d$ & 14 & 18 & 22 & 26 & 30 \\\hline
Average End Time & 20.7118 & 20.7540 & 20.6935 & 20.7876 & 20.6297 \\ \hline
Std. Dev. End Time & 12.5409 & 12.6105 & 12.5441 & 12.5777 & 12.4737 \\ \hline
Median End Time & 18 & 18 & 18 & 18 &18  \\ \hline
Min End Time & 3 & 3 & 3 & 3 & 3 \\ \hline
Max End Time & 143 & 147 & 135 & 145 & 145 \\ \hline
\end{tabular}
\caption{Descriptive Statistics of Time of Unexpected Disruptive Event for $P^1(ni)$}
\label{tbl:endtime_p1}
\end{table}\normalsize

\begin{figure}
    \centering
    \includegraphics[width=0.5\textwidth]{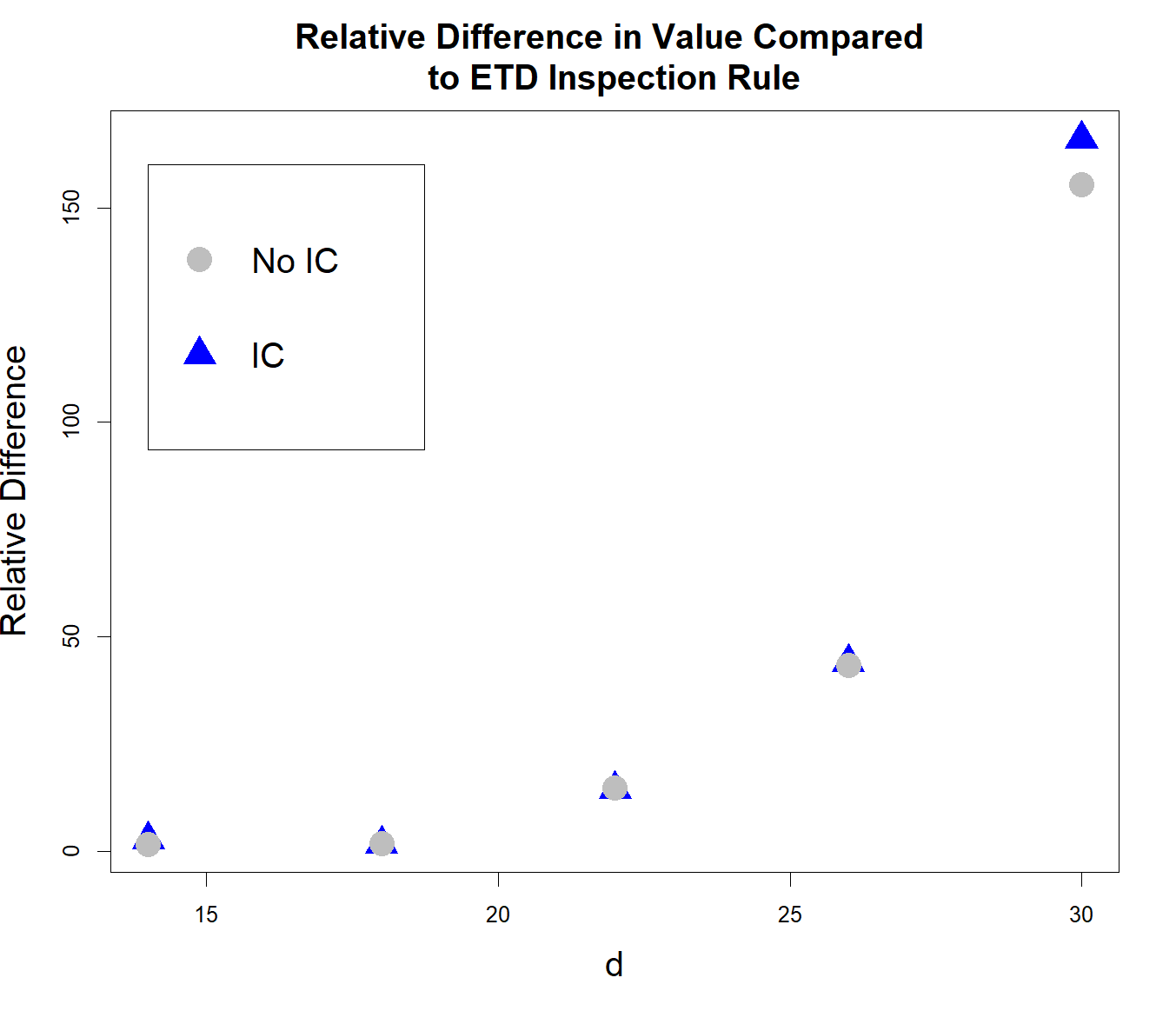}
    \caption{Relative difference in value between the ETD inspection rule and optimal inspection rules for $P^1(ni)$}
    \label{fig:reldiff}
\end{figure}

From these tables, we can see that, for $d=14$, while the ETD inspection rule catches more unexpected disruptive events than the optimal inspection rules, the value accumulated is less than that of other rules. This scenario is analogous to the case where the facility experiencing an unexpected disruptive event is not impactful enough to warrant more frequent inspections. As the penalty for a manufacturing failure increases, the optimal inspection rules inspect earlier and catch more unexpected disruptive events than ETD (and the other inspection rules), which is fixed as it does not take into account the penalty of unexpected disruptive events. The optimal inspection rules additionally accumulated more value on average than the other inspection rules, indicating that they are better able to account for the trade-offs in the problem. Figure \ref{fig:reldiff} graphically presents the relative differences in average value accumulated between the ETD inspection rule and optimal inspection rules. From this, we can see a dramatic trend in the increase in the relative difference. For smaller values of $d$, the values accumulated between the different rules are fairly similar, and the relative difference for both optimal inspection rules grows super-linearly as $d$ increases. This highlights how much more important it is for the FDA to prioritize inspecting pharmaceutical facilities that manufacture products more critical to public health.

Additionally, when accounting for the outcome of inspecting, the delay in inspecting can result in a slight increase in value even though less unexpected disruptive events are caught. In particular, when $d=14$, not accounting for this closure possibility caught approximately $7\%$ more unexpected disruptive events while accumulating less value, and when $d=18$, not accounting for this closure caught approximately $4\%$ more unexpected disruptive events while accumulating a similar value. These small differences in the number of unexpected disruptive events caught, as compared to the difference in average value accumulated, indicates that the FDA can take a calculated risk in waiting to inspect longer if the penalties for unexpected disruptive events are relatively small compared to the penalties of the facility closing for mandatory maintenance.

We next want to understand the robustness of our inspection strategy. For the next set of experiments, we randomly generate a probability matrix $\hat{P}$ where each entry $\hat{p}_{ij}$ is randomly sampled from a normal distribution with mean $p_{ij}$ and a given standard deviation $s$. We then modify $\hat{P}$ to be $\Tilde{P}$ to ensure that each row in the matrix is indeed a probability distribution and that the states $D$ and $C$ can be transitioned to from at least state $O$. We set $\Tilde{p}_{SS} = \min\{\hat{p}_{SS},0.995\}$, and $\Tilde{p}_{SS'} = \max\{0,\min\{\hat{p}_{SS'}, 1-\sum_{S'' \in \{N,V,O\}: \hat{p}_{SS''}>0} \hat{p}_{SS''}\}\}$. To ensure that each row sums to one, we do not randomly generate the transition probability for the last nonzero entry in the row $\hat{p}_{SS_F}$ and instead use $\Tilde{p}_{SS_F} = 1-\sum_{S'' \in \{N,V,O\}: \Tilde{p}_{SS''}>0} \Tilde{p}_{SS''}$. The matrix $\Tilde{P}$ is the transition probability matrix used. The inspection times considered are still the inspection times presented in Table \ref{tbl:inspect_p1}, indicating that the FDA used the expected value for transition probabilities when calculating the optimal inspection time.

We test two different standard deviations $s \in \{0.01, 0.02\}$ to see how the uncertainty around the transition probabilities impacts the number of unexpected disruptive events caught and total value accumulated. Tables \ref{tbl:value_p1} and \ref{tbl:endtime_p1_s1} present results and descriptive statistics for experiments with $s=0.01$, and Tables \ref{tbl:value_p1_s2} present results and descriptive statistics for experiments with $s=0.02$. %A standard deviation of $0.01$ results in matrices similar to $P^1(ni)$, where a standard deviation of $0.02$ results in a wide variety of matrices; for illustrative purposes, Figure \ref{fig:pnn} in Appendix \ref{app:acs} plots the probability density functions with mean $p_{NN}$. Tables \ref{tbl:value_p1} and \ref{tbl:endtime_p1_s1} present results and descriptive statistics for experiments with $s=0.01$, and Tables \ref{tbl:value_p1_s2} present results and descriptive statistics for experiments with $s=0.02$

\begin{table}[h!]
\centering
\resizebox{\textwidth}{!}{\begin{tabular}{|c|c|c|c|c|c|c|c|}
\hline
Rule          & $t_E$ & $t_V$ & $t_{VC}$ & 24 & 60 & 120 & NoIns \\ \hline\hline
\% Caught, $d=14$ & 44.940\% & -21.091\% & -27.832\% & -14.546\% & -43.533\% & -44.933\% & -  \\ \hline
Value (no IC), $d=14$ & 8.4838 & +2.132\% & - & +1.695\% & +2.208\% & +2.130\% & +2.219\% \\ \hline
Value (IC), $d=14$ & 8.3447 & - & +3.376\% & +2.193\% & +3.848\% & +3.831\% & - \\ \hline\hline
\% Caught, $d=18$ & 44.611\% & +14.713\% & +10.864\% & -14.443\% & -43.240\% & -44.596\% & -  \\ \hline
Value (no IC), $d=18$ & 6.8494 & +1.216\% & - & -3.811\% & -15.794\% & -16.283\% & -16.292\% \\ \hline
Value (IC), $d=18$ & 6.7123 & - & +1.110\% & -3.331\% & -14.143\% & -14.575\% & - \\ \hline\hline
\% Caught, $d=22$ & 44.598\% & +26.999\% & +26.999\% & -14.461\% & -43.200\% & -44.581\% & -  \\ \hline
Value (no IC), $d=22$ & 5.3079 & +12.415\% & - & -12.756\% & -43.618\% & -45.018\% & -45.042\% \\ \hline
Value (IC), $d=22$ & 5.1693 & - & +12.004\% & -12.354\% & -42.198\% & -43.546\% & - \\ \hline\hline
\% Caught, $d=26$ & 44.856\% & +35.119\% & +35.119\% & -41.820\% & -43.471\% & -44.846\% & -  \\ \hline
Value (no IC), $d=26$ & 3.8274 & +39.460\% & - & -28.562\% & -91.605\% & -94.686\% & -94.717\% \\ \hline
Value (IC), $d=26$ & 3.6879 & - & +39.969\% & -28.599\% & -91.415\% & -94.485\% & - \\ \hline\hline
\% Caught, $d=30$ & 44.372\% & +43.165\% & +43.165\% & -14.241\% & -43.001\% & -44.364\% & -  \\ \hline
Value (no IC), $d=30$ & 2.1102 & +129.703\% & - & -68.377\% & -220.609\% & -228.021\% & -228.054\% \\ \hline
Value (IC), $d=30$ & 1.9728 & - & +137.358\% & -71.300\% & -229.253\% & -236.942\% & - \\ \hline
\end{tabular}}
\caption{Percent of Unexpected Disruptive Events Caught and Average Value Accumulated with Each Inspection Rule for $\Tilde{P}^1(ni)$ with $s = 0.01$}
\label{tbl:value_p1_s1}
\end{table}

\begin{table}[h!]
\centering\small 
\begin{tabular}{|c|c|c|c|c|c|}
\hline
$d$ & 14 & 18 & 22 & 26 & 30  \\\hline
Average End Time & 20.9702 & 20.8595 & 20.8659 & 20.9291 & 20.8380 \\ \hline
Std. Dev. End Time & 12.8893 & 12.8529 & 12.8729 & 12.9012 & 12.8479 \\ \hline
Median End Time & 18 & 18 & 18 & 18 & 18 \\ \hline
Min End Time & 3 & 3 & 3 & 3 & 3 \\ \hline
Max End Time & 151 & 173 & 158 & 149 & 166 \\ \hline
\end{tabular}
\caption{Descriptive Statistics of Time of Unexpected Disruptive Event for $\Tilde{P}^1(ni)$, $s=0.01$}
\label{tbl:endtime_p1_s1}
\end{table}\normalsize

\begin{table}[h!]
\centering
\resizebox{\textwidth}{!}{\begin{tabular}{|c|c|c|c|c|c|c|c|}
\hline
Rule          & $t_E$ & $t_V$ & $t_{VC}$ & 24 & 60 & 120 & NoIns \\ \hline\hline
\% Caught, $d=14$ & 45.704\% & -20.928\% & -27.574\% & -14.213\% & -43.617\% & -45.636\% & -  \\ \hline
Value (no IC), $d=14$ & 8.7804 & +3.507\% & - & +2.840\% & +6.925\% & +7.885\% & +8.068\% \\ \hline
Value (IC), $d=14$ & 8.6422 & - & +5.335\% & +3.313\% & +8.563\% & +9.608\% & - \\ \hline\hline
\% Caught, $d=18$ & 45.722\% & +14.214\% & +10.325\% & -14.249\% & -43.665\% & -45.651\% & -  \\ \hline
Value (no IC), $d=18$ & 7.3038 & -0.998\% & - & -1.941\% & -7.977\% & -7.733\% & -7.614\% \\ \hline
Value (IC), $d=18$ & 7.1650 & - & -0.844\% & -1.439\% & -6.282\% & -5.947\% & - \\ \hline\hline
\% Caught, $d=22$ & 45.807\% & +26.309\% & +26.309\% & -14.141\% & -43.674\% & -45.728\% & -  \\ \hline
Value (no IC), $d=22$ & 5.8957 & +6.666\% & - & -8.313\% & -29.164\% & -29.527\% & -29.421\% \\ \hline
Value (IC), $d=22$ & 5.7476 & - & +6.361\% & -7.704\% & -27.443\% & -27.712\% & - \\ \hline\hline
\% Caught, $d=26$ & 45.516\% & +34.782\% & +34.792\% & -14.138\% & -43.518\% & -45.462\% & -  \\ \hline
Value (no IC), $d=26$ & 4.3086 & +30.226\% & - & -21.090\% & -67.614\% & -69.352\% & -69.134\% \\ \hline
Value (IC), $d=26$ & 4.1710 & - & +30.283\% & -20.894\% & -66.684\% & -68.343\% & - \\ \hline\hline
\% Caught, $d=30$ & 45.722\% & +42.033\% & +42.033\% & -14.252\% & -43.641\% & -45.654\% & -  \\ \hline
Value (no IC), $d=30$ & 2.9540 & +71.364\% & - & -44.259\% & -141.134\% & -145.833\% & -145.650\% \\ \hline
Value (IC), $d=30$ & 2.8141 & - & +74.180\% & -45.073\% & -143.385\% & -148.118\% & - \\ \hline
\end{tabular}}
\caption{Percent of Unexpected Disruptive Events Caught and Average Value Accumulated with Each Inspection Rule for $\Tilde{P}^1(ni)$ with $s = 0.02$}
\label{tbl:value_p1_s2}
\end{table}

\begin{table}[h!]
\centering\small 
\begin{tabular}{|c|c|c|c|c|c|}
\hline
$d$ & 14 & 18 & 22 & 36 & 30 \\\hline
Average End Time & 21.5801 & 21.5648 & 21.6401 & 21.5522 & 21.5749 \\ \hline
Std. Dev. End Time & 14.3624 & 14.1660 & 14.3501 & 14.3273 & 14.2560 \\ \hline
Median End Time & 18 & 18 & 18 & 18 & 18 \\ \hline
Min End Time & 3 & 3 & 3 & 3 & 3 \\ \hline
Max End Time & 511 & 261 & 335 & 706 & 431 \\ \hline
\end{tabular}
\caption{Descriptive Statistics of Time of Unexpected Disruptive Event for $\Tilde{P}^1(ni)$, $s=0.02$}
\label{tbl:endtime_p1_s2}
\end{table} \normalsize

When $s=0.01$, the results regarding unexpected disruptive events caught and value accumulated are similar to the case when $P$ is known with certainty, indicating that the optimal inspection rules are still able to accurately determine if the facility is critical enough to monitor more frequently. However, the value accumulated for fixed inspection times seems to increase more than the increase in the optimal inspection rules, particularly more as $T$ increases. This is likely due to instances where the randomly generated matrix is more likely to transition to the state it is currently in as compared to $P^1$, resulting in more instances with longer times to reach either state $D$ or $C$, which increases the average value. When $s=0.02$ and $d=14$, we see that the value for waiting to inspect until $T \ge 60$ accumulates more value than when implementing the optimal inspection rules. This indicates that, for less critical facilities, being less certain about the how the quality of the facility will decline over time proves to be a challenge for risk-based tools, so using simpler metrics may be more effective.

To help understand why the optimal inspection rules struggle as we become less certain about the transition probabilities, we explore what the ETD inspection rule would recommend if it were computed for each randomly generated matrix. Figure \ref{fig:p1_tETD_d18} in Appendix \ref{app:acs} presents histograms of the time periods the ETD inspection rule decides to inspect for $s=0.01$ (Figure \ref{fig:p1_tETD_d18_sd1}) and $s=0.02$ (Figure \ref{fig:p1_tETD_d18_sd2}), respectively, for the matrices generated when testing $d=18$. From these figures, it is clear to see that, for $s=0.01$, the deviation in when the ETD inspection rules would decide to inspect is relatively small, suggesting that the randomly generated transition matrices are similar enough that the optimal inspection rules still accurately capture the trade-off between waiting to inspect and preventing an unexpected disruptive event. However, when $s=0.02$, the number of matrices with a significantly different recommendation from the ETD inspection rule increases dramatically, meaning that the optimal inspection rules no longer accurately capture this trade-off. Recall that the ETD inspection rule recommends inspecting in time period $19$ for $P^1(ni)$. For the randomly generated matrices, over $90\%$ of the ETD inspection rule recommendations suggest inspecting in the range $[17,21]$ when $s=0.01$, while less than $65\%$ of the recommendations are in this range when $s=0.02$. The generated transition matrices with $s=0.02$ are also expected to reach states $D$ or $C$ significantly later than the base matrix. When $s=0.02$, over $10\%$ of the ETD inspection rule recommendations suggest inspecting in time period $24$ or later, while only $1\%$ of the recommendations suggest inspecting in time period $24$ or later when $s=0.01$. This suggests there is significant value in the FDA trying to better understand how the quality of facilities will degrade over time, as even a small decrease in standard deviation will lead to much more accurate recommendations from risk-based tools. 

We additionally wish to see how varying the likelihood of transitioning to the state the facility is currently in for the base transition probability matrix we use impacts different aspects of the model. The second matrix we test is more prone to unexpected disruptive events, as the Markov chain can transition to the state $D$ from any operational state and it is less likely to transition to the state it is currently in. Table \ref{tbl:inspect_p2} presents the recommended inspection times for each inspection rule with the varying values for $d \in \{10, 14, 18\}$. We test different values of $d$ for $P^2(ni)$ since smaller values of $d$ resulted in catching more unexpected disruptive events than $P^1(ni)$.

%\small $$ P^2(ni) = \begin{bmatrix}
%0.835 & 0.135 & 0.02 & 0.01 & 0 & 0 \\
%0 & 0.795 & 0.1425 & 0.045 & 0.0175 & 0 \\
%0 & 0 & 0.7 & 0.2 & 0.1 & 0 \\
%0 & 0 & 0 & 1 & 0 & 0 \\
%0 & 0 & 0 & 0 & 1 & 0 \\
%0 & 0 & 0 & 0 & 0 & 1 
%\end{bmatrix}. $$\normalsize

\begin{table}
\begin{minipage}[t]{0.495\linewidth}
\centering\small
\begin{tabular}{|c|c|c|c|}
\hline
$d$  & $t_V$ & $t_{VC}$ & $t_E$ \\ \hline
10 & 11 & 11 & 12 \\ \hline
14 & 7 & 7 & 12 \\ \hline
18 & 6 & 6 & 12 \\ \hline
\end{tabular}
\caption{Recommended Inspection Time Period for $P^2(ni)$}
\label{tbl:inspect_p2}\normalsize
\end{minipage}
\begin{minipage}[t]{0.495\linewidth}
\centering \small 
\begin{tabular}{|c|c|c|c|}
\hline
$d$  & $t_V$ & $t_{VC}$ & $t_E$ \\ \hline
18 & 84 & Inf & 25 \\ \hline
22 & 24 & 24 & 25 \\ \hline
26 & 17 & 17 & 25 \\ \hline
30 & 14 & 14 & 25 \\ \hline
34 & 12 & 12 & 25 \\ \hline
\end{tabular}
\caption{Recommended Inspection Time Period for $P^3(ni)$}
\label{tbl:inspect_p3}
\end{minipage}
\end{table}

Table \ref{tbl:value_p2} presents the number of unexpected disruptive events caught by each inspection rule and average value accumulated. The relative differences in the average value accumulated by the optimal inspection rules follow similar trends with $P^1(ni)$. Because the facility is more prone to experiencing an unexpected disruptive event, waiting longer to inspect proves to be more detrimental to the average value accumulated than with $P^1(ni)$, even with smaller values of $d$. Additionally, when $d$ is large enough, waiting too long can cause the average value to become increasingly negative, indicating that the facility needs to be monitored more frequently to mitigate the severely detrimental impact of an unexpected disruptive event.  %, and Table \ref{tbl:endtime_p2} presents descriptive statistics regarding the time the Markov chain terminated. 

\begin{table}[h!]
\centering
\resizebox{\textwidth}{!}{\begin{tabular}{|c|c|c|c|c|c|c|c|}
\hline
Rule          & $t_E$ & $t_V$ & $t_{VC}$ & 24 & 60 & 120 & NoIns\\ \hline\hline
\% Caught, $d=10$ & 46.193\% & +5.629\% & +5.629\% & -37.572\% & -46.177\% & -46.193\% & -  \\ \hline
Value (no IC), $d=10$ & 4.6667 & +0.351\% & - & -12.192\% & -17.444\% & -17.459\% & -17.459\% \\ \hline
Value (IC), $d=10$ & 4.4784 & - & +0.051\% & -9.503\% & -13.975\% & -13.989\% & - \\ \hline\hline
\% Caught, $d=14$ & 46.436\% & +30.634\% & +30.634\% & -37.846\% & -46.421\% & -46.436\% & - \\ \hline
Value (no IC), $d=14$ & 3.1659 & +16.618\% & - & -51.435\% & -66.252\% & -66.271\% & -66.271\% \\ \hline
Value (IC), $d=14$ & 2.9766 & - & +16.428\% & -49.845\% & -64.110\% & -64.126\% & - \\ \hline\hline
\% Caught, $d=18$ & 46.309\% & +36.809\% & +36.809\% & -37.738\% & -46.297\% & -46.309\% & - \\ \hline
Value (no IC), $d=18$ & 1.6367 & +89.931\% & - & -161.947\% & -205.205\% & -205.291\% & -205.291\% \\ \hline
Value (IC), $d=18$ & 1.4456 & - & +99.910\% & -173.187\% & -219.120\% & -219.120\% & - \\ \hline\hline
\end{tabular}}
\caption{Percent of Unexpected Disruptive Events Caught and Average Value Accumulated with Each Inspection Rule for $P^2(ni)$}
\label{tbl:value_p2}
\end{table}

%\begin{table}[h!]
%\centering \small
%\begin{tabular}{|c|c|c|c|}
%\hline
%$d$ & 14 & 18 & 22 \\ \hline
%Average End Time & 13.3437 & 13.3766 & 13.3585 \\ \hline
%Std. Dev. End Time & 7.5805 & 7.5958 & 7.5854 \\ \hline
%Median End Time & 12 & 12 & 12 \\ \hline
%Min End Time & 2 & 2 & 2 \\ \hline
%Max End Time & 75 & 81 & 71 \\ \hline
%\end{tabular}
%\caption{Descriptive Statistics of Time of Unexpected Disruptive Event for %$P^2(ni)$}
%\label{tbl:endtime_p2}
%\end{table}\normalsize

We again randomly generate the transition probabilities using the same generation procedure as before with mean values $P^2(ni)$. For brevity, these results appear in Appendix \ref{app:acs}. In these experiments, the uncertainty around the transition probabilities is less impactful than for $P^1(ni)$. While the average values accumulated decreased as the standard deviation increases, the patterns in average values accumulated and number of unexpected disruptive events caught by each inspection rule follow the same trends as when the transition matrix is known with certainty. This suggests that for facilities thought to be more prone to quality degradation, the increased frequency in monitoring suggested by the optimal inspection rules is robust to the lack of certainty about the quality degradation of the facility.

The third transition probability matrix we test is more more likely to transition to the state it is currently in than $P^1(ni)$. Table \ref{tbl:inspect_p3} presents the recommended inspection times for each inspection rule with the varying values for $d \in \{18, 22, 26, 30, 34\}$. Note that, with $d=18$, the model where inspection can result in closure determines that it is never optimal to inspect. As such, we exclude $d=14$, which will follow similar trends, and instead include $d=34$.

%\small $$ P^3(ni) = \begin{bmatrix}
%0.935 & 0.0635 & 0 & 0 & 0 & 0 \\
%0 & 0.865 & 0.085 & 0.035 & 0.015 & 0 \\
%0 & 0 & 0.775 & 0.165 & 0.06 & 0 \\
%0 & 0 & 0 & 1 & 0 & 0 \\
%0 & 0 & 0 & 0 & 1 & 0 \\
%0 & 0 & 0 & 0 & 0 & 1 
%\end{bmatrix}. $$\normalsize

Table \ref{tbl:value_p3} presents the results for $P^3(ni)$. The relative differences in the average value accumulated by the optimal inspection rules again follows similar trends as $P^1(ni)$ and $P^2(ni)$ for $d \ge 22$. The most noteworthy difference from the other transition matrices tested is that, when $d=18$ for $P^3(ni)$, it is never optimal to inspect when inspecting may result in the facility needing to close for mandatory maintenance, but still attains an average value accumulated larger than that of most other inspection rules. This suggests the super-linear behavior observed with $P^1(ni)$ and $P^2(ni)$ is actually quadratic. If the penalties for experiencing an unexpected disruptive event are relatively small compared to how quickly we expect the facility to experience an unexpected disruptive event, it is more beneficial to let the facility produce as much as possible. As the penalties increase, the trade-off between the value of letting the facility manufacture their product and inspecting before an unexpected disruptive event occurs becomes more difficult to determine, until the penalties outweigh the value of letting the facility continue operations. At this point, the relative value of inspecting earlier grows drastically as the penalties grow. The rate of this growth in relative value is larger for facilities more likely to experience an unexpected disruptive event. This suggests that the FDA can best allocate inspection resources by dedicating more resources to facilities that produce drugs that are of high importance to public health, primarily focusing on high-risk facilities.%, and Table \ref{tbl:endtime_p3} presents descriptive statistics regarding the time that the Markov chain terminated. The relative differences in the average value accumulated by the optimal inspection rules again follows similar trends as $P^1(ni)$ and $P^2(ni)$ for $d \ge 22$. The most noteworthy difference from the other transition matrices tested is that, when $d=18$ for $P^3(ni)$, it is never optimal to inspect when inspecting may result in the facility needing to close for mandatory maintenance, but still attains an average value accumulated larger than that of most other inspection rules. This suggests the super-linear behavior observed with $P^1(ni)$ and $P^2(ni)$ is actually quadratic. If the penalties for experiencing an unexpected disruptive event are relatively small compared to how quickly we expect the facility to experience an unexpected disruptive event, it is more beneficial to let the facility produce as much as possible. As the penalties increase, the trade-off between the value of letting the facility manufacture their product and inspecting before an unexpected disruptive event occurs becomes more difficult to determine, until the penalties outweigh the value of letting the facility continue operations. At this point, the relative value of inspecting earlier grows drastically as the penalties grow. The rate of this growth in relative value is larger for facilities more likely to experience an unexpected disruptive event. This suggests that the FDA can best allocate inspection resources by dedicating more resources to facilities that produce drugs that are of high importance to public health, primarily focusing on high-risk facilities. 

\begin{table}[h!]
\centering
\resizebox{\textwidth}{!}{\begin{tabular}{|c|c|c|c|c|c|c|c|}
\hline
Rule          & $t_E$ & $t_V$ & $t_{VC}$ & 24 & 60 & 120 & NoIns \\ \hline\hline
\% Caught, $d=18$ & 42.901\% & -42.044\% & -42.901\% & +2.512\% & -38.387\% & -42.829\% & - \\ \hline
Value (no IC), $d=18$ & 10.8196 & +3.490\% & - & -0.644\% & +3.6748\% & +3.506\% & +3.475\% \\ \hline
Value (IC), $d=18$ & 10.7003 & - & +4.629\% & -0.699\% & +4.698\% & +4.657\% & - \\ \hline\hline
\% Caught, $d=22$ & 42.800\% & +2.526\% & +2.526\% & +2.526\% & -38.170\% & -42.706\% & - \\ \hline
Value (no IC), $d=22$ & 9.1597 & +0.128\% & - & +0.128\% & -6.731\% & -8.249\% & -8.273\% \\ \hline
Value (IC), $d=22$ & 9.0403 & - & +0.063\% & +0.063\% & -5.929\% & -7.041\% & - \\ \hline\hline
\% Caught, $d=26$ & 42.665\% & +22.864\% & +22.864\% & +2.595\% & -38.101\% & -42.592\% & - \\ \hline
Value (no IC), $d=26$ & 7.4723 & +5.804\% & - & +1.301\% & -23.728\% & -27.314\% & -27.383\% \\ \hline
Value (IC), $d=26$ & 7.3546 & - & +5.379\% & +1.247\% & -22.713\% & -26.155\% & - \\ \hline\hline
\% Caught, $d=30$ & 42.678\% & +32.480\% & +32.480\% & +2.529\% & -38.033\% & -42.600\% & - \\ \hline
Value (no IC), $d=30$ & 5.8242 & +21.050\% & - & +2.725\% & -48.632\% & -55.469\% & -55.757\% \\ \hline
Value (IC), $d=30$ & 5.7064 & - & +20.689\% & +2.686\% & -47.825\% & -54.443\% & - \\ \hline\hline
\% Caught, $d=34$ & 43.005\% & +38.272\% & +38.272\% & +2.478\% & -38.390\% & -42.915\% & - \\ \hline
Value (no IC), $d=34$ & 4.3305 & +48.112\% & - & +4.863\% & -92.391\% & -100.327\% & -104.545\% \\ \hline
Value (IC), $d=34$ & 4.2093 & - & +48.664\% & +4.894\% & -92.517\% & -104.457\% & - \\ \hline
\end{tabular}}
\caption{Percent of Unexpected Disruptive Events Caught and Average Value Accumulated with Each Inspection Rule for $P^3(ni)$}
\label{tbl:value_p3}
\end{table}

%\begin{table}[h!]
%\centering \small 
%\begin{tabular}{|c|c|c|c|c|c|}
%\hline
%$d$ & 18 & 22 & 26 & 30 & 34 \\ \hline
%Average End Time & 26.5624 & 26.666 & 26.5696 & 26.5911 & 26.6350 \\ \hline
%Std. Dev. End Time & 16.6934 & 16.9076 & 16.7596 & 16.8441 & 16.8656 \\ %\hline
%Median End Time & 23 & 23 & 23 & 23 & 23 \\ \hline
%Min End Time & 3 & 3 & 3 & 3 & 3 \\ \hline
%Max End Time & 182 & 192 & 179 & 178 & 183 \\ \hline
%\end{tabular}
%\caption{Descriptive Statistics of Time of Unexpected Disruptive Event for $P^3(ni)$}
%\label{tbl:endtime_p3}
%\end{table}\normalsize

We again randomly generate transition probabilities using a mean of $P^3(ni)$ to test the robustness of the recommended inspection strategies. For brevity, these results also appear in Appendix \ref{app:acs}. As the base transition matrix is less prone to experiencing an unexpected disruptive event, the uncertainty around the transition probabilities proves to be a challenge in determining when is optimal to inspect. In particular, when $s=0.02$, the value of waiting longer to inspect is often more than that of inspecting earlier for the lower tested values of $d$. This suggests that, if the facility is already expected to maintain higher quality for a longer time, even when we are less certain of the true rate at which the facility's quality will degrade, it will still be acceptable to delay inspecting. However, it is still important to better understand the rate at which facility's quality will degrade, as we may expect delaying inspecting will result in more value we would actually gain.

\subsection{Penalty Parameter Sensitivity Analysis}
\label{ssec:penalty}
\begin{figure}[h!]
\centering
  \includegraphics[width=.7\linewidth]{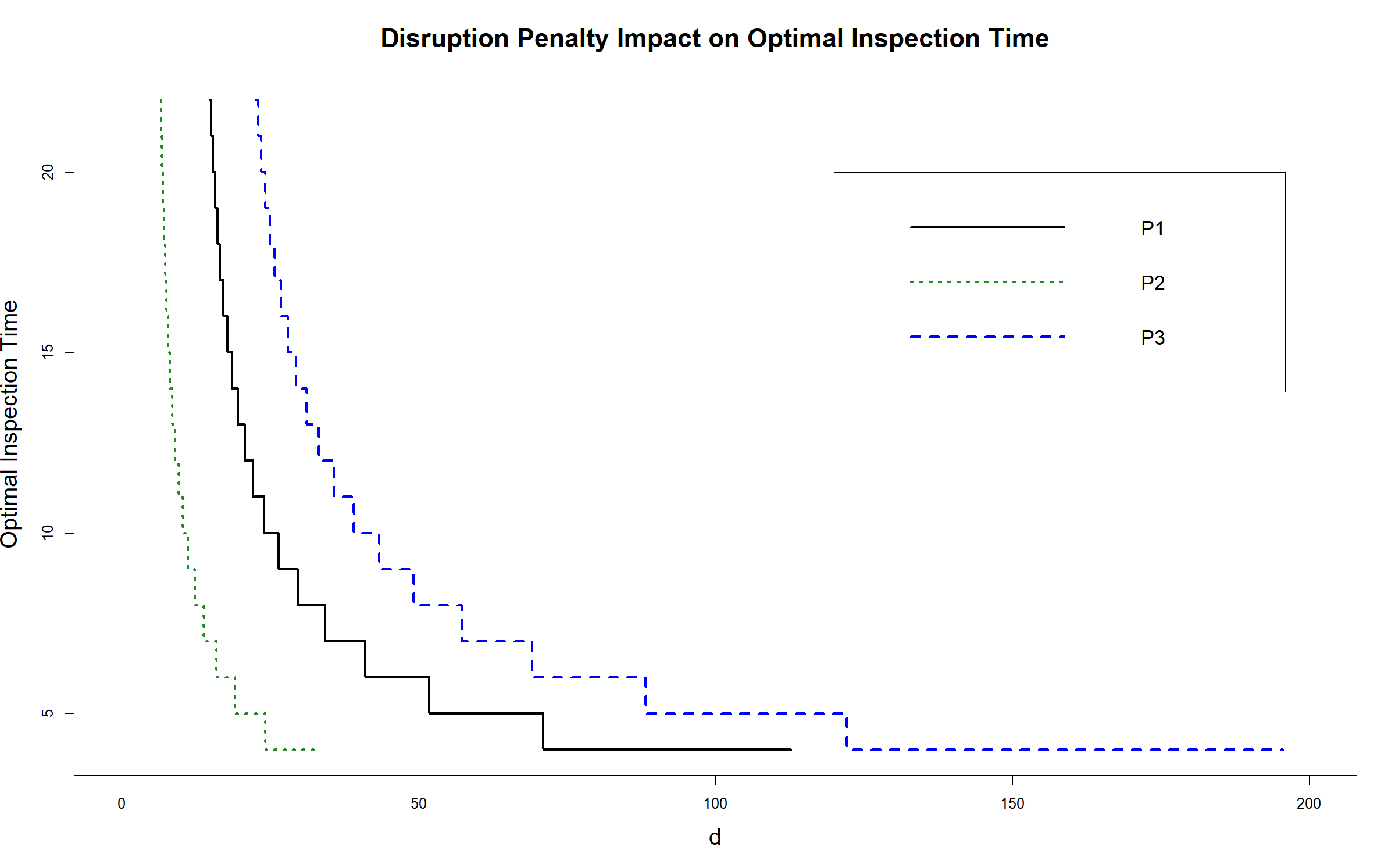}
  \label{fig:p1_opt_t}
\caption{Optimal inspection times with varying $d$}
\label{fig:opt_t}
\end{figure}

We now seek to understand how sensitive our inspection strategy is to changes in the penalty parameters. Given a fixed transition probability matrix, we test how $d$ needs to change such that the recommended inspection time changes. Based on our experiments in Section \ref{ssec:base}, we explored what the inspection time would need to be in order to catch between 50\% and 90\% of disruptive events for each transition matrix tested. We wish to identify the range of values that $d$ can take such that  the inspection periods stay the same, that is, how sensitive the inspection rules are to variations in $d$. Figure \ref{fig:opt_t} presents step plots of the ranges $d$ can lie in  such that each $t \in [4,22]$ is optimal for each transition matrix. The earliest inspection time to catch at least 90\% of unexpected disruptive events is $t=4$ with $P^2(ni)$, while the latest to catch at least 50\% is $t=22$ with $P^3(ni)$. %Tables \ref{tbl:p1_opt_t}, \ref{tbl:p2_opt_t}, and \ref{tbl:p3_opt_t} specify the ranges $[d_{L},d_{U}]$ results in catching each of the specified percentages of unexpected disruptive events.

%\begin{table}[h!]
%\centering \small 
%\begin{tabular}{|c|c|c|c|c|}
%\hline
%\% Caught & $t^*$ & $d_{L}$ & $d_{U}$ & $d_U - d_L$ \\\hline
%50\% & 17 & 16.5882 & 17.1437 & 0.5555\\ \hline
%75\% & 11 & 22.1710 & 24.0267 & 1.8557\\ \hline
%90\% & 7 & 34.2604 & 41.0004 & 6.7400\\ \hline
%\end{tabular}
%\caption{Ranges of $d$ such that certain percentages of unexpected disruptive events are caught for $P^1(ni)$}
%\label{tbl:p1_opt_t}
%\end{table} \normalsize

%\begin{table}[h!]
%\centering \small 
%\begin{tabular}{|c|c|c|c|c|}
%\hline
%\% Caught & $t^*$ & $d_{L}$ & $d_{U}$ & $d_U-d_L$ \\\hline
%50\% & 11 & 9.5578 & 10.2539 & 0.6961 \\ \hline
%75\% & 7 & 13.8172 & 15.9650 & 2.1478 \\ \hline
%90\% & 4 & 24.1638 & 33.0312 & 8.8674 \\ \hline
%\end{tabular}
%\caption{Ranges of $d$ such that certain percentages of unexpected disruptive events are caught for $P^2(ni)$}
%\label{tbl:p2_opt_t}
%\end{table} \normalsize

%\begin{table}[h!]
%\centering \small 
%\begin{tabular}{|c|c|c|c|c|}
%\hline
%\% Caught & $t^*$ & $d_{L}$ & $d_{U}$ & $d_U-d_L$ \\\hline
%50\% & 22 & 22.4775 & 22.9711 & 0.4936 \\ \hline
%75\% & 14 & 29.3898 & 31.0965 & 1.7067 \\ \hline
%90\% & 8 & 49.1557 & 57.2387 & 8.0830 \\ \hline
%\end{tabular}
%\caption{Ranges of $d$ such that certain percentages of unexpected disruptive events are caught for $P^3(ni)$}
%\label{tbl:p3_opt_t}
%\end{table} \normalsize

We observe that, for all three transition matrices, later inspection times are more sensitive to small perturbations of $d$ than earlier inspection times. This seems to be reasonable, as we believe the facility will be more likely to experience an unexpected disruptive event the later we wait to inspect, making changes in $d$ more impactful on later inspection times. Thus, if external factors were to cause a facility to become more critical to the FDA - such as if another facility producing a substitute-able product were to experience an unexpected disruptive events - then it is important for the FDA to inspect the operating facility if it has not been inspected in a long time. Additionally, to be certain that a later inspection time is optimal, it is more critical that the FDA has an accurate estimate of $d$. For earlier inspection times, the FDA does not need to have as accurate an estimation of $d$.

\section{Conclusion and Managerial Insights}
\label{sec:con}
In this work, we propose a novel POMDP approach to support the FDA in deciding when to inspect manufacturing facilities. Our model utilizes the classifications a facility can receive when inspected as the state space, and weighs the trade-offs of the risk of an unexpected disruptive event with the value gained from a facility staying in operation. Under realistic assumptions, we are able to reduce the problem of determining when is the optimal time to inspect to a series of simple computations comparing the value of inspecting now against the value of inspecting after waiting a single time period. Our numerical results show that the optimal inspection rules produced by our models are able to achieve higher average values accumulated over the time between inspections than other rules tested, especially when the penalty for a manufacturing failure occurring is large. We found that our approach is more robust to randomness in the transition probabilities for facilities that are more likely to experience unexpected disruptive events.

Our work has some key insights for FDA decision makers. First, determining when to inspect facilities independently from each other does not require long-term planning; it is sufficient to analyze the trade-offs between inspecting immediately and waiting until the next planning period to inspect. Second, there is a quadratic relationship between the penalty for a manufacturing failure occurring and the relative increase in the average value accumulated when inspection time accounts for the impact a disruptive event has on public health. The rate at which this quadratic grows is impacted by the likelihood of experiencing an unexpected disruptive event. This highlights the importance of dedicating more inspection resources to high-risk facilities that produce drugs that are more critical to public health. It is also important for the FDA to better understand how the quality of facilities degrade over time. Third, uncertainty in the penalties on public health is more impactful on determining the optimal inspection time when the FDA waits longer to inspect. The FDA will need to ensure their estimation of the penalties is accurate to verify that a later inspection time is optimal, whereas the FDA can be less careful about this estimation in situations when the optimal inspection time is shorter.

There are many avenues for future work based off of this initial model. One such avenue is to consider multiple facilities simultaneously. The FDA may be able to better mitigate for drug shortages if they are able to amend their inspection policies to account for the operation statuses for multiple facilities producing substitute-able products. Another avenue for future work is to model how manufacturers react to FDA inspections. If the FDA can more accurately anticipate whether or not manufacturers will comply with any voluntary actions suggested by the FDA, they can better determine how strictly a facility may need to be monitored and how harsh penalties for non-compliance need to be to ensure compliance with the CGMP without forcing the facility to close. Lastly, our theoretical results are a product of the dynamic structure of the value functions and assumptions on the sequences of observation probabilities, neither of which are exclusively tied to the update scheme for belief probabilities in POMDPs. A third avenue of future research is to investigate how inspections times would vary when observation probability sequences that still satisfy our assumptions are generated by different models.

\section*{Acknowledgements}
This material is based upon work partially supported by the National Science Foundation under Grant No. 2228510. 

\bibliographystyle{agsm}
\bibliography{ref}

@article{ahmad2012condition,
    title = {An overview of time-based and condition-based maintenance in industrial application},
    year = {2012},
    journal = {Computers and Industrial Engineering},
    author = {Ahmad, Rosmaini and Kamaruddin, Shahrul},
    number = {1},
    pages = {135--149},
    volume = {63},
    publisher = {Elsevier Ltd}
}

@article{ahuja2021enhancing,
    title = {Enhancing Regulatory Decision Making for Postmarket Drug Safety},
    year = {2021},
    journal = {Management Science},
    author = {Ahuja, Vishal and Alvarez, Carlos A. and Birge, John R. and Syverson, Chad},
    number = {12},
    month = {12},
    pages = {7493--7510},
    volume = {67},
    publisher = {INFORMS Inst.for Operations Res.and the Management Sciences}
}

@article{anand2012decay,
  title={Decay, shock, and renewal: Operational routines and process entropy in the pharmaceutical industry},
  author={Anand, Gopesh and Gray, John and Siemsen, Enno},
  journal={Organization Science},
  volume={23},
  number={6},
  pages={1700--1716},
  year={2012},
  publisher={INFORMS}
}

@article{babich2012managing,
    title = {Managing Opportunistic Supplier Product Adulteration: Deferred Payments, Inspection, and Combined Mechanisms},
    year = {2012},
    journal = {Manufacturing {\&} Service Operations Management},
    author = {Babich, Volodymyr and Tang, Christopher S.},
    number = {2},
    month = {4},
    pages = {301--314},
    volume = {14}
}

@article{ball2017investigator,
    title = {Do plant inspections predict future quality? The role of investigator experience},
    year = {2017},
    journal = {Manufacturing and Service Operations Management},
    author = {Ball, George and Siemsen, Enno and Shahc, Rachna},
    number = {4},
    month = {9},
    pages = {534--550},
    volume = {19},
    publisher = {INFORMS Inst.for Operations Res.and the Management Sciences}
}

@book{bremaud2020markov,
    title = {{Markov Chains}},
    year = {2020},
    author = {Br{\'{e}}maud, Pierre},
    volume = {31},
    publisher = {Springer International Publishing},
    address = {Cham, Switzerland}
}

@book{buckley2013countering,
  title={Countering the problem of falsified and substandard drugs},
  author={Buckley, Gillian J and Gostin, Lawrence O and others},
  year={2013},
  publisher={National Academies Press},
  address = {Washington, D.C., U.S.A.}
}

@article{byon2010wind,
    title = {Optimal maintenance strategies for wind turbine systems under stochastic weather conditions},
    year = {2010},
    journal = {IEEE Transactions on Reliability},
    author = {Byon, Eunshin and Ntaimo, Lewis and Ding, Yu},
    number = {2},
    month = {6},
    pages = {393--404},
    volume = {59}
}

@misc{cder2018manual,
    title = {Understanding {CDER}’s Risk-Based Site Selection Model},
    year = {2018},
    author = {{Center for Drug Evaluation and Research}},
    pages = {1--7},
    url = {https://www.fda.gov/media/116004/download},
    note = {Accessed March 31, 2023}
}

@article{chang2017risk,
    title = {Risk Assessment of Deliberate Contamination of Food Production Facilities},
    year = {2017},
    journal = {IEEE Transactions on Systems, Man, and Cybernetics: Systems},
    author = {Chang, Yanling and Erera, Alan L. and White, Chelsea C.},
    number = {3},
    month = {3},
    pages = {381--393},
    volume = {47},
    publisher = {Institute of Electrical and Electronics Engineers Inc.}
}

@article{cheng2023optimal,
    title = {Optimal joint inspection and mission abort policy for a partially observable system},
    year = {2023},
    journal = {Reliability Engineering and System Safety},
    author = {Cheng, Guoqing and Li, Ling and Shangguan, Chunxia and Yang, Nan and Jiang, Bo and Tao, Ningrong},
    month = {1},
    pages = {1 -- 11},
    volume = {229},
    publisher = {Elsevier Ltd}
}

@article{csardi2006igraph,
    title = {The i{G}raph Software Package for Complex Network Research},
    year = {2006},
    journal = {InterJournal, complex systems},
    author = {Csardi, Gabor and Nepusz, Tamas},
    number = {5},
    pages = {1--9},
    volume = {1695}
}

@misc{edney2019cheap,
    title = {America's Love Affair with Cheap Drugs has a Hidden Cost},
    year = {2019},
    booktitle = {Bloomberg News},
    author = {Edney, Anna},
    month = {1},
    url = {https://www.bloomberg.com/news/features/2019-01-29/america-s-love-affair-with-cheap-drugs-has-a-hidden-cost},
    note = {Accessed March 31, 2023}
}

@incollection{Ergun2022SupplyChains,
    title = {Supply Chain Resilience: Impact of Stakeholder Behavior and Trustworthy Information Sharing with a Case Study on Pharmaceutical Supply Chains},
    year = {2022},
    booktitle = {Tutorials in Operations Research: Emerging and Impactful Topics in Operations},
    author = {Ergun, {\"O}zlem and Griffin, Jacqueline and Chicoine, Noah and Gong, Min and Mohaddesi, Omid and Raziei, Zohreh and Harteveld, Casper and Kaeli, David and Marsella, Stacy},
    month = {10},
    pages = {133--159},
    publisher = {INFORMS},
    address = {Maryland, U.S.A.}
}

@misc{fda2018inspect,
    title = {What does {FDA} inspect?},
    year = {2018},
    author = {{U.S. Food {\&} Drug Administration}},
    month = {3},
    url = {https://www.fda.gov/about-fda/fda-basics/what-does-fda-inspect},
    note = {Accessed March 31, 2023}
}

@misc{fda2019shortages,
    title = {Drug Shortages: Root Causes and Potential Solutions},
    year = {2019},
    author = {{U.S. Food {\&} Drug Administration}},
    url = {https://www.fda.gov/media/131130/download},
    note = {Accessed June 1, 2023}
}

@misc{fda2020notifying,
    title = {Notifying {FDA} of a Permanent Discontinuance or Interruption in Manufacturing Under Section 506C of the {FD\&C} {A}ct},
    year = {2020},
    author = {{U.S. Food {\&} Drug Administration}},
    url = {https://www.fda.gov/media/136486/download},
    note = {Accessed June 1, 2023}
}

@misc{fda2020faq,
    title = {Inspections Database Frequently Asked Questions},
    year = {2020},
    author = {{U.S. Food {\&} Drug Administration}},
    url = {https://www.fda.gov/inspections-compliance-enforcement-and-criminal-investigations/inspection-references/inspections-database-frequently-asked-questions},
    note = {Accessed March 31, 2023}
}

@misc{fda2020report,
    title = {Eighth annual report on drug shortages for calendar year 2020},
    year = {2020},
    author = {{U.S. Food {\&} Drug Administration}},
    url = {https://www.fda.gov/media/150409/download},
    note = {Accessed June 1, 2023}
}

@misc{fda2022cgmp,
    title = {Current Good Manufacturing Practice ({CGMP}) Regulations},
    year = {2023},
    author = {{U.S. Food {\&} Drug Administration}},
    url = {https://www.fda.gov/drugs/pharmaceutical-quality-resources/current-good-manufacturing-practice-cgmp-regulations},
    note = {Accessed May 31, 2023}
}

@article{francas2023drivers,
  title={On the drivers of drug shortages: empirical evidence from Germany},
  author={Francas, David and Mohr, Stephan and Hoberg, Kai},
  journal={International Journal of Operations \& Production Management},
  year={2023},
  pages={1--33},
  publisher={Emerald Publishing Limited}
}

@inproceedings{ge2007optimum,
    title = {Optimum Maintenance Policy with Inspection by Semi-Markov Decision Processes},
    year = {2007},
    booktitle = {2007 39th North American Power Symposium},
    author = {Ge, Haifeng and Tomasevicz, Curtis L. and Asgarpoor, Sohrab},
    month = {9},
    pages = {541--546},
    publisher = {IEEE},
    address = {Las Cruces, NM, U.S.A.}
}

@article{ghandali2020sustainable,
    title = {A POMDP framework to find optimal policy in sustainable maintenance},
    year = {2020},
    journal = {Scientia Iranica},
    author = {Ghandali, R. and Abooie, M. H. and Nezhad, M. S.Fallah},
    number = {3},
    month = {5},
    pages = {1544--1561},
    volume = {27},
    publisher = {Sharif University of Technology}
}

@article{hein2019comparing,
    title = {{Comparing methods for clinical investigator site inspection selection: a comparison of site selection methods of investigators in clinical trials}},
    year = {2019},
    journal = {Journal of Biopharmaceutical Statistics},
    author = {Hein, Nicholas and Rantou, Elena and Schuette, Paul},
    number = {5},
    month = {9},
    pages = {860--873},
    volume = {29},
    publisher = {Taylor and Francis Inc.}
}

@misc{hsgac2023shortsuppy,
    title = {Short Suply: The Health and National Security Risks of Drug Shortages},
    year = {2023},
    author = {{U.S. Senate Committee on Homeland Security \& Governmental Affairs}},
    month = {3},
    url = {https://www.hsgac.senate.gov/wp-content/uploads/2023-06-06-HSGAC-Majority-Draft-Drug-Shortages-Report.-FINAL-CORRECTED.pdf},
    note = {Accessed June 1, 2023}
}

@article{jensen2002fda,
  title={{FDA}’s role in responding to drug shortages},
  author={Jensen, Valerie and Kimzey, Lorene M and Goldberger, Mark J},
  journal={American journal of health-system pharmacy},
  volume={59},
  number={15},
  pages={1423--1425},
  year={2002},
  publisher={Citeseer}
}

@article{jin2021food,
  title={Food safety inspection and the adoption of traceability in aquatic wholesale markets: A game-theoretic model and empirical evidence},
  author={Jin, Cang-yu and Retsef, Levi and Liang, Qiao and Renegar, Nicholas and Zhou, Jie-hong},
  journal={Journal of Integrative Agriculture},
  volume={20},
  number={10},
  pages={2807--2819},
  year={2021},
  publisher={Elsevier}
}

@misc{kimball2023temporary,
    title = {{FDA} allows temporary import of unapproved {C}hinese cancer drug to ease {U.S.} shortage},
    year = {2023},
    booktitle = {CNBC},
    author = {Kimball, Spencer},
    month = {6},
    url = {https://www.cnbc.com/2023/06/02/cancer-drug-shortage-fda-allows-import-of-unapproved-china-chemo-med.html},
    note = {Accessed June 5, 2023}
}

@article{klimberg1992improving,
    title = {Improving the Effectiveness of {FDA} Drug Inspection},
    year = {1992},
    journal = {Operations Research},
    author = {Klimberg, Ronald and Revelle, Charles and Cohon, Jared},
    number = {5},
    month = {10},
    pages = {845--855},
    volume = {40},
    publisher = {Institute for Operations Research and the Management Sciences (INFORMS)}
}

@book{krishnamurthy2016pomdp,
    title = {Partially Observed Markov Decision Processes},
    year = {2016},
    author = {Krishnamurthy, Vikram},
    month = {3},
    publisher = {Cambridge University Press},
    address = {Cambridge, U.K.}
}

@article{levi2019supply,
  title={Supply chain network analytics guiding food regulatory operational policy},
  author={Levi, Retsef and Renegar, Nicholas and Springs, Stacy and Zaman, Tauhid},
  journal={Available at SSRN 3374620},
  year={2019}
}

@article{levitin2019cost,
    title = {Cost effective scheduling of imperfect inspections in systems with hidden failures and rescue possibility},
    year = {2019},
    journal = {Applied Mathematical Modelling},
    author = {Levitin, Gregory and Xing, Liudong and Huang, Hong Zhong},
    month = {4},
    pages = {662--674},
    volume = {68},
    publisher = {Elsevier Inc.}
}

@article{liu2023sequential,
    title = {Selective maintenance and inspection optimization for partially observable systems: An interactively sequential decision framework},
    year = {2023},
    journal = {IISE Transactions},
    author = {Liu, Yu and Gao, Jian and Jiang, Tao and Zeng, Zhiguo},
    number = {5},
    pages = {463--479},
    volume = {55},
    publisher = {Taylor and Francis Ltd.}
}

@techreport{macher2006exploring,
  title={Exploring the information asymmetry gap: Evidence from {FDA} regulation},
  author={Macher, Jeffery T and Mayo, John W and Nickerson, Jack A},
  year={2006},
  url={https://www.researchgate.net/publication/228664514}
}

@article{morato2022dynamic,
    title = {Optimal inspection and maintenance planning for deteriorating structural components through dynamic Bayesian networks and Markov decision processes},
    year = {2022},
    journal = {Structural Safety},
    author = {Morato, P. G. and Papakonstantinou, K. G. and Andriotis, C. P. and Nielsen, J. S. and Rigo, P.},
    month = {1},
    pages = {1--16},
    volume = {94},
    publisher = {Elsevier B.V.}
}

@article{myers2009probability,
    title = {Probability of loss assessment of critical k-Out-of-n:G systems having a mission abort policy},
    year = {2009},
    journal = {IEEE Transactions on Reliability},
    author = {Myers, Albert},
    number = {4},
    pages = {694--701},
    volume = {58}
}

@book{nasem2022building,
    title = {Building Resilience into the Nation's Medical Product Supply Chains},
    author = {{National Academies of Sciences, Engineering, and Medicine}},
    year = {2022},
    editor = {Hopp, Wallace J. and Brown, Lisa and Shore, Carolyn},
    month = {6},
    publisher = {National Academies Press},
    address = {Washington, D.C., U.S.A.}
}

@article{Papakonstantinou2014pomdp,
    title = {Optimum inspection and maintenance policies for corroded structures using partially observable Markov decision processes and stochastic, physically based models},
    year = {2014},
    journal = {Probabilistic Engineering Mechanics},
    author = {Papakonstantinou, K. G. and Shinozuka, M.},
    pages = {93--108},
    volume = {37},
    publisher = {Elsevier Ltd}
}

@article{papakonstantinou2014inspection,
    title = {Planning structural inspection and maintenance policies via dynamic programming and Markov processes. Part I: Theory},
    year = {2014},
    journal = {Reliability Engineering and System Safety},
    author = {Papakonstantinou, K. G. and Shinozuka, M.},
    pages = {202--213},
    volume = {130},
    publisher = {Elsevier Ltd}
}

@article{papakonstantinou2014implementation,
    title = {Planning structural inspection and maintenance policies via dynamic programming and Markov processes. Part II: POMDP implementation},
    year = {2014},
    journal = {Reliability Engineering and System Safety},
    author = {Papakonstantinou, K. G. and Shinozuka, M.},
    pages = {214--224},
    volume = {130},
    publisher = {Elsevier Ltd}
}

@article{pazhayattil2022quantitative,
    title = {A Quantitative Study of {US} {FDA} Inspection Data for Drug Manufacturing Sites},
    year = {2020},
    journal = {Therapeutic Innovation and Regulatory Science},
    author = {Pazhayattil, Ajay Babu and Ingram, Marzena and Sayeed, Naheed},
    number = {4},
    month = {7},
    pages = {725--730},
    volume = {54},
    publisher = {Springer}
}

@misc{pew2012,
    title = {Heparin: A Wake-Up Call on Risks to the U.S. Drug Supply},
    year = {2012},
    author = {{Pew Health Group}},
    url = {https://www.pewtrusts.org/en/research-and-analysis/issue-briefs/2012/05/16/heparin-a-wakeup-call-on-risks-to-the-us-drug-supply},
    note = {Accessed June 1, 2023}
}

@article{qiu2019gamma,
    title = {Gamma process based optimal mission abort policy},
    year = {2019},
    journal = {Reliability Engineering and System Safety},
    author = {Qiu, Qingan and Cui, Lirong},
    month = {10},
    pages = {1--9},
    volume = {190},
    publisher = {Elsevier Ltd}
}

@article{qiu2023optimal,
    title = {Optimal Condition-Based Mission Abort Decisions},
    year = {2023},
    journal = {IEEE Transactions on Reliability},
    author = {Qiu, Qingan and Maillart, Lisa M. and Prokopyev, Oleg A. and Cui, Lirong},
    number = {1},
    month = {3},
    pages = {408--425},
    volume = {72}
}

@article{tucker2020drug,
  title={The drug shortage era: a scoping review of the literature 2001--2019},
  author={Tucker, Emily L and Cao, Yizhou and Fox, Erin R and Sweet, Burgunda V},
  journal={Clinical Pharmacology \& Therapeutics},
  volume={108},
  number={6},
  pages={1150--1155},
  year={2020},
  publisher={Wiley Online Library}
}

@article{wang2023barriers,
    title = {Manufacturing and Regulatory Barriers to Generic Drug Competition: A Structural Model Approach},
    year = {2023},
    journal = {Management Science},
    author = {Wang, Yixin (Iris) and Li, Jun and Anupindi, Ravi},
    number = {3},
    month = {3},
    pages = {1449--1467},
    volume = {69}
}

@article{ventola2011drug,
  title={The drug shortage crisis in the {U}nited {S}tates: causes, impact, and management strategies},
  author={Ventola, C Lee},
  journal={Pharmacy and Therapeutics},
  volume={36},
  number={11},
  pages={740},
  year={2011},
  publisher={MediMedia, USA}
}

@article{white2020generic,
    title = {Generic Drugs Not as Safe as {FDA} Wants You to Believe},
    year = {2020},
    journal = {Annals of Pharmacotherapy},
    author = {White, C. Michael},
    number = {3},
    month = {3},
    pages = {283--286},
    volume = {54}
}

@misc{woodcock2019testimony,
    title = {Securing the {U.S.} Drug Supply Chain: Oversight of {FDA}'s Foreign Inspection Program},
    year = {2019},
    author = {Woodcock, Janet},
    url = {https://www.fda.gov/news-events/congressional-testimony/securing-us-drug-supply-chain-oversight-fdas-foreign-inspection-program-12102019},
    note = {Accessed March 31, 2023}
}

@article{woodcock2020quality,
  title={Quality: The Often Overlooked Critical Element For Assuring Access To Safe And Effective Drugs},
  author={Woodcock, Janet and Kopcha, Michael},
  journal={Health Affairs Forefront},
  year={2020}
}

@article{wu2020monitoring,
    title = {The More Monitoring, the Better Quality? {E}mpirical Evidence from the Generic Drug Industry},
    year = {2020},
    author = {Wu, Anqi and Wang, Yixin (Iris)},
    journal={Available at SSRN 3596559}
}

@article{xu2022constrained,
    title = {A risk-aware maintenance model based on a constrained Markov decision process},
    year = {2022},
    journal = {IISE Transactions},
    author = {Xu, Jianyu and Zhao, Xiujie and Liu, Bin},
    number = {11},
    pages = {1072--1083},
    volume = {54},
    publisher = {Taylor and Francis Ltd.}
}

@article{zhang2017continuous,
    title = {Continuous-Observation Partially Observable Semi-Markov Decision Processes for Machine Maintenance},
    year = {2017},
    journal = {IEEE Transactions on Reliability},
    author = {Zhang, Mimi and Revie, Matthew},
    number = {1},
    pages = {202--218},
    volume = {66}
}

@article{zhang2019unannounced,
    title = {The Effect of Unannounced Inspection on Prevention of Drug Fraud},
    year = {2019},
    journal = {Journal of Systems Science and Systems Engineering},
    author = {Zhang, Manman and Zhang, Juliang and Cheng, T. C.E. and Hua, Guowei and Yan, Xiaojie and Liu, Yi},
    number = {1},
    month = {2},
    pages = {63--90},
    volume = {28},
    publisher = {Springer Verlag}
}

@article{zhao2021optimal,
    title = {Optimal inspection and mission abort policies for systems subject to degradation},
    year = {2021},
    journal = {European Journal of Operational Research},
    author = {Zhao, Xian and Sun, Jinglei and Qiu, Qingan and Chen, Ke},
    number = {2},
    month = {7},
    pages = {610--621},
    volume = {292},
    publisher = {Elsevier B.V.}
}

@article{zhou2022effects,
    title = {Effects of regulatory policy mixes on traceability adoption in wholesale markets: Food safety inspection and information disclosure},
    year = {2022},
    journal = {Food Policy},
    author = {Zhou, Jiehong and Jin, Yu and Liang, Qiao},
    month = {2},
    pages = {1--11},
    volume = {107},
    publisher = {Elsevier Ltd}
}

@article{zhu2021cbm,
    title = {Condition-based maintenance for multi-component systems: Modeling, structural properties, and algorithms},
    year = {2021},
    journal = {IISE Transactions},
    author = {Zhu, Zhicheng and Xiang, Yisha},
    number = {1},
    month = {1},
    pages = {88--100},
    volume = {53},
    publisher = {Taylor and Francis Ltd.}
}

@inproceedings{Zouch2011road,
    title = {Optimal resurfacing decisions for road maintenance: A POMDP perspective},
    year = {2011},
    booktitle = {2011 Proceedings-Annual Reliability and Maintainability Symposium},
    author = {Zouch, Mariem and Yeung, Thomas G. and Castanier, Bruno},
    pages = {1--6},
    publisher = {IEEE},
    address = {Lake Buena Vista, FL, U.S.A.}
}

\newpage
\appendix
\section{Summary of Notation}
\label{app:notation}
\begin{table}[H]
    \centering
    {\small %
    \begin{tabular}{|c|l|}
        \hline
        Parameter & Description of Parameter \\\hline
        $T$ & Time Horizon \\ \hline
        $d$ & Penalty for being in the ``Experienced a Manufacturing Failure" state $D$ \\ \hline
        $c$ & Penalty for being in the ``Closed For Non-mandatory Maintenance" state $C$ \\ \hline
        $\Tilde{c}$ & Penalty for the manufacturer closing for mandatory maintenance after an inspection \\ \hline
        $\alpha_S$ & Factor by which $\Tilde{c}$ is smaller than the penalty of being in state $S \in \{D,C\}$ \\ \hline
        $P$ & Transition probability matrix \\ \hline
        $p_{SS'}$ & Probability of transitioning from state $S$ to state $S'$ \\ \hline
        $p_{SIC}$ & \begin{tabular}[l]{@{}l@{}}Probability of the manufacturer closing after an inspection occurs\\ when inspected in state $S$ \end{tabular}\\ \hline
        $b^t_S$ & Belief probability that we are in state $S$ at time $t$ \\ \hline
        $P[(dr)]^t$ & Probability of experiencing a disrupting based on belief probabilities $b^t$ \\ \hline
        $P[(cr)]^t$ & \begin{tabular}[l]{@{}l@{}}Probability of experiencing a closure without intervention\\ based on belief probabilities $b^t$ \end{tabular}\\ \hline
        $P[(nr)]^t$ & \begin{tabular}[l]{@{}l@{}}Probability of no issues being reported by the manufacturer\\ based on belief probabilities $b^t$ \end{tabular}\\ \hline
        $P[(cr)|(i)]^t$ & Probability of the facility closing due to inspection based on belief probabilities $b^t$\\ \hline
        $\sigma^t_j$ & Conditional plan starting at $t$ and inspecting during time period $t+j$ \\ \hline
        $V^t(b, \sigma)$ & \begin{tabular}[l]{@{}l@{}}Cumulative expected value accrued from time period $t$ to $T$,\\ starting with belief probability vector $b$ and implementing conditional plan $\sigma$\end{tabular}
        \\ \hline
        $\Theta^i_{SS'}$ & Set of all paths from state $S$ to state $S'$ in $i-1$ time periods \\ \hline
        $\theta^i$ & Path from state $\theta_1^i$ to state $\theta^i_i$ \\ \hline
        $k_S$ & Expected reward after waiting one time period, starting from state $S$\\ \hline
        $k_{SS'}^t$ & \begin{tabular}[l]{@{}l@{}}Expected reward after waiting $t-1$ time periods, \\ starting from state $S$ and ending at state $S'$ at time $T-1$ \end{tabular}\\ \hline
        $f_{\theta^i}$ & Recursive parameter based on transition probabilities following the path $\theta^i$ \\ \hline
        $\Tilde{c}_{SS'IC}^i$ & \begin{tabular}[l]{@{}l@{}}Expected penalty from closing after an inspection\\ when inspecting after $i-1$ time periods, starting at state $S$ and inspecting in state $S'$\end{tabular}\\ \hline
    \end{tabular}
    }%
    \caption{Description of model parameters}
    \label{tab:noteset1}
\end{table}

\begin{table}[H]
    \centering
    {\small %
    \begin{tabular}{|c|l|}
        \hline
        Parameter & Description of Parameter \\\hline
        $t^*$ & optimal inspection time \\ \hline
        $\delta$ & \begin{tabular}[l]{@{}l@{}} Minimum increase in $d$ such that inspecting one time period earlier \\ than $t^*$ becomes optimal \end{tabular}\\ \hline
        $\gamma$ & \begin{tabular}[l]{@{}l@{}} Minimum increase in $c$ such that inspecting one time period earlier \\ than $t^*$ becomes optimal \end{tabular}\\ \hline
    \end{tabular}
    }%
    \caption{Description of sensitivity analysis parameters}
    \label{tab:noteset2}
\end{table}

\begin{table}[H]
    \centering
    {\small %
    \begin{tabular}{|c|l|}
        \hline
        Parameter & Description of Parameter \\\hline
        $\mu_S$ & Expected time until transitioning into an absorbing state starting at state $S$\\ \hline
        $t_V$ & Inspection time recommended by base POMDP model\\ \hline
        $t_{VC}$ & Inspection time recommended by POMDP model with closure due to inspection\\ \hline
        $t_{E}$ & Inspection time determined by the expected time until reaching an absorbing state\\ \hline
        $\Tilde{P}$ & Randomly generated transition probability matrix using normal distribution \\ \hline
        $\Tilde{p}_{ij}$ & Randomly generated probability of transition from state $i$ to state $j$ \\ \hline
        $s$ & standard deviation\\ \hline
        $d_L$ & minimum value for $d$ where inspecting in a given time period is optimal\\ \hline
        $d_U$ & maximum value for $d$ where inspecting in a given time period is optimal\\ \hline
    \end{tabular}
    }%
    \caption{Description of sensitivity analysis parameters}
    \label{tab:noteset3}
\end{table}

\section{Proofs}
\label{app:pfvar}
\textbf{Proof of Theorem \ref{thrm:equiv}:}
\begin{proof}
    We prove Equations \eqref{eq:rec} and \eqref{eq:nonrec} are equivalent by induction. First, consider $t = T-2$. For Equation \eqref{eq:rec}, we have
    \small\begin{align*}
        V^{T-2}(b^{T-2}, \sigma_{2}^{T-2}) &= 1 - P[(dr)]^{T-1} d - P[(cr)]^{T-1} c + P[(nr)]^{T-1} V^{T-1}(b^{T-1},\sigma_{1}^{T-1}) \\
        &= 1 - P[(dr)]^{T-1} d - P[(cr)]^{T-1} c + P[(nr)]^{T-1} (1 + k_N b^{T-1}_N + k_V b^{T-1}_V + k_O b^{T-1}_O) \\
        &= 1 - P[(dr)]^{T-1} d - P[(cr)]^{T-1} c + P[(nr)]^{T-1}  \\
        &+ k_N(p_{NN} b_N^{T-2}) + k_V(p_{NV} b_N^{T-2} + p_{VV} b_V^{T-2}) + k_O(p_{NO} b_N^{T-2} + p_{VO} b_V^{T-2} + p_{OO} b_O^{T-2}) \\
        &= 1 + (k_N b_N^{T-2} + k_V b_V^{T-2} + k_O b_O^{T-2}) \\
        &+ (k_O p_{NO} + k_V p_{NV} + k_N p_{NN})b_N^{T-2} + (k_O p_{VO}^{T-2} + k_V p_{VV})b_V^{T-2} + k_O p_{OO} b_O^{T-2} \\
        &= 1 + (k_O p_{NO} + k_V p_{NV} + k_N (p_{NN}+1))b_N^{T-2}\\
        &+ (k_O p_{VO} + k_V (p_{VV}+1))b_V^{T-2} + k_O (p_{OO}+1) b_O^{T-2}.
    \end{align*}\normalsize

    For Equation \eqref{eq:nonrec}, we first note that $\Theta^2 = \{(N,N), (N,V), (N,O), (V,V), (V,O), (O,O)\}$. By Equation \eqref{eq:nonrec}, we have
    \begin{align*}
        V^{T-2}(b^{T-2}, \sigma_{2}^{T-2}) &= 1 + k_N (p_{NN}+1) b_N^{T-2} + k_V p_{NV} b_N^{T-2} + k_O p_{NO} b_N^{T-2} \\
        &+ k_V (p_{VV}+1) b_V^{T-2} + k_O p_{VO} b_V^{T-2} + k_O (p_{OO}+1) b_O^{T-2} \\
        &= 1 + (k_O p_{NO} + k_V p_{NV} + k_N (p_{NN}+1))b_N^{T-2} \\
        &+ (k_O p_{VO} + k_V (p_{VV}+1))b_V^{T-2} + k_O (p_{OO}+1) b_O^{T-2}.
    \end{align*}

    Thus, for $t=T-2$, Equations \eqref{eq:rec} and \eqref{eq:nonrec} are equivalent. Next, suppose that Equations \eqref{eq:rec} and \eqref{eq:nonrec} are equivalent for $t=T-j, \ldots, T-2$. We now show they are also equivalent for $t=T-(j+1)$.

    First, note that for any $\theta^{j+1} \in \Theta^{j+1}_{SS'}$, there is a $\bar{\theta}^j \in \Theta^j_{\theta^{j+1}_2 S'}$ such that $\theta^{j+1} = (\theta^{j+1}_1,\bar{\theta}^j)$. Thus,
    \begin{align*}
        & V^{T-(j+1)}(b^{T-(j+1)}, \sigma_{j+1}^{T-(j+1)}) \\
        &= 1 - P[(dr)]^{T-j} d - P[(cr)]^{T-j} c + P[(nr)]^{T-j} V^{T-j}(b^{T-j},\sigma_{j}^{T-j}) \\
        &= 1 - P[(dr)]^{T-j} d - P[(cr)]^{T-j} c + P[(nr)]^{T-j} (1 + \sum_{S \in \{N,V,O\}} \sum_{\substack{S' \in \{N,V,O\} \\S' \text{ is reachable from }S}} k_{SS'}^j b^{T-j}_S) \\
        &= 1 + k_N b_N^{T-(j+1)} + k_V b_V^{T-(j+1)} + k_O b_O^{T-(j+1)}  + \sum_{S' \in \{N,V,O\}} k_{NS'}^j (p_{NN} b_N^{T-(j+1)}) \\
        & + \sum_{S \in \{V,O\}} k_{VS'}^j (p_{NV} b_N^{T-(j+1)} + p_{VV} b_V^{T-(j+1)}) + k_{OO}^j (p_{NO} b_N^{T-(j+1)} + p_{VO} b_V^{T-(j+1)} + p_{OO} b_O^{T-(j+1)}) \\
        &= 1 + (k_N + \sum_{S' \in \{N,V,O\}} k_{NS'}^j p_{NN} + \sum_{S \in \{V,O\}} k_{VS'}^j p_{NV} + k_{OO}^j p_{NO}) b_N^{T-(j+1)} \\
        & + (k_V + \sum_{S \in \{V,O\}} k_{VS'}^j p_{VV} + k_{OO}^j p_{VO}) b_V^{T-(j+1)} + (k_O + k_{OO}^j p_{OO}) b_O^{T-(j+1)}.
    \end{align*}

    For $S, S' \in \{N,V,O\}$ where $S \ne S'$, we have that 
    \begin{align*}
        k_{SS'}^{j+1} &= k_{S'} \sum_{\theta^{j+1} \in \Theta^{j+1}_{S S'}} f_{\theta^{j+1}} \\
        &= k_{S'} \sum_{\substack{\theta^{j+1}_2 \in \{N,V,O\} \\ \theta^{j+1}_2 \text{ is reachable from } S}} p_{S \theta^{j+1}_2} \left(\sum_{\bar{\theta}^j \in \Theta^j_{\theta^{j+1}_2 S'}} f_{\bar{\theta}^j} \right) \\
        &= \sum_{\substack{\theta^{j+1}_2 \in \{N,V,O\} \\ \theta^{j+1}_2 \text{ is reachable from } S}} p_{S \theta^{j+1}_2} \left(k_S'\sum_{\bar{\theta}^j \in \Theta^j_{\theta^{j+1}_2 S'}} f_{\bar{\theta}^j} \right) \\
        &= \sum_{\substack{\theta^{j+1}_2 \in \{N,V,O\} \\ \theta^{j+1}_2 \text{ is reachable from } S}} k_{\bar{\theta}^j_2 S'} p_{S \theta^{j+1}_2}.
    \end{align*}

    Likewise, when $S = S'$, we have that
    \begin{align*}
        k_{SS}^{j+1} &= k_S f_{(S, S, \ldots, S)}^{j+1} = k_S (p_{SS} f_{(S, \ldots, S)}^j +1) = k_S + k_{SS}^j p_{SS}
    \end{align*}

    Thus, 
    \begin{align*}
        &1 + (k_N + \sum_{S' \in \{N,V,O\}} k_{NS'}^i p_{NN} + \sum_{S \in \{V,O\}} k_{VS'}^i p_{NV} + k_{OO}^i p_{NO}) b_N^{T-(j+1)} \\
        & + (k_V + \sum_{S \in \{V,O\}} k_{VS'}^j p_{VV} + k_{OO}^j p_{VO}) b_V^{T-(j+1)} + (k_O + k_{OO}^j p_{OO}) b_O^{T-(j+1)} \\
        &= 1 + \sum_{S \in \{N,V,O\}} \sum_{\substack{S' \in \{N,V,O\} \\S' \text{ is reachable from }S}} k_{SS'}^{j+1} b_S^{T-(j+1)}.
    \end{align*}

    Thus, for $t=T-(j+1)$, Equations \eqref{eq:rec} and \eqref{eq:nonrec} are equivalent.
\end{proof}

\noindent\textbf{Proof of Lemma \ref{lm:timeinc}:}
\begin{proof}
    Suppose that $(P[(dr)]^t)$ and $(P[(cr)]^t)$ are non-decreasing sequences and that $(P[(nr)]^t)$ is a non-increasing sequence. We prove that if $V^t(b^t, \sigma_j^t) \ge V^t(b^t, \sigma_{j+1}^t)$ implies $V^{t+1}(b^{t+1}, \sigma_j^{t+1}) \ge V^{t+1}(b^{t+1}, \sigma_{j+1}^{t+1})$ by induction. First consider $j=0$. Suppose that $V^t(b^t, (i)) \ge V^t(b^t, (ni,i))$, and recall that $V^t(b^t, (i))=1$ for all $t$. We have that
    \begin{align*}
        V^{t+1}(b^{t+1},(i)) = 1 &\ge V^{t}(b^{t},(ni,i)) \\
        & = 1 - P[(dr)]^{t+1} d - P[(cr)]^{t+1} c + P[(nr)]^{t+1} \\
        & \ge 1 - P[(dr)]^{t+2} d - P[(cr)]^{t+2} c + P[(nr)]^{t+2}\\
        & = V^{t+1}(b^{t+1},(ni,i))
    \end{align*}

    Next, suppose $V^{t}(b^{t}, \sigma_{j}^{t}) \ge V^{t}(b^{t},\sigma_{j+1}^{t})$ implies $V^{t+1}(b^{t+1}, \sigma_{j}^{t+1}) \ge V^{t+1}(b^{t+1},\sigma_{j+1}^{t+1})$ for $0 \le j \le l < T-t$. We now show that if $V^{t}(b^{t}, \sigma_{j}^{t}) \ge V^{t}(b^{t},\sigma_{j+1}^{t})$, then $V^{t+1}(b^{t+1}, \sigma_{j}^{t+1}) \ge V^{t+1}(b^{t+1},\sigma_{j+1}^{t+1})$ for $j = l+1$. Suppose $V^{t}(b^{t}, \sigma_{l+1}^{t}) \ge V^{t}(b^{t},\sigma_{l+2}^{t})$, and note that 
    \begin{align*}
        & V^{t}(b^{t}, \sigma_{l+1}^{t}) \ge V^{t}(b^{t},\sigma_{l+2}^{t})\\
        & 1 - P[(dr)]^{t+1} d - P[(cr)]^{t+1} c + P[(nr)]^{t+1}(V^{t+1}(b^{t+1},\sigma_l^{t+1})) \\
        &\ge 1 - P[(dr)]^{t+1} d - P[(cr)]^{t+1} c + P[(nr)]^{t+1}(V^{t+1}(b^{t+1},\sigma_{l+1}^{t+1})) \\
        & V^{t+1}(b^{t+1},\sigma_l^{t+1})) \ge V^{t+1}(b^{t+1},\sigma_{l+1}^{t+1}).
    \end{align*}

    Since $V^{t+1}(b^{t+1},\sigma_l^{t+1})) \ge V^{t+1}(b^{t+1},\sigma_{l+1}^{t+1})$, then $V^{t+2}(b^{t+2},\sigma_l^{t+2})) \ge V^{t+2}(b^{t+2},\sigma_{l+1}^{t+2})$ by our inductive hypothesis. Thus, 
    \begin{align*}
        V^{t+1}(b^{t+1}, \sigma_{l+1}^{t+1})  & = 1 - P[(dr)]^{t+2}d - P[(cr)]^{t+2}c + P[(nr)]^{t+2} (V^{t+2}(b^{t+2},\sigma_l^{t+2})) \\
        & \ge 1 - P[(dr)]^{t+2}d - P[(cr)]^{t+2}c + P[(nr)]^{t+2} (V^{t+2}(b^{t+2},\sigma_{l+1}^{t+2})) \\
        & = V^{t+1}(b^{t+1}, \sigma_{l+2}^{t+1}).
    \end{align*}

    We thus have that if $V^{t}(b^{t}, \sigma_{j}^{t}) \ge V^{t}(b^{t},\sigma_{j+1}^{t})$ implies that $V^{t+1}(b^{t+1}, \sigma_{j}^{t+1}) \ge V^{t+1}(b^{t+1},\sigma_{j+1}^{t+1})$, then $V^{t}(b^{t}, \sigma_{j+1}^{t}) \ge V^{t}(b^{t},\sigma_{j+2}^{t})$ implies that $V^{t+1}(b^{t+1}, \sigma_{j+1}^{t+1}) \ge V^{t+1}(b^{t+1},\sigma_{j+2}^{t+1})$.

    Similarly, we prove that if $V^t(b^t, \sigma_j^t) \le V^t(b^t, \sigma_{j+1}^t)$, then $V^{t-1}(b^{t-1}, \sigma_j^{t-1}) \le V^{t-1}(b^{t-1}, \sigma_{j+1}^t)$ by induction on $j$. First consider $j=0$. Suppose that $V^t(b^t, (i)) \le V^t(b^t, (ni,i))$. 
    We have that
    \begin{align*}
        V^{t-1}(b^{t+1},(i)) = 1 &\le V^{t}(b^{t},(ni,i)) \\
        & = 1 - P[(dr)]^{t+1} d - P[(cr)]^{t+1} c + P[(nr)]^{t+1} \\
        & \le 1 - P[(dr)]^{t} d - P[(cr)]^{t} c + P[(nr)]^{t}\\
        & = V^{t-1}(b^{t-1},(ni,i))
    \end{align*}
    
    Next, suppose that $V^{t}(b^{t}, \sigma_{j}^{t}) \le V^{t}(b^{t},\sigma_{j+1}^{t})$ implies $V^{t-1}(b^{t-1}, \sigma_{j}^{t-1}) \ge V^{t-1}(b^{t-1},\sigma_{j+1}^{t-1})$ for $0 \le j \le l < T-t$. We know show $V^{t}(b^{t}, \sigma_{j}^{t}) \le V^{t}(b^{t},\sigma_{j+1}^{t})$ implies $V^{t-1}(b^{t-1}, \sigma_{j}^{t-1}) \ge V^{t-1}(b^{t-1},\sigma_{j+1}^{t-1})$ for $j = l+1$. Suppose $V^{t}(b^{t}, \sigma_{l+1}^{t}) \le V^{t}(b^{t},\sigma_{l+2}^{t})$, and note that 
    \begin{align*}
        & V^{t}(b^{t}, \sigma_{l+1}^{t}) \le V^{t}(b^{t},\sigma_{l+2}^{t})\\
         & 1 - P[(dr)]^{t+1} d - P[(cr)]^{t+1} c + P[(nr)]^{t+1}(V^{t+1}(b^{t+1},\sigma_l^{t+1})) \\
         &\le 1 - P[(dr)]^{t+1} d - P[(cr)]^{t+1} c + P[(nr)]^{t+1}(V^{t+1}(b^{t+1},\sigma_{l+1}^{t+1})) \\
         &V^t(b^{t+1},\sigma_l^{t+1})) \le V^{t+1}(b^{t+1},\sigma_{l+1}^{t+1}).
    \end{align*}

    Since $V^{t+1}(b^{t+1},\sigma_l^{t+1})) \le V^{t+1}(b^{t+1},\sigma_{l+1}^{t+1})$, then $V^t(b^{t},\sigma_l^{t})) \le V^{t}(b^{t},\sigma_{l+1}^{t})$ by our inductive hypothesis. Thus, 
    \begin{align*}
        V^{t-1}(b^{t-1}, \sigma_{l+1}^{t-1})  & = 1 - P[(dr)]^{t}d - P[(cr)]^{t}c + P[(nr)]^{t} (V^{t}(b^{t},\sigma_l^{t})) \\
        & \le 1 - P[(dr)]^{t}d - P[(cr)]^{t}c + P[(nr)]^{t} (V^{t}(b^{t},\sigma_{l+1}^{t})) \\
        & = V^{t-1}(b^{t-1}, \sigma_{l+2}^{t-1}).
    \end{align*}

    We thus have that if $V^t(b^t, \sigma_j^t) \le V^t(b^t, \sigma_{j+1}^t)$, then $V^{t-1}(b^{t-1}, \sigma_j^{t-1}) \le V^{t-1}(b^{t-1}, \sigma_{j+1}^t)$, then $V^t(b^t, \sigma_{j+1}^t) \le V^t(b^t, \sigma_{j+2}^t)$, then $V^{t-1}(b^{t-1}, \sigma_{j+1}^{t-1}) \le V^{t-1}(b^{t-1}, \sigma_{j+2}^t)$.
\end{proof}

\noindent\textbf{Proof of Theorem \ref{thrm:nii}:}
\begin{proof}
    Suppose $(P[(dr)]^t)$ and $(P[(cr)]^t)$ are non-decreasing sequences, $(P[(nr)]^t)$ is a non-increasing sequence and there exists a $t$ such that $V^t(b^t,(i)) \ge V^t(b^t, (ni,i))$, and recall that $V^t(b^t,(i))=1$. We show that $V^t(b^t,\sigma_{j}^t) \ge V^t(b^t,\sigma_{j+1}^t)$ for $1 \le j < T-t$ be induction. First consider $j=1$. By Lemma \ref{lm:timeinc}, we have that $V^{t+1}(b^{t+1},(i)) \ge V^{t+1}(b^{t+1},(ni,i))$. We have that
    \begin{align*}
        V^t(b^t,(ni,ni,i)) & = 1 - P[(dr)]^{t+1} d - P[(cr)]^{t+1} c + P[(nr)]^{t+1} (V^{t+1}(b^{t+1},(ni,i))) \\
        & \le 1 - P[(dr)]^{t+1} d - P[(cr)]^{t+1} c + P[(nr)]^{t+1} (V^{t+1}(b^{t+1},(i))) \\
        & \le 1 - P[(dr)]^{t+1} d - P[(cr)]^{t+1} c + P[(nr)]^{t+1} \\
        & = V^t(b^t, (ni,i))
    \end{align*}
    Thus, $V^t(b^t, (ni,i)) \ge V^t(b^t,(ni,ni,i))$.

   Now suppose $V^t(b^t,\sigma_{j}^t) \ge V^t(b^t,\sigma_{j+1}^t)$ for $1 \le j \le l < T-t$. We now show $V^t(b^t,\sigma_{j}^t) \ge V^t(b^t,\sigma_{j+1}^t)$ for $j=l+1$. Suppose $V^t(b^t,\sigma_{l}^t) \ge V^t(b^t,\sigma_{l+1}^t)$. We know that $V^{t+1}(b^{t+1},\sigma_{l}^{t+1}) \ge V^{t+1}(b^{t+1},\sigma_{l+1}^{t+1})$ by Lemma \ref{lm:timeinc}. We have that 
    \begin{align*}
        V^t(b^t,\sigma_{l+1}^t) &= 1 - P[(dr)]^{t+1} d - P[(cr)]^{t+1} c + P[(nr)]^{t+1} (V^{t+1}(b^{t+1}, \sigma_{l}^{t+1})) \\
        & \ge 1 - P[(dr)]^{t+1} d - P[(cr)]^{t+1} c + P[(nr)]^{t+1} (V^{t+1}(b^{t+1}, \sigma_{l+1}^{t+1})) \\
        & = V^t(b^t,\sigma_{l+2}^t)
    \end{align*}

    Thus, if $V^t(b^t,\sigma_{l}^t) \ge V^t(b^t,\sigma_{l+1}^t)$, then $V^t(b^t,\sigma_{l+1}^t) \ge V^t(b^t,\sigma_{l+2}^t)$.
\end{proof}

\noindent\textbf{Proof of Theorem \ref{eq:nonrecAlt}:}
\begin{proof}
    We prove Equations \eqref{eq:rec} and \eqref{eq:nonrecAlt} are equivalent by induction. First, consider $t=T-2$. For Equation \ref{eq:rec}, 
    \small\begin{align*}
        V^{T-2}(b^{T-2}, \sigma_{2}^{T-2}) &= 1 - P[(dr)]^{T-1} d - P[(cr)]^{T-1} c + P[(nr)]^{T-1} V^{T-1}(b^{T-1},\sigma_{1}^{T-1}) \\
        &= 1 - P[(dr)]^{T-1} d - P[(cr)]^{T-1} c \\
        &+ P[(nr)]^{T-1} (1 + (k_N - \Tilde{c}(p_{NN}p_{NIC} + p_{NV}p_{VIC} + p_{NO}p_{OIC}))b_N^{T-1} \\
        &+ (k_V - \Tilde{c}(p_{VV}p_{VIC}+p_{VO}p_{OIC}))b_V^{T-1} + (k_O -\Tilde{c}p_{OO}p_{OIC})b_O^{T-1}) \\
        &= 1 - P[(dr)]^{T-1} d - P[(cr)]^{T-1} c + P[(nr)]^{T-1} \\
        &+ (k_N - \Tilde{c}(p_{NN}p_{NIC} + p_{NV}p_{VIC} + p_{NO}p_{OIC}))(p_{NN}b_N^{T-2}) \\
        &+ (k_V - \Tilde{c}(p_{VV}p_{VIC}+p_{VO}p_{OIC}))(p_{NV} b_N^{T-2} + p_{VV} b_V^{T-2})\\ 
        &+(k_O -\Tilde{c}p_{OO}p_{OIC})(p_{NO} b_N^{T-2} + p_{VO} b_V^{T-2} + p_{OO} b_O^{T-2}) \\
        &= 1 + (k_N b_N^{T-2} + k_V b_V^{T-2} + k_O b_O^{T-2}) \\
        &+ (k_O p_{NO} + k_V p_{NV} + k_N p_{NN} - \Tilde{c}(p_{NO} p_{OO} p_{OIC} + p_{NV} p_{VO} p_{OIC} + p_{NV} p_{VV} p_{VIC} \\ 
        &+ p_{NN} p_{NO} p_{OIC} + p_{NN} p_{NV} p_{VIC} + p_{NN} p_{NN} p_{NIC}))b_N^{T-2} \\
        & + (k_O p_{VO} + k_V p_{VV} -\Tilde{c}(p_{VO}p_{OO}p_{OIC} + p_{VV} p_{VO} p_{OIC} + p_{VV} p_{VV} p_{VIC}))b_V^{T-2} \\
        &+ (k_O p_{OO} - \Tilde{c}p_{OO}p_{OO} p_{OIC}) b_O^{T-2} \\
        &= 1 + (k_O p_{NO} + k_V p_{NV} + k_N (p_{NN}+1) - \Tilde{c}(p_{NO} p_{OO} p_{OIC} + p_{NV} p_{VO} p_{OIC} \\
        &+ p_{NV} p_{VV} p_{VIC} + p_{NN} p_{NO} p_{OIC} + p_{NN} p_{NV} p_{VIC} + p_{NN} p_{NN} p_{NIC}))b_N^{T-2}\\
        &+ (k_O p_{VO} + k_V (p_{VV}+1) -\Tilde{c}(p_{VO}p_{OO}p_{OIC} + p_{VV} p_{VO} p_{OIC} + p_{VV} p_{VV} p_{VIC}))b_V^{T-2} \\
        &+ (k_O (p_{OO}+1) - \Tilde{c}p_{OO}p_{OO} p_{OIC})b_O^{T-2}.
    \end{align*}\normalsize

    For Equation \eqref{eq:nonrecAlt}, we first note that $\Theta^2 = \{(N,N), (N,V), (N,O), (V,V), (V,O), (O,O)\}$ and 
    \begin{align*}
        \Theta^3 = \{\,
            & (N,N,N), (N,N,V), (N,N,O), (N,V,V), (N,V,O),\\
            & (N,O,O), (V,V,V), (V,V,O), (V,O,O), (O,O,O)\,\}\enspace.
    \end{align*}
    By Equation \eqref{eq:nonrecAlt}, we have 
    \begin{align*}
        V^{T-2}(b^{T-2}, \sigma_{2}^{T-2}) &= 1 + k_N (p_{NN}+1) b_N^{T-2} + k_V p_{NV} b_N^{T-2} + k_O p_{NO} b_N^{T-2} \\
        &+ k_V (p_{VV}+1) b_V^{T-2} + k_O p_{VO} b_V^{T-2} + k_O (p_{OO}+1) b_O^{T-2} \\
        & - \Tilde{c}(p_{NO} p_{OO} p_{OIC} + p_{NV} p_{VO} p_{OIC} + p_{NV} p_{VV} p_{VIC} \\ 
        &+ p_{NN} p_{NO} p_{OIC} + p_{NN} p_{NV} p_{VIC} + p_{NN} p_{NN} p_{NIC}))b_N^{T-2}\\
        & -\Tilde{c}(p_{VO}p_{OO}p_{OIC} + p_{VV} p_{VO} p_{OIC} + p_{VV} p_{VV} p_{VIC}))b_V^{T-2} - \Tilde{c}p_{OO}p_{OO} p_{OIC}b_O^{T-2}\\
        &= 1 + (k_O p_{NO} + k_V p_{NV} + k_N (p_{NN}+1) - \Tilde{c}(p_{NO} p_{OO} p_{OIC} + p_{NV} p_{VO} p_{OIC} \\
        &+ p_{NV} p_{VV} p_{VIC} + p_{NN} p_{NO} p_{OIC} + p_{NN} p_{NV} p_{VIC} + p_{NN} p_{NN} p_{NIC}))b_N^{T-2}\\
        &+ (k_O p_{VO} + k_V (p_{VV}+1) -\Tilde{c}(p_{VO}p_{OO}p_{OIC} + p_{VV} p_{VO} p_{OIC} + p_{VV} p_{VV} p_{VIC}))b_V^{T-2} \\
        &+ (k_O (p_{OO}+1) - \Tilde{c}p_{OO}p_{OO} p_{OIC})b_O^{T-2}.
    \end{align*}

    Thus, for $t=T-2$, Equations \eqref{eq:rec} and \eqref{eq:nonrec} are equivalent.

    Suppose that Equations \eqref{eq:rec} and \eqref{eq:nonrec} are equivalent for $t=T-j, \ldots, T-2$. We now show they are also equivalent for $t=T-(j+1)$.

    Recall that, for all $j$, for any $\theta^{j+1} \in \Theta^{j+1}_{SS'}$, there is a $\bar{\theta}^j \in \Theta^j_{\theta^{j+1}_2 S'}$ such that $\theta^{j+1} = (\theta^{j+1}_1,\bar{\theta}^j)$. Thus,
    \begin{align*}
        &V^{T-(j+1)}(b^{T-(j+1)}, \sigma_{j+1}^{T-(j+1)}) \\
        &= 1 - P[(dr)]^{T-j} d - P[(cr)]^{T-j}c + P[(nr)]^{T-j} V^{T-j}(b^{T-j},\sigma_{j}^{T-j}) \\
        &= 1 - P[(dr)]^{T-j} d - P[(cr)]^{T-j} c \\
        &+ P[(nr)]^{T-j} (1 + \sum_{S \in \{N,V,O\}} \sum_{S \in \{N,V,O\}} \sum_{\substack{S' \in \{N,V,O\} \\S' \text{ is reachable from }S}} (k_{SS'}^j - \Tilde{c}_{SS'IC}^{j+1}) b_S^{T-j}) \\
        &= 1 + k_N b_N^{T-(j+1)} + k_V b_V^{T-(j+1)} + k_O b_O^{T-(j+1)} + \sum_{S' \in \{N,V,O\}} (k_{NS'}^j - \Tilde{c}_{NS'IC}^{j+1}) (p_{NN} b_N^{T-(j+1)}) \\
        & + \sum_{S \in \{V,O\}} (k_{VS'}^j - \Tilde{c}_{VS'IC}^{j+1}) (p_{NV} b_N^{T-(j+1)} + p_{VV} b_V^{T-(j+1)}) \\
        &+ (k_{OO}^j -\Tilde{c}_{OOIC}^{j+1}) (p_{NO} b_N^{T-(j+1)} + p_{VO} b_V^{T-(j+1)} + p_{OO} b_O^{T-(j+1)}) \\
        &= 1 + (k_N + \sum_{S' \in \{N,V,O\}} (k_{NS'}^j - \Tilde{c}_{NS'IC}^{j+1}) p_{NN} \\
        & +\sum_{S \in \{V,O\}} (k_{VS'}^j - \Tilde{c}_{VS'IC}^{j+1}) p_{NV} + (k_{OO}^j - \Tilde{c}_{OOIC}^{j+1})) p_{NO}) b_N^{T-(j+1)} \\
        & + (k_V + \sum_{S \in \{V,O\}} (k_{VS'}^j - \Tilde{c}_{VS'IC}^{j+1})) p_{VV} + (k_{OO}^j - \Tilde{c}_{OOIC}^{j+1}) p_{VO}) b_V^{T-(j+1)} \\
        &+ (k_O + (k_{OO}^j - \Tilde{c}_{OOIC}^{j+1})) p_{OO}) b_O^{T-(j+1)}.
    \end{align*}

    For any $S,S' \in \{N,V,O\}$, it is clear that
    \begin{align*}
        \Tilde{c}_{SS'IC}^{j+2} &= \Tilde{c} p_{S'IC} \sum_{\theta^{j+2} \in \Theta^{j+2}_{SS'}}\left(\prod_{l=1}^{j+1}p_{\theta_{l}^{j+1} \theta_{l+1}^{j+1}}\right)\\
        &= \Tilde{c} p_{S'IC} \sum_{\substack{\theta^{j+1}_2 \in \{N,V,O\}\\ \theta^{j+1}_2 \text{ is reachable from } S}} \sum_{\bar{\theta}^{j+1} \in \Theta^{j+1}_{\theta^{j+1}_2 S'}}\left(p_{S \theta^{j+1}_2} \prod_{l=2}^{j+1}p_{\theta_{l}^{j+1} \theta_{l+1}^{j+1}}\right) \\
        & \sum_{\substack{\theta^{j+1}_2 \in \{N,V,O\}\\ \theta^{j+1}_2 \text{ is reachable from } S}} p_{S \theta^{j+1}_2} \left(\Tilde{c} p_{S'IC}\sum_{\bar{\theta}^{j+1} \in \Theta^{j+1}_{\theta^{j+1}_2 S'}}\prod_{l=2}^{j+1}p_{\theta_{l}^{j+1} \theta_{l+1}^{j+1}}\right) \\
        & = \sum_{\substack{\theta^{j+1}_2 \in \{N,V,O\}\\ \theta^{j+1}_2 \text{ is reachable from } S}} p_{S \theta^{j+1}_2} \Tilde{c}_{\theta^{j+1}_2 S'IC}^{j+1}
    \end{align*}
    As in Theorem \ref{thrm:equiv}, we have that when $S \ne S'$, $k_{SS'}^{j+1} = \sum_{\substack{\theta^{j+1}_2 \in \{N,V,O\} \\ \theta^{j+1}_2 \text{ is reachable from } S}} k_{\bar{\theta}^j_2 S'} p_{S \theta^{j+1}_2}$, and when $S = S'$, $k_{SS}^{j+1} = k_S + k_{SS}^j p_{SS}$. Thus, 
    \small\begin{align*}
        &1 + (k_N + \sum_{S' \in \{N,V,O\}} (k_{NS'}^j - \Tilde{c}_{NS'IC}^{j+1}) p_{NN}  +\sum_{S \in \{V,O\}} (k_{VS'}^j - \Tilde{c}_{VS'IC}^{j+1}) p_{NV} + (k_{OO}^j - \Tilde{c}_{OOIC}^{j+1})) p_{NO}) b_N^{T-(j+1)} \\
        & + (k_V + \sum_{S \in \{V,O\}} (k_{VS'}^j - \Tilde{c}_{VS'IC}^{j+1})) p_{VV} + (k_{OO}^j - \Tilde{c}_{OOIC}^{j+1}) p_{VO}) b_V^{T-(j+1)} \\
        &+ (k_O + (k_{OO}^j - \Tilde{c}_{OOIC}^{j+1})) p_{OO}) b_O^{T-(j+1)} \\
        &= 1 + \sum_{S \in \{N,V,O\}} \sum_{\substack{S' \in \{N,V,O\} \\S' \text{ is reachable from }S}} (k_{SS'}^{j+1} - \Tilde{c}_{SS'IC}^{j+2}) b_S^{T-(j+1)}.
    \end{align*}\normalsize

    Thus, for $t = T-(j+1)$, Equations \eqref{eq:rec} and \eqref{eq:nonrecAlt} are equivalent.
\end{proof}

\noindent\textbf{Proof of Lemma \ref{lm:timeincAlt}:}
\begin{proof}
    Suppose $(P[(dr)]^t)$, $(P[(cr)]^t)$, $(P[(cr)|(i)]^t)$ and $(P[(nr)]^{t} P[(cr)|(i)]^{t})$ are non-decreasing sequences, $(P[(nr)]^t)$ is a non-increasing sequence and that $\alpha_d (P[(dr)]^{t+1} - P[(dr)]^t) + \alpha_c (P[(cr)]^{t+1} - P[(cr)]^t) \ge (P[(cr)|(i)]^{t+1} - P[(cr)|(i)]^t)$ for all $t$. Suppose for given $j \le t$, $V^t(b^t, \sigma_{j}^t) \ge V^t(b^t,\sigma_{j+1}^t)$. We prove $V^{t+1}(b^{t+1}, \sigma_{j}^{t+1}) \ge V^{t+1}(b^{t+1},\sigma_{j+1}^{t+1})$ by induction on $j$.

    First consider the case where $j=0$, i.e., there exists a $t$ such that $V^t(b^t, (i)) \ge V^t(b^t, (ni,i))$. We show $V^{t+1}(b^{t+1}, (i)) \ge V^{t+1}(b^{t+1}, (ni,i))$ by showing $V^t(b^t, (ni,i)) - V^{t+1}(b^{t+1}, (ni,i)) \ge V^t(b^t, (i)) - V^{t+1}(b^{t+1}, (i))$.
    \begin{align*}
        V^t(b^t, (ni,i)) - V^{t+1}(b^{t+1}, (ni,i)) & = 1 - P[(dr)]^{t+1}d - P[(cr)]^{t+1}c + P[(nr)]^{t+1} V^{t+1}(b^{t+1},(i)) \\
        & - \left(1 - P[(dr)]^{t+2}d - P[(cr)]^{t+2}c + P[(nr)]^{t+2} V^{t+2}(b^{t+2},(i))\right) \\
        & = (P[(dr)]^{t+2}-P[(dr)]^{t+1})d + (P[(cr)]^{t+2}-P[(cr)]^{t+1})c \\
        &+ P[(nr)]^{t+1} - P[(nr)]^{t+2} \\
        &+ \left(P[(nr)]^{t+2} P[(cr)|(i)]^{t+2} - P[(nr)]^{t+1} P[(cr)|(i)]^{t+1}\right)\Tilde{c} \\
        & \ge \left((P[(dr)]^{t+1}-P[(dr)]^{t})\alpha_d + (P[(cr)]^{t+1}-P[(cr)]^{t})\alpha_c \right)\Tilde{c} \\
        & \ge (P[(cr)|(i)]^{t+1} - P[(cr)|(i)]^t) \Tilde{c} \\
        & = (1 - \Tilde{c} P[(cr)|(i)]^t) - (1 - \Tilde{c} P[(cr)|(i)]^{t+1}) \\
        & = V^t(b^t, (i)) - V^{t+1}(b^{t+1}, (i)).
    \end{align*}

    Thus, 
    \begin{align*}
        V^t(b^t, (ni,i)) - V^{t+1}(b^{t+1}, (ni,i)) &\ge V^t(b^t, (i)) - V^{t+1}(b^{t+1}, (i)) \\
        V^{t+1}(b^{t+1}, (i)) - V^{t+1}(b^{t+1}, (ni,i)) &\ge V^{t}(b^{t}, (i)) - V^{t}(b^{t}, (ni,i)) \\
        V^{t+1}(b^{t+1}, (i)) - V^{t+1}(b^{t+1}, (ni,i)) &\ge 0 \\
        V^{t+1}(b^{t+1}, (i)) &\ge V^{t+1}(b^{t+1}, (ni,i)).
    \end{align*}
    
    Next, suppose $V^{t}(b^{t}, \sigma_{j}^{t}) \ge V^{t}(b^{t},\sigma_{j+1}^{t})$ implies $V^{t+1}(b^{t+1}, \sigma_{j}^{t+1}) \ge V^{t+1}(b^{t+1},\sigma_{j+1}^{t+1})$ for $0 \le j \le l < T-t$. We now show that if $V^{t}(b^{t}, \sigma_{j}^{t}) \ge V^{t}(b^{t},\sigma_{j+1}^{t})$, then $V^{t+1}(b^{t+1}, \sigma_{j}^{t+1}) \ge V^{t+1}(b^{t+1},\sigma_{j+1}^{t+1})$ for $j = l+1$. Suppose $V^{t}(b^{t}, \sigma_{l+1}^{t}) \ge V^{t}(b^{t},\sigma_{l+2}^{t})$, and note that 
    \begin{align*}
         &V^{t+1}(b^{t+1}, \sigma_{l+1}^{t+1}) \ge V^{t+1}(b^{t+1},\sigma_{j+1}^{t+1})\\
         &1 - P[(dr)]^{t+1} d - P[(cr)]^{t+1} c + P[(nr)]^{t+1}(V^{t+1}(b^{t+1},\sigma_l^{t+1})) \\
         &\ge 1 - P[(dr)]^{t+1} d - P[(cr)]^{t+1} c + P[(nr)]^{t+1}(V^{t+1}(b^{t+1},\sigma_{l+1}^{t+1})) \\
         &V^{t+1}(b^{t+1},\sigma_l^{t+1})) \ge V^{t+1}(b^{t+1},\sigma_{l+1}^{t+1}).
    \end{align*}

    Since $V^{t+1}(b^{t+1},\sigma_l^{t+1})) \ge V^{t+1}(b^{t+1},\sigma_{l+1}^{t+1})$, then $V^{t+2}(b^{t+2},\sigma_l^{t+2})) \ge V^{t+2}(b^{t+2},\sigma_{l+1}^{t+2})$ by our inductive hypothesis.
    \begin{align*}
        V^{t+1}(b^{t+1}, \sigma_{l+1}^{t+1})  & = 1 - P[(dr)]^{t+2}d - P[(cr)]^{t+2}c + P[(nr)]^{t+2} (V^{t+2}(b^{t+2},\sigma_l^{t+2})) \\
        & \ge 1 - P[(dr)]^{t+2}d - P[(cr)]^{t+2}c + P[(nr)]^{t+2} (V^{t+2}(b^{t+2},\sigma_{l+1}^{t+2})) \\
        & = V^{t+1}(b^{t+1}, \sigma_{l+2}^{t+1}).
    \end{align*}

    We thus have that if $V^{t}(b^{t}, \sigma_{j}^{t}) \ge V^{t}(b^{t},\sigma_{j+1}^{t})$ implies that $V^{t+1}(b^{t+1}, \sigma_{j}^{t+1}) \ge V^{t+1}(b^{t+1},\sigma_{j+1}^{t+1})$, then $V^{t}(b^{t}, \sigma_{j+1}^{t}) \ge V^{t}(b^{t},\sigma_{j+2}^{t})$ implies that $V^{t+1}(b^{t+1}, \sigma_{j+1}^{t+1}) \ge V^{t+1}(b^{t+1},\sigma_{j+2}^{t+1})$.
\end{proof}

\noindent\textbf{Proof of Theorem \ref{thrm:niiAlt}:}
\begin{proof}
    Suppose $(P[(dr)]^t)$, $(P[(cr)]^t)$, $(P[(cr)|(i)]^t)$ and $(P[(nr)]^{t} P[(cr)|(i)]^{t})$ are non-decreasing sequences, $(P[(nr)]^t)$ is a non-increasing sequence and that $\alpha_d (P[(dr)]^{t+1} - P[(dr)]^t) + \alpha_c (P[(cr)]^{t+1} - P[(cr)]^t) \ge (P[(cr)|(i)]^{t+1} - P[(cr)|(i)]^t)$ for all $t$. Additionally, suppose there exists a $t$ such that $V^t(b^t,(i)) \ge V^t(b^t, (ni,i))$. We show that $V^t(b^t, \sigma_{j}^t) \ge V^t(b^t,\sigma_{j+1}^t)$ for $1 \le j \le T-t$ be induction. First consider $j=1$. By Lemma \ref{lm:timeincAlt}, we also know that $V^{t+1}(b^{t+1}, (i)) \ge V^{t+1}(b^{t+1},(ni,i))$. We have
    \begin{align*}
        V^t(b^t,(ni,ni,i)) & = 1 - P[(dr)]^{t+1} d - P[(cr)]^{t+1} c + P[(nr)]^{t+1} (V^{t+1}(b^{t+1},(ni,i))) \\
        & \le 1 - P[(dr)]^{t+1} d - P[(cr)]^{t+1} c + P[(nr)]^{t+1}(V^{t+1}(b^{t+1},(i))) \\
        & = V^t(b^t, (ni,i))
    \end{align*}
    Thus, $V^t(b^t, (ni,i)) \ge V^t(b^t,(ni,ni,i))$.

    Now suppose $V^t(b^t,\sigma_{j}^t) \ge V^t(b^t,\sigma_{j+1}^t)$ for $1 \le j \le l < T-t$. We now show $V^t(b^t,\sigma_{j}^t) \ge V^t(b^t,\sigma_{j+1}^t)$ for $j=l+1$. Suppose $V^t(b^t,\sigma_{l}^t) \ge V^t(b^t,\sigma_{l+1}^t)$. We know that $V^{t+1}(b^{t+1},\sigma_{l}^{t+1}) \ge V^{t+1}(b^{t+1},\sigma_{l+1}^{t+1})$ by Lemma \ref{lm:timeincAlt}. We have that 
    \begin{align*}
        V^t(b^t,\sigma_{l+1}^t) &= 1 - P[(dr)]^{t+1} d - P[(cr)]^{t+1} c + P[(nr)]^{t+1} (V^{t+1}(b^{t+1}, \sigma_{l}^{t+1})) \\
        & \ge 1 - P[(dr)]^{t+1} d - P[(cr)]^{t+1} c + P[(nr)]^{t+1} (V^{t+1}(b^{t+1}, \sigma_{l+1}^{t+1})) \\
        & = V^t(b^t,\sigma_{l+2}^t)
    \end{align*}
    Thus, if $V^t(b^t,\sigma_{l}^t) \ge V^t(b^t,\sigma_{l+1}^t)$, then $V^t(b^t,\sigma_{l+1}^t) \ge V^t(b^t,\sigma_{l+2}^t)$.
\end{proof}

\section{Sensitivity of Optimal Inspection Time on Penalty Parameters}
\label{app:sens}
We now explore how varying $d$ and $c$ impacts $t^*$. For the entirety of this section, we suppose there exists a $t^*$ where the optimal decision is to inspect. We perform this analysis on the model presented in Section \ref{sec:model}, as the same logic can be applied to the model presented in Section \ref{sec:modelvar}. Additionally, we assume $b^1 = (1,0,0,0,0,0)$, as the same logic can be applied to any other starting $b^1$.

We first explore the impact of an increase in $d$ and $c$ on $t^*$, such that we would want to inspect earlier. Let $\delta$ be the increase in the manufacturing failure penalty and $\gamma$ be the increase in penalty for closing for non-mandatory maintenance. Let $\bar{k}_S = 1 - (d+\delta+1)p_{SD} - (c+\gamma+1)p_{SC}$. We want to explore what values $\delta$ and $\gamma$ need to take such that $V^1(b^1,\sigma^1_{t^*-1}) \ge V^1(b^1,\sigma^1_{t^*})$. The following is a simplified inequality on $\delta$ and $\gamma$.
\footnotesize{\begin{align*}
    &V^1(b^1,\sigma^1_{t^*-1}) - V^1(b^1,\sigma^1_{t^*}) \ge 0\\
    &\sum_{S' \in \{N,V,O\}} (k_{NS'}^{t^*-1} - k_{NS'}^{t^*}) \ge 0 \\
    &\sum_{S' \in \{N,V,O\}} \bar{k}_{S'}\left(\sum_{\theta^{t^*-1} \in \Theta_{NS'}^{t^*-1}} f_{\theta^{t^*-1}} - \sum_{\theta^{t^*} \in \Theta_{NS'}^{t^*}} f_{\theta^{t^*}}\right) \ge 0 \\
    &\sum_{S' \in \{N,V,O\}} (k_{S'} - p_{S'D} \delta -  p_{S'C}\gamma)\left(\sum_{\theta^{t^*-1} \in \Theta_{NS'}^{t^*-1}} f_{\theta^{t^*-1}} - \sum_{\theta^{t^*} \in \Theta_{NS'}^{t^*}} f_{\theta^{t^*}}\right) \ge 0 \\
    &\left(\sum_{S' \in \{N,V,O\}}\left(\sum_{\theta^{t^*-1} \in \Theta_{NS'}^{t^*-1}} f_{\theta^{t^*-1}} - \sum_{\theta^{t^*} \in \Theta_{NS'}^{t^*}} f_{\theta^{t^*}}\right)p_{S'D}\right) \delta \\
    &+ \left(\sum_{S' \in \{N,V,O\}}\left(\sum_{\theta^{t^*-1} \in \Theta_{NS'}^{t^*-1}} f_{\theta^{t^*-1}} - \sum_{\theta^{t^*} \in \Theta_{NS'}^{t^*}} f_{\theta^{t^*}}\right)p_{S'C}\right)\gamma \\
    & \le \sum_{S' \in \{N,V,O\}} k_{S'}\left(\sum_{\theta^{t^*-1} \in \Theta_{NS'}^{t^*-1}} f_{\theta^{t^*-1}} - \sum_{\theta^{t^*} \in \Theta_{NS'}^{t^*}} f_{\theta^{t^*}}\right)
\end{align*}}

\normalsize
Note that as $p_{S'D}$ and $p_{S'C}$ increase, $k_{S'}$ decreases. This indicates that, the more prone to experiencing unexpected disruptive events, the more sensitive our optimal inspection time is to a change in the penalties from unexpected disruptive events occurring.

We next consider what values $d$ and $c$ can take such that a pre-specified time period $t^*$ is the optimal inspection time. By Corollary \ref{cor:structure}, we know that $V^1(b^1,\sigma^1_{t^*-1}) \le V^1(b^1,\sigma^1_{t^*})$ and $V^1(b^1,\sigma^1_{t^*}) \ge V^1(b^1,\sigma^1_{t^*+1})$ is sufficient to imply that $V^1(b^1,\sigma^1_{t^*}) \ge V^1(b^1,\sigma^1_{t})$ for all $t$. Thus, we only need to compare the value functions for $t^*-1$, $t^*$ and $t^*+1$ to determine bounds on $d$ and $c$.

First, to ensure that we should wait to inspect no later than $t^*$, we devise conditions on $d$ and $c$ based on $V^1(b^1,\sigma^1_{t^*}) \ge V^1(b^1,\sigma^1_{t^*+1})$. 
\footnotesize{\begin{align*}
    &V^1(b^1,\sigma^1_{t^*}) - V^1(b^1,\sigma^1_{t^*+1}) \ge 0\\
    &\sum_{S' \in \{N,V,O\}} (k_{NS'}^{t^*} - k_{NS'}^{t^*+1}) \ge 0 \\
    &\sum_{S' \in \{N,V,O\}} k_{S'}\left(\sum_{\theta^{t^*} \in \Theta_{NS'}^{t^*}} f_{\theta^{t^*}} - \sum_{\theta^{t^*+1} \in \Theta_{NS'}^{t^*+1}} f_{\theta^{t^*+1}}\right) \ge 0 \\
    &\sum_{S' \in \{N,V,O\}} (1 - (d+1)p_{S'D} - (c+1)p_{S'C}) \left(\sum_{\theta^{t^*} \in \Theta_{NS'}^{t^*}} f_{\theta^{t^*}} - \sum_{\theta^{t^*+1} \in \Theta_{NS'}^{t^*+1}} f_{\theta^{t^*+1}}\right) \ge 0 \\
    &\left(\sum_{S' \in \{N,V,O\}}\left(\sum_{\theta^{t^*} \in \Theta_{NS'}^{t^*}} f_{\theta^{t^*}} - \sum_{\theta^{t^*+1} \in \Theta_{NS'}^{t^*+1}} f_{\theta^{t^*+1}}\right)p_{S'D}\right) d \\
    &+ \left(\sum_{S' \in \{N,V,O\}}\left(\sum_{\theta^{t^*} \in \Theta_{NS'}^{t^*}} f_{\theta^{t^*}} - \sum_{\theta^{t^*+1} \in \Theta_{NS'}^{t^*+1}} f_{\theta^{t^*+1}}\right)p_{S'C}\right) c \\
    &\le \sum_{S' \in \{N,V,O\}} (1-p_{S'D} - p_{S'C}) \left(\sum_{\theta^{t^*} \in \Theta_{NS'}^{t^*}} f_{\theta^{t^*}} - \sum_{\theta^{t^*+1} \in \Theta_{NS'}^{t^*+1}} f_{\theta^{t^*+1}}\right).
\end{align*}} \normalsize
If $d$ and $c$ satisfy the above inequality, then it is better to inspect in time period $t^*$ than it would be delay inspection.

Likewise, we can derive conditions on $d$ and $c$ such that we should inspect no earlier than $t^*$ based on $V^1(b^1,\sigma^1_{t^*}) \ge V^1(b^1,\sigma^1_{t^*-1})$. Using the same derivation as above, we arrive at:
\footnotesize{\begin{align*}
    &\left(\sum_{S' \in \{N,V,O\}}\left(\sum_{\theta^{t^*} \in \Theta_{NS'}^{t^*}} f_{\theta^{t^*}} - \sum_{\theta^{t^*-1} \in \Theta_{NS'}^{t^*-1}} f_{\theta^{t^*-1}}\right)p_{S'D}\right) d \\
    &+ \left(\sum_{S' \in \{N,V,O\}}\left(\sum_{\theta^{t^*} \in \Theta_{NS'}^{t^*}} f_{\theta^{t^*}} - \sum_{\theta^{t^*+-} \in \Theta_{NS'}^{t^*-1}} f_{\theta^{t^*-1}}\right)p_{S'C}\right) c \\
    &\le \sum_{S' \in \{N,V,O\}} (1-p_{S'D} - p_{S'C}) \left(\sum_{\theta^{t^*} \in \Theta_{NS'}^{t^*}} f_{\theta^{t^*}} - \sum_{\theta^{t^*-1} \in \Theta_{NS'}^{t^*-1}} f_{\theta^{t^*-1}}\right).
\end{align*}}
\normalsize
If $d$ and $c$ satisfy the above inequality, then it is better to inspect in time period $t^*$ than it would be inspect earlier. Thus, if $d$ and $c$ lie in the polytope defined by the two inequalities, then it is optimal to inspect in time period $t^*$.

\section{Additional Computational Study}
\label{app:acs}
\begin{figure}[h!]
\centering
\begin{subfigure}{.48\textwidth}
  \centering
  \includegraphics[width=.95\linewidth]{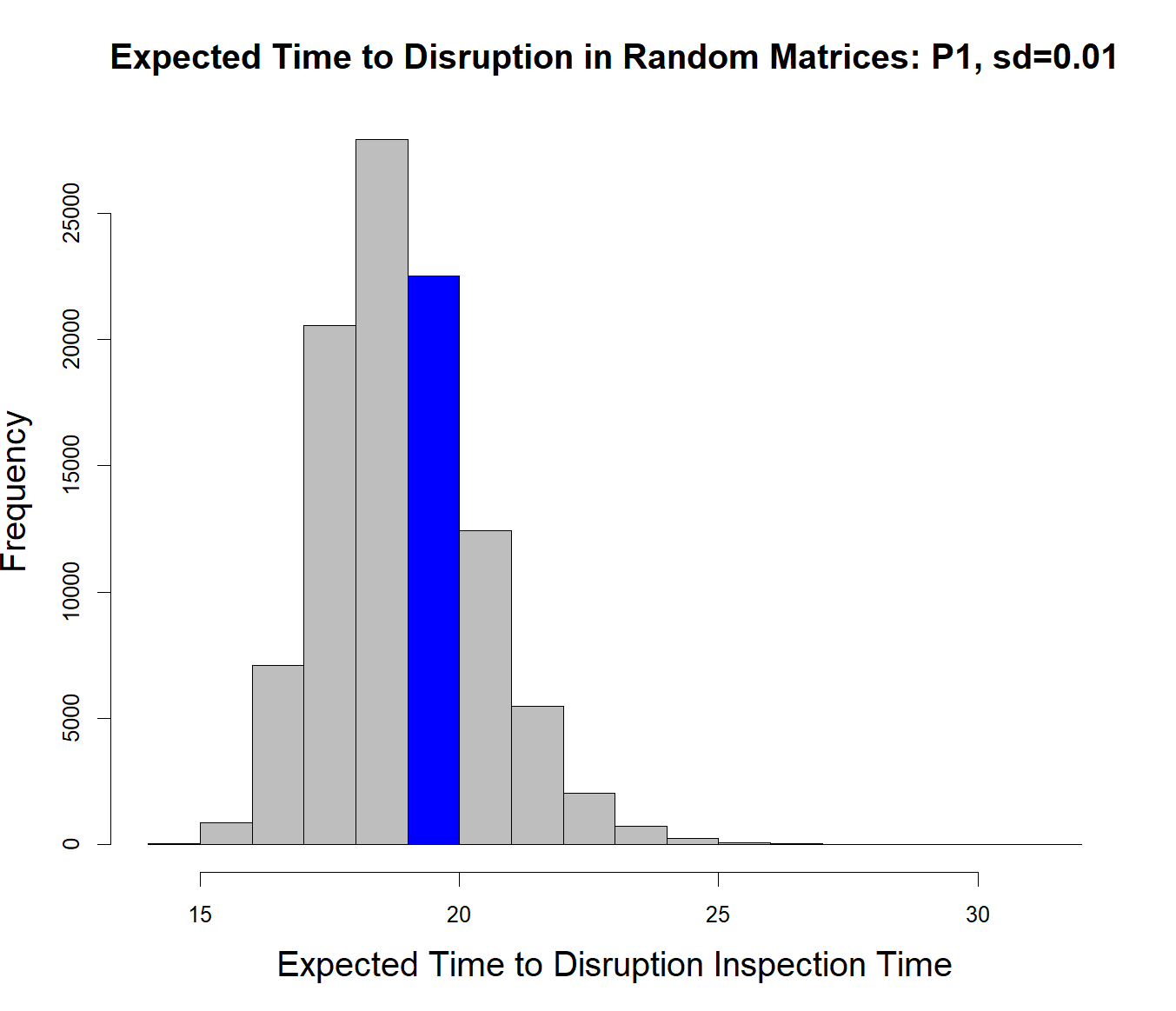}
  \caption{$d=18$, $s=0.01$}
  \label{fig:p1_tETD_d18_sd1}
\end{subfigure}%
\begin{subfigure}{.48\textwidth}
  \centering
  \includegraphics[width=.95\linewidth]{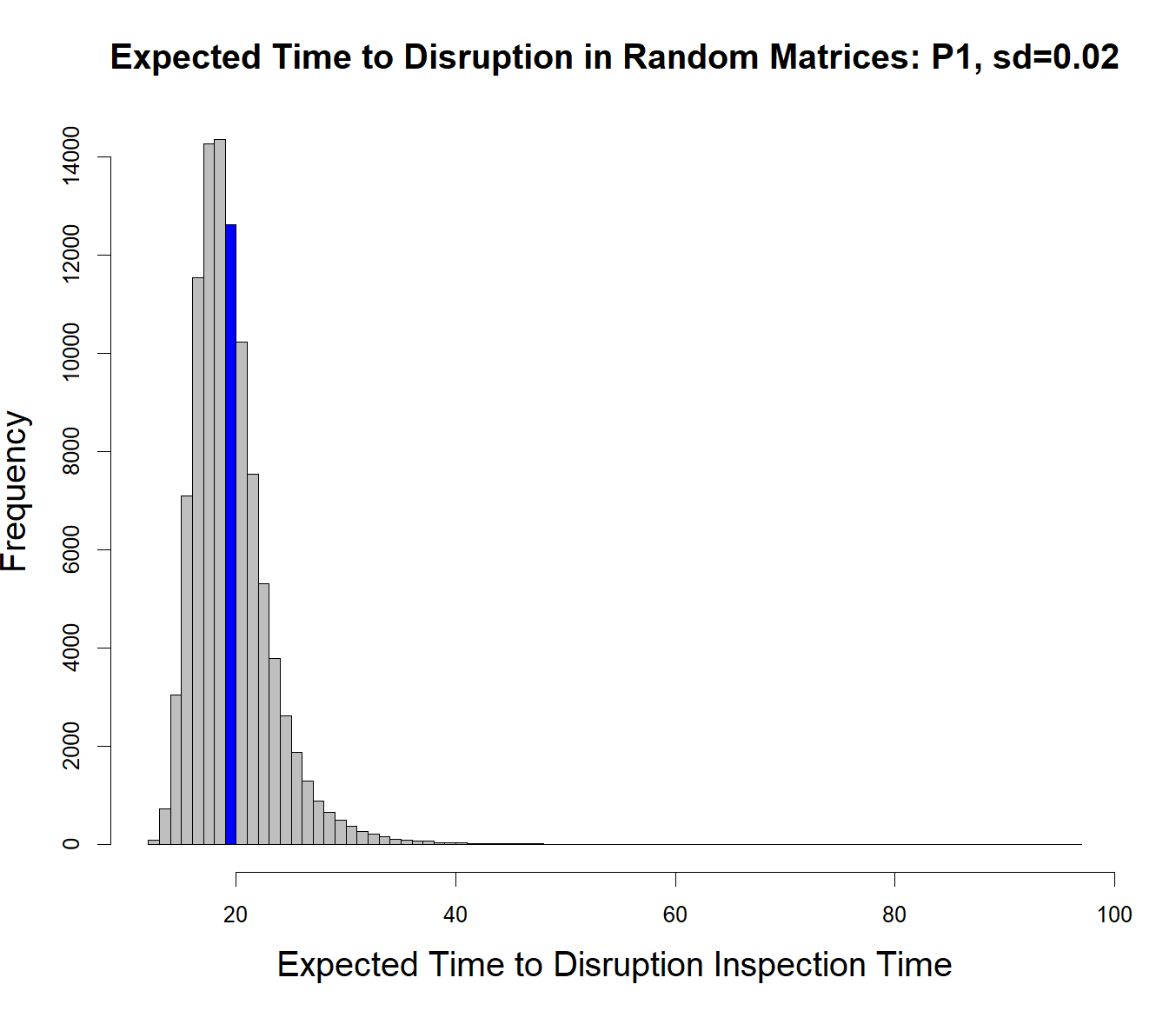}
  \caption{$d=18$, $s=0.02$}
  \label{fig:p1_tETD_d18_sd2}
\end{subfigure}
\caption{Histograms of ETD inspection rule recommendations for randomly generated matrices}
\label{fig:p1_tETD_d18}
\end{figure}

%\begin{figure}[h!]
%\centering
%\begin{subfigure}{.48\textwidth}
%  \centering
%  \includegraphics[width=.9\linewidth]{P1NN_dist.png}
%  \caption{$P^1(ni)$}
%  \label{fig:p1nn}
%\end{subfigure}%
%\begin{subfigure}{.48\textwidth}
%  \centering
%  \includegraphics[width=.9\linewidth]{P2NN_dist.png}
%  \caption{$P^2(ni)$}
%  \label{fig:p2nn}
%\end{subfigure}%
%
%\begin{subfigure}{.48\textwidth}
%  \centering
%  \includegraphics[width=.9\linewidth]{P3NN_dist.png}
%  \caption{$P^3(ni)$}
%  \label{fig:p3nn}
%\end{subfigure}
%\caption{Probability density functions with varying standard deviations and means from each transition matrix transitioning from state $N$ to state $N$}
%\label{fig:pnn}
%\end{figure}

\begin{table}[h!]
\centering
\resizebox{\textwidth}{!}{\begin{tabular}{|c|c|c|c|c|c|c|c|}
\hline
Rule          & $t_E$ & $t_V$ & $t_{VC}$ & 24 & 60 & 120 & NoIns \\ \hline\hline
\% Caught, $d=10$ & 46.530\% & +5.601\% & +5.601\% & -37.880\% & -46.513\% & -46.530\% & -  \\ \hline
Value (no IC), $d=10$ & 4.6743 & +0.128\% & - & -12.062\% & -16.856\% & -16.871\% & -16.871\% \\ \hline
Value (IC), $d=10$ & 4.4856 & - & -0.145\% & -9.352\% & -13.258\% & -13.374\% & - \\ \hline\hline
\% Caught, $d=14$ & 46.070\% & +30.244\% & +30.244\% & -37.488\% & -46.047\% & -46.070\% & -  \\ \hline
Value (no IC), $d=14$ & 3.0582 & +16.470\% & - & -51.592\% & -66.555\% & -66.578\% & -66.578\% \\ \hline
Value (IC), $d=14$ & 2.8713 & - & +16.386\% & -50.009\% & -64.382\% & -64.403\% & - \\ \hline\hline
\% Caught, $d=18$ & 46.601\% & +36.149\% & +36.149\% & -37.787\% & -46.579\% & -46.601\% & -  \\ \hline
Value (no IC), $d=18$ & 1.6318 & +84.894\% & - & --160.526\% & -205.050\% & -205.142\% & -205.142\% \\ \hline
Value (IC), $d=18$ & 1.4439 & - & +94.376\% & -171.535\% & -218.727\% & -218.824\% & - \\ \hline
\end{tabular}}
\caption{Percent of Unexpected Disruptive Events Caught and Average Value Accumulated with Each Inspection Rule for $\Tilde{P}^2(ni)$ with $s = 0.01$}
\label{tbl:value_p2_s1}
\end{table}

%\begin{table}[h!]
%\centering \small
%\begin{tabular}{|c|c|c|c|}
%\hline
%$d$ & 10 & 14 & 18 \\\hline
%Average End Time & 13.3842 & 13.3194 & 13.4198 \\ \hline
%Std. Dev. End Time & 7.6792 & 7.6661 & 7.6970 \\ \hline
%Median End Time & 12 & 12 & 12 \\ \hline
%Min End Time & 2 & 2 & 2 \\ \hline
%Max End Time & 71 & 85 & 81 \\ \hline
%\end{tabular}
%\caption{Descriptive Statistics of Time of Unexpected Disruptive Event for %$\Tilde{P}^2(ni)$, $s=0.01$}
%\label{tbl:endtime_p2_s1}
%\end{table}\normalsize

\begin{table}[h!]
\centering
\resizebox{\textwidth}{!}{\begin{tabular}{|c|c|c|c|c|c|c|c|}
\hline
Rule          & $t_E$ & $t_V$ & $t_{VC}$ & 24 & 60 & 120 & NoIns \\ \hline\hline
\% Caught, $d=10$ & 45.903\% & +5.542\% & +5.542\% & -36.848\% & -45.872\% & -45.903\% & -  \\ \hline
Value (no IC), $d=10$ & 4.5517 & +0.204\% & - & -10.058\% & -14.838\% & -14.856\% & -14.856\% \\ \hline
Value (IC), $d=10$ & 4.3683 & - & -0.094\% & -7.330\% & -11.265\% & -11.281\% & - \\ \hline\hline
\% Caught, $d=14$ & 46.091\% & +29.706\% & +29.706\% & -36.881\% & -46.063\% & -46.091\% & -  \\ \hline
Value (no IC), $d=14$ & 3.0497 & +13.605\% & - & -46.191\% & -61.091\% & -61.157\% & -61.157\% \\ \hline
Value (IC), $d=14$ & 2.8672 & - & +13.382\% & -44.357\% & -58.622\% & -58.684\% & - \\ \hline\hline
\% Caught, $d=18$ & 46.103\% & +35.745\% & +35.745\% & -37.068\% & -46.065\% & -46.103\% & -  \\ \hline
Value (no IC), $d=18$ & 1.5778 & +80.980\% & - & -156.832\% & -200.615\% & -200.780\% & -200.780\% \\ \hline
Value (IC), $d=18$ & 1.3936 & - & +90.342\% & -167.580\% & -213.928\% & -214.100\% & - \\ \hline
\end{tabular}}
\caption{Percent of Unexpected Disruptive Events Caught and Average Value Accumulated with Each Inspection Rule for $\Tilde{P}^2(ni)$ with $s = 0.02$}
\label{tbl:value_p2_s2}
\end{table}

%\begin{table}[h!]
%\centering \small 
%\begin{tabular}{|c|c|c|c|}
%\hline
%$d$ & 10 & 14 & 18 \\\hline
%Average End Time & 13.3362 & 13.3004 & 13.3820 \\ \hline
%Std. Dev. End Time & 7.8783 & 7.9553 & 7.9092 \\ \hline
%Median End Time & 12 & 12 & 12 \\ \hline
%Min End Time & 2 & 2 & 2 \\ \hline
%Max End Time & 80 & 80 & 90 \\ \hline
%\end{tabular}
%\caption{Descriptive Statistics of Time of Unexpected Disruptive Event for %$\Tilde{P}^2(ni)$, $s=0.02$}
%\label{tbl:endtime_p2_s2}
%\end{table}\normalsize

\begin{table}[h!]
\centering
\resizebox{\textwidth}{!}{\begin{tabular}{|c|c|c|c|c|c|c|c|}
\hline
Rule          & $t_E$ & $t_V$ & $t_{VC}$ & 24 & 60 & 120 & NoIns \\ \hline\hline
\% Caught, $d=18$ & 43.217\% & -41.990\% & -43.217\% & +2.406\% & -38.030\% & -43.059\% & -  \\ \hline
Value (no IC), $d=18$ & 10.9907 & +7.139\% & - & -0.812\% & +6.541\% & +7.414\% & +7.466\% \\ \hline
Value (IC), $d=18$ & 10.8724 & - & +8.636\% & -0.871\% & +7.561\% & +8.579\% & - \\ \hline\hline
\% Caught, $d=22$ & 43.432\% & +2.478\% & +2.478\% & +2.478\% & -38.283\% & -43.286\% & -  \\ \hline
Value (no IC), $d=22$ & 9.4174 & -0.096\% & - & -0.096\% & -4.358\% & -5.149\% & -5.156\% \\ \hline
Value (IC), $d=22$ & 9.2988 & - & -0.152\% & -0.251\% & -3.309\% & -3.942\% & - \\ \hline\hline
\% Caught, $d=26$ & 43.114\% & +22.641\% & +22.641\% & +2.487\% & -37.955\% & -42.974\% & -  \\ \hline
Value (no IC), $d=26$ & 7.7032 & +4.618\% & - & +0.946\% & -18.996\% & -21.822\% & -21.912\% \\ \hline
Value (IC), $d=26$ & 7.5860 & - & +4.163\% & +0.895\% & -17.944\% & -20.620\% & - \\ \hline\hline
\% Caught, $d=30$ & 43.186\% & +32.031\% & +32.031\% & +2.416\% & -37.940\% & -43.037\% & - \\ \hline
Value (no IC), $d=30$ & 6.0871 & +17.488\% & - & +2.249\% & -41.802\% & -47.436\% & -47.532\% \\ \hline
Value (IC), $d=30$ & 5.9699 & - & +17.136\% & +2.206\% & -40.927\% & -46.413\% & - \\ \hline\hline
\% Caught, $d=34$ & 43.371\% & +38.178\% & +38.178\% & +2.466\% & -38.086\% & -43.219\% & - \\ \hline
Value (no IC), $d=34$ & 4.6578 & +41.288\% & - & +4.386\% & -78.301\% & -89.282\% & -89.542\% \\ \hline
Value (IC), $d=34$ & 4.5378 & - & +41.558\% & +4.392\% & -78.071\% & -89.006\% & - \\ \hline
\end{tabular}}
\caption{Percent of Unexpected Disruptive Events Caught and Average Value Accumulated with Each Inspection Rule for $\Tilde{P}^3(ni)$ with $s = 0.01$}
\label{tbl:value_p3_s1}
\end{table}

%\begin{table}[h!]
%\centering \small 
%\begin{tabular}{|c|c|c|c|c|c|}
%\hline
%$d$ & 18 & 22 & 26 & 30 & 34 \\\hline
%Average End Time & 27.0599 & 27.0225 & 26.9900 & 27.0438 & 27.0988 \\ \hline
%Std. Dev. End Time & 17.6907 & 17.5354 & 17.5645 & 17.6867 & 17.7312 \\ %\hline
%Median End Time & 23 & 23 & 23 & 23 & 23 \\ \hline
%Min End Time & 3 & 3 & 3 & 3 & 3 \\ \hline
%Max End Time & 239 & 235 & 205 & 242 & 253 \\ \hline
%\end{tabular}
%\caption{Descriptive Statistics of Time of Unexpected Disruptive Event for %$\Tilde{P}^3(ni)$, $s=0.01$}
%\label{tbl:endtime_p3_s1}
%\end{table}\normalsize

\begin{table}[h!]
\centering
\resizebox{\textwidth}{!}{\begin{tabular}{|c|c|c|c|c|c|c|c|}
\hline
Rule          & $t_E$ & $t_V$ & $t_{VC}$ & 24 & 60 & 120 & NoIns \\ \hline\hline
\% Caught, $d=18$ & 45.272\% & -42.632\% & -45.272\% & +2.332\% & -37.909\% & -44.490\% & -  \\ \hline
Value (no IC), $d=18$ & 11.7223 & +18.104\% & - & -1.059\% & +14.699\% & +20.384\% & +23.580\% \\ \hline
Value (IC), $d=18$ & 11.6058 & - & +24.820\% & -1.107\% & +15.702\% & +21.582\% & - \\ \hline\hline
\% Caught, $d=22$ & 45.140\% & +2.359\% & +2.359\% & +2.359\% & -37.824\% & -44.270\% & -  \\ \hline
Value (no IC), $d=22$ & 10.2258 & -0.628\% & - & -0.628\% & +6.477\% & +11.255\% & +15.580\% \\ \hline
Value (IC), $d=22$ & 10.1099 & - & -0.684\% & -0.684\% & +7.324\% & +12.519\% & - \\ \hline\hline
\% Caught, $d=26$ & 44.755\% & +21.713\% & +21.713\% & +2.394\% & -37.573\% & -43.939\% & -  \\ \hline
Value (no IC), $d=26$ & 8.6925 & -1.137\% & - & +0.2222\% & -5.215\% & -0.978\% & +2.625\% \\ \hline
Value (IC), $d=26$ & 8.5782 & - & -1.587\% & +0.162\% & -3.566\% & +0.329\% & - \\ \hline\hline
\% Caught, $d=30$ & 45.051\% & +30.695\% & +30.695\% & +2.415\% & -37.748\% & -44.285\% & - \\ \hline
Value (no IC), $d=30$ & 7.2452 & +4.766\% & - & +1.159\% & -19.512\% & -18.117\% & -14.008\% \\ \hline
Value (IC), $d=30$ & 7.1269 & - & +4.296\% & +1.107\% & -18.429\% & -16.774\% & - \\ \hline\hline
\% Caught, $d=34$ & 45.031\% & +36.747\% & +36.747\% & +2.396\% & -37.813\% & -44.270\% & - \\ \hline
Value (no IC), $d=34$ & 5.8057 & +19.708\% & - & +2.689\% & -43.952\% & -45.051\% & -40.493\% \\ \hline
Value (IC), $d=34$ & 5.6868 & - & +19.526\% & +2.641\% & -42.701\% & -43.919\% & - \\ \hline
\end{tabular}}
\caption{Percent of Unexpected Disruptive Events Caught and Average Value Accumulated with Each Inspection Rule for $\Tilde{P}^3(ni)$ with $s = 0.02$}
\label{tbl:value_p3_s2}
\end{table}

%\begin{table}[h!]
%\centering \small
%\begin{tabular}{|c|c|c|c|c|c|}
%\hline
%$d$ & 18 & 22 & 26 & 30 & 34 \\\hline
%Average End Time & 29.3162 & 29.3753 & 29.1401 & 29.1671 & 29.1348 \\ \hline
%Std. Dev. End Time & 24.6145 & 25.8545 & 24.3033 & 24.3668 & 24.3818 \\ %\hline
%Median End Time & 24 & 24 & 23 & 23 & 23 \\ \hline
%Min End Time & 3 & 3 & 3 & 3 & 3 \\ \hline
%Max End Time & 979 & 1048 & 1013 & 1059 & 1169 \\ \hline
%\end{tabular}
%\caption{Descriptive Statistics of Time of Unexpected Disruptive Event for $\Tilde{P}^3(ni)$, $s=0.02$}
%\label{tbl:endtime_p3_s2}
%\end{table}\normalsize
\end{document}